\documentclass[a4paper,twoside]{amsart}
\usepackage{xypic,amscd,amssymb,latexsym}
\makeatletter
\def\@setcopyright{}
\makeatother
\renewcommand{\-}{\hbox{-}}
\newcommand{\X}{\raisebox{0.43ex}{$\chi$}}
\newcommand{\eins}{\mbox{\rm 1\hspace{-0.24em}l}}
\renewcommand{\mod}{\operatorname{mod}\nolimits}
\newcommand{\Sub}{\operatorname{Sub}\nolimits}
\newcommand{\add}{\operatorname{add}\nolimits}
\newcommand{\Sets}{\operatorname{Sets}\nolimits}
\newcommand{\Hom}{\operatorname{Hom}\nolimits}
\newcommand{\End}{\operatorname{End}\nolimits}
\newcommand{\Aut}{\operatorname{Aut}\nolimits}
\newcommand{\Irr}{\operatorname{Irr}\nolimits}
\newcommand{\Inn}{\operatorname{Inn}\nolimits}
\newcommand{\Out}{\operatorname{Out}\nolimits}
\renewcommand{\Im}{\operatorname{Im}\nolimits}
\newcommand{\Ker}{\operatorname{Ker}\nolimits}
\newcommand{\Coker}{\operatorname{Coker}\nolimits}
\newcommand{\rad}{\operatorname{rad}\nolimits}
\newcommand{\soc}{\operatorname{soc}\nolimits}
\newcommand{\ann}{\operatorname{ann}\nolimits}
\newcommand{\Tor}{\operatorname{Tor}\nolimits}
\newcommand{\Tr}{\operatorname{Tr}\nolimits}
\newcommand{\GL}{\operatorname{GL}\nolimits}
\newcommand{\GF}{\operatorname{GF}\nolimits}
\newcommand{\Mat}{\operatorname{Mat}\nolimits}
\newcommand{\PSL}{\operatorname{PSL}\nolimits}
\newcommand{\SL}{\operatorname{SL}\nolimits}
\newcommand{\Fi}{\operatorname{Fi}\nolimits}
\newcommand{\Co}{\operatorname{Co}\nolimits}
\newcommand{\J}{\operatorname{J}\nolimits}
\newcommand{\kar}{\operatorname{char}\nolimits}
\newcommand{\Ext}{\operatorname{Ext}\nolimits}
\newcommand{\op}{{\operatorname{op}}}
\newcommand{\tr}{{\operatorname{tr}}}
\newcommand{\Ab}{\operatorname{Ab}\nolimits}
\newcommand{\id}{\operatorname{id}\nolimits}
\newcommand{\comp}{\operatorname{\scriptstyle\circ}}
\newcommand{\E}{\operatorname{\mathcal E}\nolimits}
\newcommand{\I}{\operatorname{\mathcal I}\nolimits}
\newcommand{\M}{\operatorname{\mathcal M}\nolimits}
\renewcommand{\P}{\operatorname{\mathcal P}\nolimits}
\newcommand{\R}{\operatorname{\mathfrak{R}}\nolimits}
\newcommand{\Y}{\operatorname{\mathcal Y}\nolimits}
\newcommand{\w}{\operatorname{\mathcal W}\nolimits}
\newcommand{\G}{\Gamma}
\renewcommand{\L}{\Lambda}
\newcommand{\la}{\lambda}
\renewcommand{\r}{\operatorname{\underline{r}}\nolimits}
\renewcommand{\a}{\underline{\mathrm a}}
\renewcommand{\b}{\underline{\mathrm b}}
\newcommand{\m}{\underline{\mathrm m}}
\newcommand{\dwnsimeq}{\Big\downarrow\lefteqn{\wr}}
\newcommand{\upsimeq}{\Big\uparrow\lefteqn{\wr}}
\newcommand{\subneq}{\hspace*{2pt}\mathrel
{{\raisebox{-.3ex}{\hbox{$\scriptscriptstyle\not$}}}\!\!\!\subseteq}}
\newcommand{\semi}{\mathbin{\vcenter{\hbox{$\scriptscriptstyle|$}}\;\!\!\!\times}}
\newcommand{\longequals}{\Relbar\joinrel=}
\newcommand{\A}{\alpha}
\newcommand{\B}{\beta}
\newcommand{\C}{\gamma}
\newcommand{\D}{\delta}
\newcommand{\T}{\tau}
\newcommand{\Z}{\xi}

\newtheorem{lemma}{Lemma}[section]
\newtheorem{proposition}[lemma]{Proposition}

\newtheorem{theorem}[lemma]{Theorem}
\newtheorem{remark}[lemma]{Remark}

\newtheorem{cclass}[lemma]{}
\newtheorem{ctab}[lemma]{}

\renewcommand{\arraystretch}{1.5}

\begin{document}

\thispagestyle{empty}
\title{ CONSTRUCTION OF FISCHER'S SPORADIC GROUP $\Fi_{24}'$ INSIDE $\GL_{8671}(13)$}

\author{Hyun Kyu Kim}
\address{Department of Mathematics, Yale University, New Haven,
  CT. 06511, USA}

\author{Gerhard O. Michler}
\address{Department of Mathematics, Cornell University, Ithaca,
  N.Y. 14853, USA}

\begin{abstract}
In this article we construct an irreducible simple subgroup
$\mathfrak G = \langle \mathfrak q, \mathfrak y, \mathfrak t,
\mathfrak w\rangle$ of $\GL_{8671}(13)$ from an irreducible
subgroup $T$ of $\GL_{11}(2)$ isomorphic to Mathieu's simple group
$\M_{24}$ by means of Algorithm 2.5 of \cite{michler1}. We also
use the first author's similar construction of Fischer's sporadic
simple group $G_1 = \Fi_{23}$ described in \cite{kim1}. He starts
from an irreducible subgroup $T_1$ of $\GL_{11}(2)$ contained in
$T$ which is isomorphic to $\M_{23}$. In \cite{hall} J. Hall and
L. S. Soicher published a nice presentation of Fischer's original
$3$-transposition group $\Fi_{24}$ \cite{fischer1}. It is used
here to show that $\mathfrak G$ is isomorphic to the simple
commutator subgroup $\Fi_{24}'$  of $\Fi_{24}$. We also determine
a faithful permutation representation of $\mathfrak G$ of degree
$306936$ with stabilizer $\mathfrak G_1 = \langle \mathfrak q,
\mathfrak y, \mathfrak w\rangle \cong \Fi_{23}$. It enabled MAGMA
to calculate the character table of $\mathfrak G$ automatically.

Furthermore, we prove that $\mathfrak G$ has two conjugacy classes
of involutions $\mathfrak z$ and $\mathfrak u$ such that
$C_{\mathfrak G}(\mathfrak u) = \langle \mathfrak q, \mathfrak y,
\mathfrak t\rangle \cong 2Aut(\Fi_{22})$. Moreover, we determine a
presentation of $\mathfrak H = C_{\mathfrak G}(\mathfrak z)$ and a
faithful permutation representation of degree $258048$ for which
we document a stabilizer.
\end{abstract}

\maketitle

\section{Introduction}

In $1971$ B. Fischer \cite{fischer1} discovered $3$ sporadic
simple groups by characterizing all finite groups $G$ which can be
generated by a conjugacy class $D =z^G$ of $3$-transpositions.
This means that the product of $2$ elements of $D$ in $G$ has order $1$,
$2$ or $3$. The largest of these $3$ sporadic groups turned out to be
the commutator subgroup $\Fi_{24}'$ of the $3$-transposition group
$\Fi_{24}$. In \cite{fischer} Fischer constructed for each of the
$3$ groups $\Fi_k$ a graph $\mathcal G_k$ and showed that $\Fi_k$
is isomorphic to the automorphism group $\rm Aut(\mathcal G_k)$,
$k \in \{22, 23, 24\}$. Using the structure of $\mathcal G_{24}$
J. Hall and L. Soicher determined a nice presentation of
$\Fi_{24}$, see \cite{hall} and \cite{praeger}, p. 124.

The results of this article are part of our joint research project
{\em Simultaneous construction of the sporadic simple groups of
Conway, Fischer and Janko}. Its goal is to provide uniform
existence proofs for the sporadic sporadic simple groups discovered by Conway, Fischer and Janko by
means of Algorithm 2.5 of \cite{michler1} constructing finite
simple groups from irreducible subgroups $T$ of $\GL_n(2)$. In
\cite{kim} we constructed Conway's and Fischer's sporadic groups $\Co_2$ and $\Fi_{22}$ simultaneously
from the irreducible subgroup $\M_{22}$ in $\GL_{10}(2)$. In \cite{kim1} the
first author applied the same methods to the
irreducible subgroup $\M_{23}$ of $\GL_{11}(2)$ and realized
$\Fi_{23}$ as an irreducible subgroup of $\GL_{782}(17)$.

In Lemma \ref{l. M24-extensions} of Section 2 we construct a
presentation of the unique non split extension $E$ of $\M_{24}$ by
the natural $GF(2)$-vector space $V$ of dimension $11$ by means of
Holt's Algorithm \cite{holt5} implemented in MAGMA. By Fischer's
work \cite{fischer1} it is known that this extension group has a
Sylow $2$-subgroup $S$ which is isomorphic to the ones of his
sporadic simple group $\Fi_{24}'$. Lemma \ref{l. M24-extensions}
also states that $E$ has a unique conjugacy class of $2$-central
involutions $z$. Furthermore, it has a subgroup $N_1$ of index
$24$ and a non $2$-central involution $u$ such that $N_1/V \cong
\M_{23}$ and $E = \langle N_1, C_E(u)\rangle$.

In our first attempt we applied Algorithm 2.5 of \cite{michler1}
to $E$ and constructed a finitely presented group $H$ containing a
Sylow $2$-subgroup $S_1$ having a maximal elementary abelian
normal subgroup $A$ such that $N_H(A) \cong D = C_E(z)$ and $|H :
N_H(A)|$ is odd. Furthermore, we calculated the character tables
of $E$, $H$ and $D$ and applied Algorithm 7.4.8 of \cite{michler}
to show that the free product $Q = H*_DE$ with amalgamated
subgroup $D$ has $939,080,024,064$ irreducible representations of
minimal dimension $8671$ over $GF(13)$. In view of our time
constraints we decided not to start Thompson's amalgamation
process described in Theorem 7.2.2 of \cite{michler} in order to
find an irreducible representation of $Q$ of degree $8671$ which
has a Sylow $2$-subgroup isomorphic to $S_1 \cong S$. Instead we
now construct such a matrix representation $\mathfrak G$ of $Q$
using the first author's work \cite{kim1}.

There he applied Algorithm 2.5 of \cite{michler1} to $N_1$ and
obtained a simple subgroup $\mathfrak G_1$ of $\GL_{782}(17)$
which he showed to be isomorphic to Fischer's sporadic group
$\Fi_{23}$. He also determined its character table, a faithful
permutation representation $PG_1$ of degree $31671$  and a
presentation of the centralizer $\mathfrak H_1 = C_{\mathfrak
G_1}(\mathfrak z_1)$ of a $2$-central involution $\mathfrak u_1$
of $\mathfrak G_1$. Since $\mathfrak H_1 \cong 2\Fi_{22}$ and
$C_{N_1}(u)$ have isomorphic Sylow $2$-subgroups by \cite{kim1}, we
determine in Lemma \ref{l. aut2Fi22} of Section 2 a presentation
for $A_1 = 2Aut(\Fi_{22})$, its character table, a system of
representatives of its conjugacy classes and a faithful
permutation representation $PA_1$. Thus we can show that $A_1$ and
$C_E(u)$ have isomorphic Sylow $2$-subgroups where $u$ is the non
$2$-central involution of $E$ mentioned above.

Lemma \ref{l. gensFi24} of Section 3 asserts that the amalgam
$\mathfrak G_1 \leftarrow \mathfrak H_1 \rightarrow A_1$ has $8$
distinct compatible pairs of semi-simple characters of $\mathfrak
G_1$ and $A_1$ of degree $8671$. All these compatible pairs have
the same semi-simple character $\tau = \tau_{\bf 3} + \tau_{\bf
4}$ of $\mathfrak G_1$. Its two irreducible constituents
$\tau_{\bf 3}$ and $\tau_{\bf 4}$ have respective degrees $3588$
and $5083$. Each compatible pair $(\tau,\chi_{i})$ determines a
subgroup $\mathfrak K_i$ of $\GL_{8671}(13)$. In the following two
sections we construct these $8$ matrix groups.

Since $\mathfrak G_1$ does not have any suitable subgroups whose
permutation characters of $\mathfrak G_1$ contain $\tau_{\bf 3}$
or $\tau_{\bf 4}$ as irreducible constituents we first construct a
pair $mH$ and $mE$ of subgroups of the permutation representation
$PG_1$ of $\mathfrak G_1$ such that $$PG_1 = \langle mH, mE
\rangle \quad \mbox{and} \quad mD = mH \cap mE \cong G_2(3) \times
Sym(3).$$ These subgroups enable us to construct the irreducible
representations of $PG_1$ in $\GL_{3588}(13)$ and $\GL_{5083}(13)$
corresponding to $\tau_{\bf 3}$ or $\tau_{\bf 4}$, respectively,
by means of Thompson's methods described in Theorem 7.2.2 of
\cite{michler}. The character tables of $mH$ and $mE$ are Tables
\ref{Fi_24 ct mH} and \ref{Fi_24 ct mE} of the Appendix,
respectively. They were automatically computed by MAGMA using
$PG_1$.

By the main result of \cite{kim1} we have $3$ permutations $py$,
$pq$ and $pw$ and a $2$-central involution $pz_1$ of $PG_1$ such
that $PG_1 = \langle py,pq,pw\rangle \cong \mathfrak G_1$ and
$C_{PG_1}(pz_1) = \langle py, pq\rangle \cong \mathfrak H_1$. We
found these two subgroups $mH$ and $mE$ of $PG_1$ by means of the
MAGMA command $\verb"LowIndexSubgroups(PG_1: n)"$ for $n = 137632$
and $n = 148642560$. An application of the first author's program
$\verb"GetShortGens(PG_1,U)"$ provides short words $pu$, $pv$,
$pr$ and $ps$ in terms of the generators $py$, $pq$ and $pw$ of
$PG_1$ such that $mD = \langle pu,pv\rangle$, $mH = \langle mD, pr
\rangle$ and $mE = \langle mD,ps \rangle$. For the construction of
the target matrix group $G$ it is necessary to get short words of
the old generators $py$, $pq$ and $pw$ in terms of the new
generators of $PG_1$. This is also done in Lemma \ref{l.
reps13Fi23} of Section 4. Thus we obtain the following diagram of
permutation groups.

$$ \diagram
& PG_1 = \langle py,pq,pw \rangle = \langle pu,pv,pr,ps \rangle& \\
mH = \langle pu,pv,pr \rangle \urto \drto & & mE = \langle pu,pv,ps \rangle \ulto \dlto \\
& mD = \langle pu,pv \rangle &
\enddiagram
$$

The restrictions of the characters $\tau_3$ and $\tau_4$ of $PG_1$
to the two subgroups $mH$ and $mE$ are determined in Lemma \ref{l.
reps13Fi23} of Section 4. Each of their irreducible constituents
is an irreducible constituent of a permutation character of $mH$
or $mE$. In Propositions \ref{prop. rep3588Fi23} and \ref{prop.
rep5083Fi23} we apply the first author's efficient implementation
of Algorithm 5.7.1 of \cite{michler} and get the various matrix
representations corresponding to the irreducible constituents of
the permutation characters of $mH$ and $mE$ determined in Lemma
\ref{l. reps13Fi23}. They enable us to construct the correct
blocked matrices $u$, $v$, $r$ and $s$ of the generators $pu$,
$pv$, $pr$ and $ps$ of $PG_1$ in $\GL_{8671}(13)$ corresponding to
the restrictions of the irreducible characters $\tau_3$ and
$\tau_4$ of $PG_1$ to $mH$, $mE$ and $mD$. Let $y$, $q$ and $w$ be
the matrices of $\GL_{8671}(13)$ obtained by inserting the
matrices $u$, $v$, $r$ and $s$ into the words of $py$, $pq$ and
$pw$ in terms of $pu$, $pv$, $pr$ and $ps$. Then $G_1 = \langle
y,q,w \rangle \cong PG_1$, and $H_1 = \langle y, q\rangle \cong
PH_1$.

Since the restriction of $\tau_{3} + \tau_{4}$ to $mD$ is not
multiplicity-free we have to solve a difficult amalgamation
problem in order to get the corresponding irreducible
representations of $PG_1$. By Lemma \ref{l. reps13Fi23} we know
that the generator $pr$ of $mH$ is an involution which commutes in
the known group $PG_1$ with the involution $pf =
(pu^3pvpupspv)^9$. Using this information and the structure of the
blocked matrices $r$, $s$, $u$ and $v$ we are able to solve the
amalgamation problems in Propositions \ref{prop. rep3588Fi23} and
\ref{prop. rep5083Fi23} by calculating the solutions of well
determined systems of algebraic equations with coefficients in
$GF(13)$.

As $H_1$ is isomorphic to a normal subgroup of index $2$ in $A_1$
we use Clifford's Theorem to determine exactly $8$ semi-simple
representations of $A_1$ corresponding to the $8$ compatible pairs
having the same restriction to $PH_1$ as $\tau_{\bf 3} \oplus
\tau_{\bf 4}$, see Proposition \ref{prop. repsA1}. In Remark
\ref{r. 8pairs} we show that for exactly one compatible pair of
Lemma \ref{l. gensFi24}(e) we can construct a matrix $t \in
\GL_{8671}(13)$ such that the matrix subgroup $G = \langle y,q,w,t
\rangle$ may have a Sylow $2$-subgroup isomorphic to the ones of
$E$. For this matrix $t$ we have the following amalgam of matrix
groups:

\vspace{2mm}

$$
\diagram
& G = \langle y, q, w, t\rangle& \\
\Fi_{23} \cong G_1 = \langle y, q,  w \rangle \urto & & 2\Fi_{22}:2 \cong \mathfrak A_1 = \langle y, q, t \rangle \ulto \\
& 2\Fi_{22} \cong H_1 = \langle y, q \rangle \ulto \urto &
\enddiagram
$$

\vspace{2mm}

In Section 6 we show that this matrix group $G$ is isomorphic to
the commutator subgroup $P'$ of finitely presented group $P$ of
Hall and Soicher \cite{hall}, see Lemma \ref{l. nicepres.Fi_24}.
Thus $G \cong \Fi_{24}'$. In particular, the subgroup $G = \langle
q, y, w, t\rangle$ of $\GL_{8671}(13)$ is a simple group of order
$2^{21} \cdot 3^{16} \cdot 5^2 \cdot 7^3 \cdot 11 \cdot 13 \cdot
17 \cdot 23 \cdot 29$ and it has a faithful permutation
representation of degree $306936$ with stabilizer $G_1 = \langle
q, y, w \rangle$, see Theorem \ref{thm. existenceFi_24}.

In Section 7 we determine generators and a presentation of $H =
C_G(z)$ of a $2$-central involution $z$ of $G$, see
Proposition \ref{prop. 2-centralFi_24}. Furthermore, we construct a
faithful permutation representation of degree $258048$ of $H$ with
a documented stabilizer. It has been used to calculate a system of
representatives of its $167$ conjugacy classes and its character
table, see Tables \ref{Fi_24 cc H} and \ref{Fi_24 ct H}.,
respectively.

In Section 8 we show that $G$ has 2 conjugacy classes of
involutions. They are represented by $u = [(y(y^5t)^7]^{14}$ and
$z = (xyw)^8$. Their centralizers are $C_G(u) = \langle
q,y,t\rangle = A_1$ and $C_G(z) = H$. Using their character tables
we calculate the group order of $G$ using Thompson's group
order formula and Theorem 6.1.4 of \cite{michler2}, see
Proposition \ref{prop. Fi24'order}.

In the appendix we collect all the systems of representatives of
conjugacy classes in terms of the given generators of the local
subgroups of $G$ which have been used to construct this matrix
group $G \cong \Fi_{24}'$. We also state the character tables of
these subgroups. The four generating matrices of the simple
subgroup $\mathfrak G = \langle \mathfrak q, \mathfrak y,
\mathfrak t, \mathfrak w\rangle$ of $\GL_{8671}(13)$ can be
downloaded from the first author's website\\
$\verb"http://www.math.yale.edu/~hk47/Fi24/index.html"$.

Concerning our notation and terminology we refer to the books
\cite{carter} and \cite{michler}. The computer algebra system
MAGMA is described in Cannon-Playoust \cite{magma}.

Acknowledgements: The first author kindly acknowledges financial
support by the Hunter R. Rawlings III Cornell Presidential
Research Scholars Program for participation in the research
project {\em Simultaneous construction of the sporadic simple
groups of Conway, Fischer and Janko}. This collaboration has also
been supported by the grant NSF/SCREMS DMS-0532/06.

\newpage

\section{An extension of Mathieu's group $\M_{24}$ }

In \cite{fischer1} B. Fischer stated that his simple group
$\Fi_{24}'$ has a subgroup $E$ which is isomorphic to a non
split extension of the Mathieu group $\M_{24}$ by a well
determined $11$-dimensional irreducible $\M_{24}$-module $V$ over $GF(2)$ such that
$E$ contains a Sylow $2$-subgroup $S$ of his simple group
$\Fi_{24}'$. In this section we construct such an extension by
means of Holt's Algorithm \cite{holt5} implemented in MAGMA. It
follows that $E$ is uniquely determined by $\M_{24}$ up to
isomorphism.

A faithful permutation representation of degree $24$ of $\M_{24}$ is stated in Lemma 8.2.2 of
\cite{michler}. The irreducible $2$-modular representations of the
Mathieu group $\M_{24}$ were determined by G. James \cite{james}.
Here only the $2$ non isomorphic simple modules $V_i$, $i = 1,2$,
of dimension $11$ over $F = GF(2)$ will be used. Todd's
permutation representations of the Mathieu groups are stated in
Lemma 8.2.2 of \cite{michler}. Therefore all conditions of Holt's
Algorithm \cite{holt5} implemented in MAGMA are satisfied. It
constructs all split and non split extensions of $\M_{24}$ by
$V_1$ and $V_2$ up to isomorphism. Here only the presentation of
the non split extension of $\M_{24}$ by $V_2$ is stated.

\begin{lemma}\label{l. M24-extensions}
Let $\M_{24} = \langle a,b,c,d,t,g,h,i,j,k \rangle$ be the
finitely presented group of Definition 8.2.1 of \cite{michler}.
Then the following statements hold:

\begin{enumerate}

\item[\rm(a)] The first irreducible representation $V_1$ of
$\M_{24}$ is described by the following matrices:
{\renewcommand{\arraystretch}{0.5}
$$
a_1 = \left( \begin{array}{*{11}{c@{\,}}c}
1 & 0 & 0 & 0 & 0 & 0 & 0 & 0 & 0 & 0 & 0\\
0 & 1 & 0 & 0 & 0 & 0 & 0 & 0 & 1 & 0 & 0\\
0 & 0 & 1 & 0 & 0 & 0 & 0 & 0 & 1 & 0 & 0\\
0 & 0 & 0 & 1 & 0 & 0 & 0 & 0 & 1 & 0 & 0\\
0 & 0 & 0 & 0 & 0 & 1 & 0 & 0 & 0 & 0 & 0\\
0 & 0 & 0 & 0 & 1 & 0 & 0 & 0 & 0 & 0 & 0\\
0 & 0 & 0 & 0 & 0 & 0 & 0 & 0 & 0 & 1 & 0\\
0 & 0 & 0 & 0 & 0 & 0 & 0 & 0 & 1 & 0 & 1\\
0 & 0 & 0 & 0 & 0 & 0 & 0 & 0 & 1 & 0 & 0\\
0 & 0 & 0 & 0 & 0 & 0 & 1 & 0 & 0 & 0 & 0\\
0 & 0 & 0 & 0 & 0 & 0 & 0 & 1 & 1 & 0 & 0
\end{array} \right),\quad
b_1 = \left( \begin{array}{*{11}{c@{\,}}c}
1 & 0 & 0 & 0 & 0 & 0 & 0 & 1 & 1 & 0 & 1\\
0 & 1 & 0 & 0 & 0 & 0 & 0 & 1 & 1 & 0 & 1\\
0 & 0 & 1 & 0 & 0 & 0 & 0 & 1 & 1 & 0 & 1\\
0 & 0 & 0 & 1 & 0 & 0 & 0 & 0 & 1 & 0 & 0\\
0 & 0 & 0 & 0 & 0 & 0 & 1 & 1 & 0 & 0 & 0\\
0 & 0 & 0 & 0 & 0 & 0 & 0 & 0 & 1 & 1 & 1\\
0 & 0 & 0 & 0 & 1 & 0 & 0 & 0 & 0 & 0 & 1\\
0 & 0 & 0 & 0 & 0 & 0 & 0 & 0 & 0 & 0 & 1\\
0 & 0 & 0 & 0 & 0 & 0 & 0 & 0 & 1 & 0 & 0\\
0 & 0 & 0 & 0 & 0 & 1 & 0 & 1 & 1 & 0 & 0\\
0 & 0 & 0 & 0 & 0 & 0 & 0 & 1 & 0 & 0 & 0
\end{array} \right),
$$
}

{\renewcommand{\arraystretch}{0.5}
$$
c_1 = \left( \begin{array}{*{11}{c@{\,}}c}
1 & 0 & 0 & 0 & 0 & 0 & 1 & 0 & 0 & 1 & 0\\
0 & 1 & 0 & 0 & 0 & 0 & 1 & 0 & 0 & 1 & 0\\
0 & 0 & 1 & 0 & 0 & 0 & 1 & 0 & 1 & 1 & 0\\
0 & 0 & 0 & 1 & 0 & 0 & 0 & 0 & 1 & 0 & 0\\
0 & 0 & 0 & 0 & 0 & 0 & 1 & 1 & 1 & 0 & 0\\
0 & 0 & 0 & 0 & 0 & 0 & 0 & 0 & 0 & 1 & 1\\
0 & 0 & 0 & 0 & 0 & 0 & 0 & 0 & 1 & 1 & 0\\
0 & 0 & 0 & 0 & 1 & 0 & 0 & 0 & 0 & 1 & 0\\
0 & 0 & 0 & 0 & 0 & 0 & 0 & 0 & 1 & 0 & 0\\
0 & 0 & 0 & 0 & 0 & 0 & 1 & 0 & 1 & 0 & 0\\
0 & 0 & 0 & 0 & 0 & 1 & 1 & 0 & 1 & 0 & 0
\end{array} \right),\quad
d_1 = \left( \begin{array}{*{11}{c@{\,}}c}
1 & 0 & 0 & 0 & 0 & 0 & 1 & 0 & 0 & 1 & 0\\
0 & 1 & 0 & 0 & 0 & 1 & 1 & 0 & 0 & 0 & 1\\
0 & 0 & 1 & 0 & 0 & 1 & 1 & 1 & 0 & 0 & 0\\
0 & 0 & 0 & 1 & 0 & 1 & 1 & 1 & 0 & 0 & 0\\
0 & 0 & 0 & 0 & 0 & 0 & 0 & 1 & 1 & 1 & 0\\
0 & 0 & 0 & 0 & 0 & 0 & 1 & 0 & 0 & 0 & 1\\
0 & 0 & 0 & 0 & 0 & 0 & 0 & 1 & 0 & 1 & 1\\
0 & 0 & 0 & 0 & 0 & 1 & 0 & 0 & 0 & 1 & 0\\
0 & 0 & 0 & 0 & 1 & 1 & 1 & 1 & 0 & 1 & 1\\
0 & 0 & 0 & 0 & 0 & 0 & 1 & 1 & 0 & 0 & 1\\
0 & 0 & 0 & 0 & 0 & 1 & 0 & 1 & 0 & 1 & 1
\end{array} \right),
$$
}

{\renewcommand{\arraystretch}{0.5}
$$
t_1 = \left( \begin{array}{*{11}{c@{\,}}c}
1 & 0 & 0 & 0 & 0 & 0 & 1 & 0 & 0 & 1 & 0 \\
0 & 1 & 0 & 0 & 0 & 1 & 0 & 0 & 1 & 1 & 0 \\
0 & 0 & 1 & 0 & 0 & 0 & 0 & 1 & 0 & 1 & 0 \\
0 & 0 & 0 & 1 & 0 & 0 & 0 & 0 & 1 & 0 & 0 \\
0 & 0 & 0 & 0 & 1 & 0 & 1 & 1 & 1 & 1 & 0 \\
0 & 0 & 0 & 0 & 0 & 1 & 0 & 1 & 1 & 0 & 0 \\
0 & 0 & 0 & 0 & 0 & 1 & 1 & 0 & 1 & 0 & 0 \\
0 & 0 & 0 & 0 & 0 & 0 & 1 & 1 & 1 & 0 & 0 \\
0 & 0 & 0 & 0 & 0 & 1 & 1 & 1 & 1 & 0 & 0 \\
0 & 0 & 0 & 0 & 0 & 1 & 0 & 1 & 0 & 1 & 1 \\
0 & 0 & 0 & 0 & 0 & 1 & 1 & 0 & 0 & 1 & 0
\end{array} \right),\quad
g_1 = \left( \begin{array}{*{11}{c@{\,}}c}
1 & 0 & 0 & 0 & 0 & 0 & 0 & 0 & 0 & 1 & 1 \\
0 & 1 & 0 & 0 & 0 & 0 & 0 & 0 & 1 & 0 & 0 \\
0 & 0 & 1 & 0 & 0 & 0 & 0 & 0 & 0 & 0 & 1 \\
0 & 0 & 0 & 0 & 1 & 0 & 0 & 0 & 1 & 1 & 1 \\
0 & 0 & 0 & 1 & 0 & 0 & 0 & 0 & 1 & 1 & 1 \\
0 & 0 & 0 & 0 & 0 & 1 & 0 & 0 & 1 & 1 & 1 \\
0 & 0 & 0 & 0 & 0 & 0 & 1 & 0 & 1 & 1 & 0 \\
0 & 0 & 0 & 0 & 0 & 0 & 0 & 1 & 1 & 0 & 1 \\
0 & 0 & 0 & 0 & 0 & 0 & 0 & 0 & 1 & 0 & 0 \\
0 & 0 & 0 & 0 & 0 & 0 & 0 & 0 & 0 & 1 & 0 \\
0 & 0 & 0 & 0 & 0 & 0 & 0 & 0 & 0 & 0 & 1
\end{array} \right),
$$
}

{\renewcommand{\arraystretch}{0.5}
$$
h_1 = \left( \begin{array}{*{11}{c@{\,}}c}
1 & 0 & 0 & 0 & 0 & 0 & 1 & 0 & 0 & 1 & 0 \\
0 & 1 & 0 & 0 & 0 & 0 & 0 & 0 & 0 & 1 & 0 \\
0 & 0 & 0 & 1 & 0 & 0 & 0 & 0 & 0 & 1 & 0 \\
0 & 0 & 1 & 0 & 0 & 0 & 0 & 0 & 0 & 1 & 0 \\
0 & 0 & 0 & 0 & 1 & 0 & 1 & 0 & 0 & 0 & 0 \\
0 & 0 & 0 & 0 & 0 & 1 & 0 & 0 & 0 & 1 & 0 \\
0 & 0 & 0 & 0 & 0 & 0 & 1 & 0 & 0 & 0 & 0 \\
0 & 0 & 0 & 0 & 0 & 0 & 1 & 0 & 0 & 0 & 1 \\
0 & 0 & 0 & 0 & 0 & 0 & 1 & 0 & 1 & 1 & 0 \\
0 & 0 & 0 & 0 & 0 & 0 & 0 & 0 & 0 & 1 & 0 \\
0 & 0 & 0 & 0 & 0 & 0 & 1 & 1 & 0 & 0 & 0
\end{array} \right),\quad
i_1 = \left( \begin{array}{*{11}{c@{\,}}c}
1 & 0 & 0 & 0 & 0 & 0 & 0 & 0 & 0 & 0 & 0 \\
0 & 0 & 1 & 0 & 0 & 0 & 1 & 1 & 0 & 0 & 0 \\
0 & 1 & 0 & 0 & 0 & 1 & 1 & 0 & 1 & 0 & 0 \\
0 & 0 & 0 & 1 & 0 & 1 & 1 & 1 & 0 & 0 & 0 \\
0 & 0 & 0 & 0 & 1 & 1 & 0 & 1 & 1 & 0 & 0 \\
0 & 0 & 0 & 0 & 0 & 0 & 0 & 1 & 1 & 0 & 0 \\
0 & 0 & 0 & 0 & 0 & 1 & 0 & 1 & 0 & 0 & 0 \\
0 & 0 & 0 & 0 & 0 & 0 & 1 & 1 & 1 & 0 & 0 \\
0 & 0 & 0 & 0 & 0 & 1 & 1 & 1 & 1 & 0 & 0 \\
0 & 0 & 0 & 0 & 0 & 1 & 0 & 0 & 1 & 0 & 1 \\
0 & 0 & 0 & 0 & 0 & 1 & 1 & 0 & 0 & 1 & 0
\end{array} \right),
$$
} {\renewcommand{\arraystretch}{0.5}
$$
j_1 = \left( \begin{array}{*{11}{c@{\,}}c}
0 & 1 & 0 & 0 & 0 & 1 & 0 & 0 & 0 & 0 & 0 \\
1 & 0 & 0 & 0 & 0 & 1 & 0 & 0 & 0 & 0 & 0 \\
0 & 0 & 1 & 0 & 0 & 1 & 0 & 0 & 1 & 0 & 0 \\
0 & 0 & 0 & 1 & 0 & 0 & 0 & 0 & 1 & 0 & 0 \\
0 & 0 & 0 & 0 & 1 & 0 & 0 & 0 & 1 & 0 & 0 \\
0 & 0 & 0 & 0 & 0 & 1 & 0 & 0 & 0 & 0 & 0 \\
0 & 0 & 0 & 0 & 0 & 1 & 1 & 0 & 1 & 0 & 0 \\
0 & 0 & 0 & 0 & 0 & 1 & 0 & 1 & 0 & 0 & 0 \\
0 & 0 & 0 & 0 & 0 & 0 & 0 & 0 & 1 & 0 & 0 \\
0 & 0 & 0 & 0 & 0 & 0 & 0 & 0 & 1 & 0 & 1 \\
0 & 0 & 0 & 0 & 0 & 0 & 0 & 0 & 1 & 1 & 0
\end{array} \right)\quad \text{and}~\quad
k_1 = \left( \begin{array}{*{11}{c@{\,}}c}
1 & 1 & 1 & 1 & 1 & 1 & 1 & 0 & 1 & 1 & 1 \\
0 & 1 & 0 & 0 & 0 & 1 & 0 & 1 & 0 & 0 & 0 \\
0 & 0 & 1 & 0 & 0 & 1 & 1 & 1 & 0 & 0 & 0 \\
0 & 0 & 0 & 1 & 0 & 1 & 1 & 1 & 0 & 0 & 0 \\
0 & 0 & 0 & 0 & 1 & 0 & 1 & 0 & 0 & 0 & 0 \\
0 & 0 & 0 & 0 & 0 & 0 & 0 & 1 & 0 & 0 & 0 \\
0 & 0 & 0 & 0 & 0 & 0 & 1 & 0 & 0 & 0 & 0 \\
0 & 0 & 0 & 0 & 0 & 1 & 0 & 0 & 0 & 0 & 0 \\
0 & 0 & 0 & 0 & 0 & 1 & 1 & 1 & 1 & 0 & 0 \\
0 & 0 & 0 & 0 & 0 & 0 & 1 & 1 & 0 & 0 & 1 \\
0 & 0 & 0 & 0 & 0 & 1 & 1 & 0 & 0 & 1 & 0
\end{array} \right).
$$
}

\item[\rm(b)] The second irreducible representation $V_2$ of
$\M_{24}$ is described by the transpose inverse matrices of the
generating matrices of $\M_{24}$ defining $V_1$:

$a_2 = [a_1^{-1}]^{T}$, $b_2 = [b_1^{-1}]^{T}$, $c_2 =
[c_1^{-1}]^{T}$, $d_2 = [d_1^{-1}]^{T}$, $t_2 = [t_1^{-1}]^{T}$,
$g_2 = [g_1^{-1}]^{T}$, $h_2 = [h_1^{-1}]^{T}$, $i_2 =
[i_1^{-1}]^{T}$, $j_2 = [j_1^{-1}]^{T}$, and $k_2 =
[k_1^{-1}]^{T}$.

\item[\rm(c)] $dim_F[H^2(\M_{24},V_1)] = 0$ and
$dim_F[H^2(\M_{24},V_2)] = 1$.

In particular, there is a uniquely determined non split extension
$E$ of $\M_{24}$ by $V_2$.

\item[\rm(d)] The non split extension $$E = E(Fi_{24}') = \langle
a,b,c,d,t,g,h,i,j,k,v_1,v_2,v_3,v_4,v_5,v_6,v_8,v_8,v_9,v_{10},v_{11}
\rangle$$ of $\M_{24}$ by $V_2$  has a set ${\mathcal R}(E)$ of
defining relations consisting of ${\mathcal R_1}(V_2\rtimes
\M_{24})$ and the following relations:

\begin{eqnarray*}
&&v_i^2 = 1 \quad \mbox{and} \quad v_kv_j = v_jv_k \quad \mbox{for all} \quad 1 \le i, j, k \le 11.\\
&&a^{-1}v_1 av_1^{-1}v_2^{-1}v_5^{-1}v_6^{-1}v_7^{-1}v_8^{-1}v_{10}^{-1} = a^{-1}v_2 av_{10}^{-1} = a^{-1}v_3 av_4^{-1} =1,\\
&&a^{-1}v_4 av_3^{-1} = a^{-1}v_5 av_6^{-1} =a^{-1}v_6 av_5^{-1} =a^{-1}v_7 av_8^{-1} =a^{-1}v_8 av_7^{-1} = 1,\\
&&a^{-1}v_9 av_9^{-1} = a^{-1}v_{10} av_2^{-1} =a^{-1}v_{11} av_{11}^{-1} = b^{-1}v_1 bv_6^{-1} =b^{-1}v_2 bv_8^{-1} = 1,\\
&&b^{-1}v_3 bv_1^{-1}v_2^{-1}v_4^{-1}v_6^{-1}v_{10}^{-1}v_{11}^{-1} = b^{-1}v_4 bv_1^{-1}v_3^{-1}v_6^{-1}v_7^{-1}v_8^{-1}v_{11}^{-1} = 1, \\
&&b^{-1}v_5 bv_1^{-1}v_2^{-1}v_5^{-1}v_6^{-1}v_7^{-1}v_8^{-1}v_{10}^{-1} = b^{-1}v_6 bv_1^{-1} =b^{-1}v_7 bv_{10}^{-1} =1,\\
&&b^{-1}v_8 bv_2^{-1} = b^{-1}v_9 bv_9^{-1} =b^{-1}v_{10} bv_7^{-1} =b^{-1}v_{11} bv_{11}^{-1} = c^{-1}v_1 cv_{10}^{-1} =1,\\
&&c^{-1}v_2 cv_1^{-1}v_2^{-1}v_5^{-1}v_6^{-1}v_7^{-1}v_8^{-1}v_{10}^{-1} = c^{-1}v_3 cv_2^{-1}v_4^{-1}v_5^{-1}v_8^{-1}v_9^{-1}v_{10}^{-1} = 1, \\
&&c^{-1}v_4 cv_2^{-1}v_3^{-1}v_6^{-1}v_7^{-1}v_9^{-1}v_{10}^{-1} = c^{-1}v_5 cv_8^{-1} =c^{-1}v_6 cv_7^{-1} =c^{-1}v_7 cv_6^{-1} = 1, \\
&&c^{-1}v_8 cv_5^{-1} =c^{-1}v_9 cv_9^{-1} =c^{-1}v_{10} cv_1^{-1} = c^{-1}v_{11} cv_{11}^{-1} =1,\\
&&d^{-1}v_1 dv_3^{-1}v_5^{-1}v_7^{-1}v_8^{-1}v_9^{-1}v_{10}^{-1}v_{11}^{-1} =d^{-1}v_2 dv_1^{-1}v_3^{-1}v_6^{-1}v_7^{-1}v_8^{-1}v_{11}^{-1} =1,\\
&&d^{-1}v_3 dv_7^{-1} =d^{-1}v_4 dv_8^{-1} = d^{-1}v_5 dv_2^{-1}v_3^{-1}v_6^{-1}v_7^{-1}v_9^{-1}v_{10}^{-1} = 1, \\
&&d^{-1}v_6 dv_2^{-1}v_4^{-1}v_5^{-1}v_8^{-1}v_9^{-1}v_{10}^{-1} = d^{-1}v_7 dv_3^{-1} =d^{-1}v_8 dv_4^{-1} =1,\\
\end{eqnarray*}
\begin{eqnarray*}
&&d^{-1}v_9 dv_9^{-1} = d^{-1}v_{10} dv_1^{-1}v_2^{-1}v_4^{-1}v_6^{-1}v_{10}^{-1}v_{11}^{-1} = d^{-1}v_{11} dv_{11}^{-1} =1,\\
&&t^{-1}v_1 tv_6^{-1}v_{11}^{-1} = t^{-1}v_2 tv_{10}^{-1}v_{11}^{-1} =t^{-1}v_3 tv_7^{-1}v_{11}^{-1} = 1,\\
&&t^{-1}v_4 tv_2^{-1}v_4^{-1}v_5^{-1}v_8^{-1}v_9^{-1}v_{10}^{-1}v_{11}^{-1} = t^{-1}v_5 tv_8^{-1}v_{11}^{-1} =1,\\
&&t^{-1}v_6tv_2^{-1}v_3^{-1}v_6^{-1}v_7^{-1}v_9^{-1}v_{10}^{-1}v_{11}^{-1} = 1,\\
&&t^{-1}v_7 tv_1^{-1}v_2^{-1}v_5^{-1}v_6^{-1}v_7^{-1}v_8^{-1}v_{10}^{-1}v_{11}^{-1} = 1,\\
&&t^{-1}v_8 tv_1^{-1}v_3^{-1}v_6^{-1}v_7^{-1}v_8^{-1} =t^{-1}v_9 tv_{11}^{-1} =1,\\
&&t^{-1}v_{10} tv_3^{-1}v_5^{-1}v_7^{-1}v_8^{-1}v_9^{-1}v_{10}^{-1} =t^{-1}v_{11} tv_9^{-1}v_{11}^{-1} = g^{-1}v_1 gv_4^{-1}v_6^{-1} =1,\\
&&g^{-1}v_2 gv_4^{-1}v_{10}^{-1} = g^{-1}v_3 gv_3^{-1}v_4^{-1} =g^{-1}v_4 gv_4^{-1} = g^{-1}v_5 gv_4^{-1}v_8^{-1} =1,\\
&&g^{-1}v_6 gv_1^{-1}v_4^{-1} = g^{-1}v_7 gv_1^{-1}v_2^{-1}v_4^{-1}v_5^{-1}v_6^{-1}v_7^{-1}v_8^{-1}v_{10}^{-1} = 1,\\
&&g^{-1}v_8 gv_4^{-1}v_5^{-1} =g^{-1}v_9 gv_2^{-1}v_5^{-1}v_8^{-1}v_9^{-1}v_{10}^{-1} = g^{-1}v_{10} gv_2^{-1}v_4^{-1} =1,\\
&&g^{-1}v_{11} gv_1^{-1}v_2^{-1}v_6^{-1}v_{10}^{-1}v_{11}^{-1} =h^{-1}v_1 hv_1^{-1}v_9^{-1} = h^{-1}v_2 hv_9^{-1}v_{10}^{-1} =1,\\
&&h^{-1}v_3 hv_7^{-1}v_9^{-1} = h^{-1}v_4 hv_8^{-1}v_9^{-1} =h^{-1}v_5 hv_2^{-1}v_4^{-1}v_5^{-1}v_8^{-1}v_{10}^{-1} =1, \\
&&h^{-1}v_6 hv_2^{-1}v_3^{-1}v_6^{-1}v_7^{-1}v_{10}^{-1} =h^{-1}v_7 hv_3^{-1}v_9^{-1} = 1,\\
&&h^{-1}v_8 hv_4^{-1}v_9^{-1} = h^{-1}v_9 hv_9^{-1} =h^{-1}v_{10} hv_2^{-1}v_9^{-1} = 1,\\
&&h^{-1}v_{11} hv_9^{-1}v_{11}^{-1} = i^{-1}v_1 iv_2^{-1}v_3^{-1}v_6^{-1}v_7^{-1}v_{10}^{-1} = 1,\\
&&i^{-1}v_2 iv_3^{-1}v_5^{-1}v_7^{-1}v_8^{-1}v_{10}^{-1}v_{11}^{-1} = i^{-1}v_3 iv_5^{-1}v_9^{-1} = 1,\\
&&i^{-1}v_4 iv_4^{-1}v_9^{-1} = i^{-1}v_5 iv_3^{-1}v_9^{-1} =i^{-1}v_6 iv_6^{-1}v_9^{-1} = 1,\\
&&i^{-1}v_7iv_1^{-1}v_3^{-1}v_6^{-1}v_7^{-1}v_8^{-1}v_9^{-1}v_{11}^{-1} = i^{-1}v_8iv_1^{-1}v_2^{-1}v_5^{-1}v_6^{-1}v_7^{-1}v_8^{-1}v_9^{-1}v_{10}^{-1} = 1,\\
&&i^{-1}v_9 iv_9^{-1} =i^{-1}v_{10} iv_9^{-1}v_{10}^{-1} =i^{-1}v_{11} i v_9^{-1}v_{11}^{-1} =j^{-1}v_1 jv_6^{-1} =1,\\
&&j^{-1}v_2 jv_1^{-1}v_2^{-1}v_5^{-1}v_6^{-1}v_7^{-1}v_8^{-1}v_{10}^{-1} = 1,\\
&&j^{-1}v_3 jv_3^{-1}v_5^{-1}v_7^{-1}v_8^{-1}v_9^{-1}v_{10}^{-1}v_{11}^{-1} = 1,\\
&&j^{-1}v_4 jv_4^{-1} =j^{-1}v_5 jv_8^{-1} = j^{-1}v_6 jv_1^{-1} =j^{-1}v_7 jv_{10}^{-1} =1,\\
&&j^{-1}v_8 jv_5^{-1} = j^{-1}v_9 jv_{11}^{-1} = j^{-1}v_{10} jv_7^{-1} = j^{-1}v_{11} jv_9^{-1} = 1,\\
&&k^{-1}v_1 kv_3^{-1}v_5^{-1}v_7^{-1}v_8^{-1}v_{10}^{-1}v_{11}^{-1} = 1,\\
&&k^{-1}v_2 kv_2^{-1}v_3^{-1}v_6^{-1}v_7^{-1}v_{10}^{-1} = k^{-1}v_3 kv_7^{-1}v_9^{-1} = 1,\\
&&k^{-1}v_4 kv_4^{-1}v_9^{-1} = k^{-1}v_5 kv_1^{-1}v_3^{-1}v_6^{-1}v_7^{-1}v_8^{-1}v_9^{-1}v_{11}^{-1} = 1,\\
&&k^{-1}v_6 kv_9^{-1}v_{10}^{-1} = k^{-1}v_7 kv_3^{-1}v_9^{-1} =k^{-1}v_8 kv_8^{-1}v_9^{-1} =k^{-1}v_9 kv_9^{-1} = 1, \\
&&k^{-1}v_{10} kv_6^{-1}v_9^{-1} =k^{-1}v_{11} kv_9^{-1}v_{11}^{-1} = a^2v_5^{-1}v_6^{-1}v_7^{-1}v_8^{-1} =1,\\
&&b^2v_1^{-1}v_3^{-1}v_4^{-1}v_6^{-1}v_7^{-1}v_8^{-1}v_{11}^{-1} = c^2v_1^{-1}v_5^{-1}v_6^{-1}v_7^{-1}v_8^{-1}v_{10}^{-1} = 1, \\
&&d^2v_1^{-1}v_2^{-1}v_3^{-1}v_6^{-1}v_7^{-1}v_8^{-1}v_{11}^{-1} = b^{-1} a ba^{-1}v_1^{-1}v_3^{-1}v_4^{-1}v_6^{-1}v_7^{-1}v_{10}^{-1} = 1,\\
&&c^{-1} a c a^{-1}v_2^{-1}v_{10}^{-1} = d^{-1} a d a^{-1}v_5^{-1}v_6^{-1}v_9^{-1} = 1, \\
&&c^{-1} b cb^{-1}v_1^{-1}v_3^{-1}v_4^{-1}v_8^{-1}v_9^{-1} = d^{-1} b d b^{-1}v_2^{-1}v_4^{-1} = 1,\\
&&d^{-1} c d c^{-1}v_2^{-1}v_3^{-1}v_5^{-1}v_{10}^{-1}v_{11}^{-1} = t^3v_1^{-1}v_3^{-1}v_8^{-1} = 1,\\
&&t^{-1} a t d^{-1} c^{-1}v_2^{-1}v_3^{-1}v_4^{-1}v_5^{-1}v_7^{-1}v_8^{-1}v_9^{-1}v_{11}^{-1} = 1,\\
&&t^{-1} b t d^{-1} a^{-1}v_2^{-1}v_3^{-1}v_5^{-1}v_6^{-1}v_9^{-1}v_{11}^{-1} =1,\\
&&t^{-1} c t d^{-1} b^{-1}v_1^{-1}v_2^{-1}v_4^{-1}v_5^{-1}v_8^{-1}v_9^{-1}v_{11}^{-1} = 1,\\
&&t^{-1} d t c^{-1} b^{-1}a^{-1}v_2^{-1}v_3^{-1}v_5^{-1}v_7^{-1}v_8^{-1}v_9^{-1}v_{11}^{-1} = 1,\\
&& g^2v_5^{-1}v_8^{-1} = g a g a g av_2^{-1}v_3^{-1}v_4^{-1}v_7^{-1}v_8^{-1}v_{11}^{-1} = 1,\\
&&g b g b g bv_1^{-1}v_3^{-1}v_4^{-1}v_8^{-1}v_9^{-1}v_{10}^{-1} = g c g c g cv_1^{-1}v_2^{-1}v_5^{-1}v_6^{-1}v_9^{-1} = 1,\\
\end{eqnarray*}
\begin{eqnarray*}
&&g t g tv_4^{-1}v_{11}^{-1} =h^2v_4^{-1}v_8^{-1}v_9^{-1} = h^{-1} a h a^{-1}v_2^{-1}v_3^{-1}v_4^{-1}v_9^{-1}v_{10}^{-1} = 1,\\
&&h^{-1} b h d^{-1} b^{-1} a^{-1}v_4^{-1}v_5^{-1}v_7^{-1}v_8^{-1}v_{10}^{-1} = 1, \\
&&h^{-1} c h c^{-1} a^{-1}v_1^{-1}v_3^{-1}v_4^{-1}v_5^{-1}v_7^{-1}v_{10}^{-1}v_{11}^{-1} = 1,\\
&&h^{-1} d h d^{-1}v_2^{-1}v_4^{-1}v_8^{-1}v_9^{-1}v_{10}^{-1} = 1,\\
&&h^{-1} t h tv_1^{-1}v_2^{-1}v_3^{-1}v_4^{-1}v_5^{-1}v_6^{-1}v_7^{-1} v_{10}^{-1}v_{11}^{-1} =1,\\
&&g h g h g hv_1^{-1}v_3^{-1}v_4^{-1}v_6^{-1}v_8^{-1}v_9^{-1} = 1,\\
&&i^2v_1^{-1}v_5^{-1}v_6^{-1}v_8^{-1}v_9^{-1}v_{11}^{-1} = i^{-1} a i d^{-1}c^{-1}v_3^{-1}v_5^{-1}v_7^{-1}v_9^{-1}v_{10}^{-1} = 1,\\
&&i^{-1} b i d^{-1} a^{-1}v_2^{-1}v_3^{-1}v_4^{-1}v_6^{-1}v_9^{-1}v_{11}^{-1} = 1,\\
&&i^{-1} c i d^{-1} c^{-1} b^{-1} a^{-1}v_3^{-1}v_5^{-1}v_7^{-1}v_9^{-1}v_{10}^{-1}v_{11}^{-1} = 1, \\
&&i^{-1} d i d^{-1} c^{-1} b^{-1}v_1^{-1}v_3^{-1}v_4^{-1}v_9^{-1}v_{10}^{-1}v_{11}^{-1} = 1,\\
&&i^{-1} t i tv_1^{-1}v_5^{-1}v_6^{-1}v_8^{-1} = i^{-1} g i g^{-1} t^{-1}v_1^{-1}v_2^{-1}v_3^{-1}v_8^{-1}v_{10}^{-1}v_{11}^{-1} = 1,\\
&&h i h i h iv_1^{-1}v_2^{-1}v_4^{-1}v_5^{-1}v_6^{-1}v_7^{-1}v_{10}^{-1}v_{11}^{-1} = j^2v_9^{-1}v_{11}^{-1} = 1,\\
&&j^{-1} a j c^{-1} b^{-1} a^{-1}v_1^{-1}v_3^{-1}v_5^{-1}v_6^{-1}v_7^{-1}v_8^{-1}v_9^{-1} = 1,\\
&&j^{-1} b j b^{-1}v_9^{-1}v_{11}^{-1} = j^{-1} c j c^{-1}v_1^{-1}v_5^{-1}v_6^{-1}v_7^{-1}v_8^{-1}v_{10}^{-1} = 1,\\
&&j^{-1} d j d^{-1}c^{-1}v_2^{-1}v_5^{-1}v_8^{-1} = j^{-1} t j tv_3^{-1}v_5^{-1}v_8^{-1}v_{10}^{-1} =1,\\
&&j^{-1} g jg^{-1}v_2^{-1}v_4^{-1}v_5^{-1}v_7^{-1}v_8^{-1}v_9^{-1}v_{11}^{-1} = j^{-1} h j h^{-1} t^{-1}v_3^{-1}v_4^{-1}v_5^{-1}v_6^{-1}v_8^{-1}v_{11}^{-1} = 1,\\
&&i j i j i jv_1^{-1}v_3^{-1}v_4^{-1}v_5^{-1}v_7^{-1}v_8^{-1} =k^2v_1^{-1}v_5^{-1}v_6^{-1}v_8^{-1}v_9^{-1}v_{11}^{-1} = 1,\\
&&k^{-1} a k d^{-1}a^{-1}v_4^{-1}v_6^{-1}v_8^{-1}v_9^{-1} = k^{-1} b k d^{-1} c^{-1}v_1^{-1}v_2^{-1}v_7^{-1}v_9^{-1} = 1,\\
&&k^{-1} c k d^{-1} b^{-1}v_2^{-1}v_3^{-1}v_5^{-1}v_8^{-1}v_9^{-1}v_{10}^{-1}v_{11}^{-1} = 1, \\
&&k^{-1} d k d^{-1}v_2^{-1}v_3^{-1}v_5^{-1}v_7^{-1} = k^{-1} t k tv_2^{-1}v_5^{-1}v_6^{-1}v_8^{-1}v_9^{-1}v_{11}^{-1} = 1,\\
&&k^{-1} g k g^{-1}t^{-1}v_2^{-1}v_3^{-1}v_4^{-1}v_5^{-1}v_6^{-1}v_9^{-1}v_{10} = k^{-1} h k h^{-1}v_2^{-1}v_4^{-1}v_8^{-1}v_{10}^{-1} = 1,\\
&&k^{-1} i k i^{-1}v_2^{-1}v_7^{-1}v_8^{-1}v_{10}^{-1}v_{11}^{-1} = j k j k j kv_2^{-1}v_5^{-1}v_8^{-1}v_9^{-1}v_{10}^{-1} = 1. \\
\end{eqnarray*}

\item[\rm(e)] $E$ has a faithful permutation representation $PE$
of degree $1518$ with stabilizer $T = \langle g, h, i, (dg)^5,
(dhjk)^3, (ijkj)^2, (dhjidg)^3\rangle$.

\item[\rm(f)] $z = a^2$ is a $2$-central involution of $E$ with
centralizer $D = C_{E}(z) = \langle x, y \rangle$ of order
$2^{21}\cdot3^3\cdot5$ where $x = a(agik)^3$ and $y = d(cgihj)^4$
have orders $4$ and $6$, respectively and $z = (xy^3)^{12}$.
Furthermore, $E = \langle x, y, e\rangle$, where $e = g$ has order
$4$.

\item[\rm(g)] $E = \langle x,y,e\rangle$ has $73$ conjugacy
classes. A system of their representatives is given in Table
\ref{Fi_24 cc E}.

\item[\rm(h)] Table \ref{Fi_24 ct_E} is the character table of
$E$.

\item[\rm(i)] $V_2$ is the unique maximal elementary abelian
normal subgroup of each Sylow $2$-subgroup $S$ of the extension
group $E$.

\item[\rm(j)] $C_E(V_2) = V_2$.

\item[\rm(k)] $N_1=\langle a,b,c,d,t,g,h,i,j,V_2 \rangle$ is a non split extension of
$M_{23}$ by $V_2$.

\item[\rm(l)] $u = (agt)^5$ is an involution of $N_1$ generating the center $Z(R_2)$ of $C_{R_1}(u)$.

\item[\rm(m)] $E = \langle N_1, C_E(u)\rangle$.
\end{enumerate}
\end{lemma}

\begin{proof}
(a) The $2$ irreducible $F\M_{24}$-modules $V_i$, $i =1,2$, occur
as composition factors with multiplicity $1$ in the permutation
module $(1_{\M_{23}})^{\M_{24}}$ and can easily be constructed
using the faithful permutation representation of $\M_{24}$ stated
in (a) and the Meataxe algorithm implemented in MAGMA. The
corresponding matrices of the generators of $\M_{24}$ with respect
to the first irreducible representation of $\M_{24}$ are stated in
(a). The second irreducible representation $V_2$ of $\M_{24}$ is
dual to $V_1$ and so it is defined by the equations given in (b).

(c) The cohomological dimensions $d_i = dim_F[H^2(\M_{24},V_i)]$,
$i =1,2$, have been calculated by means of MAGMA using Holt's
Algorithm 7.4.5 of \cite{michler}, the presentation of $\M_{24}$
of Definition 8.2.1 of \cite{michler} and all the data stated in
(a) and (b). It follows that $d_1 = 0$ and $d_2 = 1$.

(d) The presentation of $E$ has been obtained by means of step 3
of Holt's Algorithm 7.4.5 of \cite{michler} and MAGMA.

(e) Using a stand-alone program due Paul Young we found a faithful
permutation representation $pE$ of $E$ of degree $24288$ with
stabilizer $\langle (i^{-1}j^{-1}(b g)^2,
(i^{-1}bth^{-1})^4\rangle $. Applying the MAGMA command
$\verb"DegreeReduction(pE)"$ we obtained the faithful permutation
representation $PE$ of degree $1518$. The given generators of its
stabilizer $T$ were obtained by means of the first author's
program\\
$\verb"GetShortGens(PE,BasicStabilizer(PE,2))"$.

(f) Using MAGMA and the faithful permutation representation $PE$
of $E$ the reader easily verifies that the centralizer $C_E(z)$ of
$z = a^2$ has order $2^{21}\cdot3^3\cdot5$. Hence $z$ is a
$2$-central involution of $E$ by (d). The words of the generators
$x$, $y$ of $D$ were calculated by means of $PE$, MAGMA and the
first author's program $\verb"GetShortGens(PE,PD)"$. Another check
with MAGMA and $PE$ verifies that $E = \langle D,g\rangle$.

(g) Using Kratzer's Algorithm 5.3.18 of \cite{michler}, the
faithful permutation representation $PE$ and MAGMA we observed
that $E$ has $73$ conjugacy classes. Their representatives are
given in Table \ref{Fi_24 cc E}.

(h) The character table of $E$ was automatically computed by MAGMA
using $PE$.

The remaining $4$ statements can be checked with MAGMA and the faithful permutation representation $PE$.
\end{proof}

\section{The $2$-fold cover of the automorphism group $\rm Aut(\Fi_{22})$ }

Applying Algorithm 2.5 of \cite{michler1} to an extension group isomorphic to the subgroup $N_1$ of $E = \langle x,y,e\rangle$ described in Lemma \ref{l. M24-extensions}(k) H. Kim realized Fischer's
second sporadic simple group $\Fi_{23}$ as an irreducible subgroup
$G_1$ of $\GL_{782}(17)$ in his senior thesis \cite{kim1}. He showed that the centralizer $H_1 = C_{G_1}(u)$ of a $2$-central involution $u$ of $G_1$ is isomorphic to the $2$-fold cover $2Fi_{22}$ of Fischer's smallest sporadic simple group $\Fi_{22}$. Furthermore, he constructed a faithful permutation representation $PG_1$ of $G_1$ of degree $31671$. In
this section we use these results to construct the $2$-fold cover $A_1$ of the automorphism group $Aut(H_1)$ and show that $H_1$ is its commutator subgroup. Thus we obtain an amalgam $A_1
\leftarrow H_1 \rightarrow G_1$ such that $A_1$ and $D_1 = C_E(z_1)$ have isomorphic Sylow $2$-subgroups where $z_1$ is the involution $e^2$ of $E$. Using the character tables of the $3$ groups of the amalgam we also show that it has $8$ compatible
pairs of semi-simple characters of degree $8671$.

\begin{lemma}\label{l. aut2Fi22}
Let $A_1$ be the $2$-fold cover of the automorphism group
$Aut(\Fi_{22})$ of Fischer's simple group $\Fi_{22}$ and let $H_1
= A_1'$ be its derived subgroup. Let $E = \langle x,y,e\rangle$ be the non split extension of $\M_{24}$ by its simple $GF(2)$-module constructed in Lemma \ref{l. M24-extensions}. Then the following assertions
hold:

\begin{enumerate}
\item[\rm(a)] $H_1 = \langle a,b,c,d,e,f,g,h,i,z \rangle$ has the
following set $\mathcal R(H_1)$ of defining relations:
\begin{eqnarray*}
&&a^2 = b^2 = c^2 = d^2 =  e^2 = f^2 = g^2 = h^2 = i^2 = 1,\\
&&(ab)^3 = 1, (bc)^3 = z, (c d)^3 = (d e)^3 = 1, (e f)^3 = (f g)^3 = z, \\
&&(a c)^2 = (a d)^2 = (a e)^2 = (a f)^2 = (a g)^2 = (a h)^2 = (a i)^2 = 1,\\
&&(b d)^2 = (b e)^2 = (b f)^2 = (b g)^2 = (b h)^2 = (b i)^2 = 1,\\
&&(c e)^2 = (c f)^2 = (c g)^2 = (c h)^2 = (c i)^2 = 1,\\
&&(d f)^2 = (d g)^2 = (e g)^2 = (e h)^2 = (e i)^2 = 1,\\
&& (d h)^3 = (h i)^3 = (d i)^2 = (f h)^2 = (f i)^2 = (g h)^2 = (g i)^2 = 1,\\
&&(d c b d e f d h i)^{10} = (a b c d e f h)^9 = (b c d e f g h)^9 = 1,\\
&& z^2 = (z,a) = (z,b) = (z,c)= (z,d)=(z,e)=(z,f)= 1,\\
&& (z,g)=(z,h)=(z,i) = 1.
\end{eqnarray*}

\item[\rm(b)] $A_1 = \langle H_1, t\rangle$ has a set $\mathcal
R(A_1)$ of defining relations consisting of $\mathcal R(H_1)$ and
the following relations:

\begin{eqnarray*}
&&t^2 = 1,\quad (z,t) = 1, \quad a^tg = 1, \quad b^tf = c^te = (dt)^2 = (ht)^2 = (it)^2 = z.\\
\end{eqnarray*}

\item[\rm(c)] $A_1 = 2Aut(\Fi_{22})$ has a faithful permutation
representation $PA_1$ of degree $56320$ with stabilizer $\langle
bz, c, d, e, fz, g, h, i \rangle$.

\item[\rm(d)] A system of representatives $a_i$ of the $150$
conjugacy classes of $A_1$ and the corresponding centralizers
orders $|C_A(a_i)|$ are given in Table \ref{Fi_24 cc A_1}.

\item[\rm(e)] The character table of $A_1$ is given in Table
\ref{Fi_24 ct A_1}.

\item[\rm(f)] The group $A_1$ and the centralizer $C_E(z_1)$ of the involution $z_1 = e^2$ of $E$ have isomorphic Sylow $2$-subgroups of order $2^{19}$.
\end{enumerate}
\end{lemma}

\begin{proof}
(a) The given presentation of $H_1$ is a restatement of
Proposition 6.2.3 of \cite{michler2} due to H. Kim, see
\cite{kim1}.

(b) By that result we also know that $H_1$ has a faithful
permutation representation $PH_1$ of degree $28160$ with
stabilizer $U = \langle bz, c, d, e, fz, g, h, i \rangle$. Using
it and the MAGMA command $\verb"AutomorphismGroup(H_1)"$ we see
that $|Aut(H_1)| = |H_1|$. As $H_1/\langle z \rangle \cong
\Fi_{22}$ has the same presentation as $\Fi_{22}$ given in
\cite{praeger}, p. 110 we can quote the presentation of
$Aut(\Fi_{22})$ given in \cite{praeger}, p. 111, where $z$ is
replaced by $1$. Now (b) follows from (a) and $2^6$ MAGMA
calculations with $PH_1$ checking whether $1$ or $z$ has to be on
the right hand side of the six new relations stated in (b) different from $t^2 = 1$ and $[z,t] = 1$. It
follows that there is exactly one solution.

(c) This statement has been verified by means of the MAGMA
command\\ $\verb"CosetAction(A_1,U)"$.

(d) The system of representatives of the conjugacy classes of
$A_1$ has been calculated by means of the permutation
representations $PA_1$ of $A_1 = 2Aut(H_1)$, MAGMA and Kratzer's
Algorithm 5.3.18 of \cite{michler}.

(e) The character table of $A_1$ has been calculated by means of
$PA_1$ and MAGMA.

(f) Let $PE$ be the faithful permutation representation of $E$ constructed in Lemma \ref{l. M24-extensions}.
Let $C = C_E(u)$ for the involution $u = e^2$ of $E = \langle x,y,e\rangle$. Now (f) can be verified by using the permutation representations $PA_1$ and $PE$ together with the Cannon-Holt isomorphism test implemented in MAGMA.
\end{proof}

\begin{lemma}\label{l. gensFi24}
Keep the notation of Lemma \ref{l. aut2Fi22}. Let $G_1  = \langle
x,y,q,w\rangle \cong \Fi_{23}$ be the simple subgroup of
$\GL_{782}(17)$ of order $2^{18}\cdot 3^{13}\cdot
5^2\cdot7\cdot11\cdot13\cdot17\cdot23$ with centralizer $H_1 =
\langle x,y,q\rangle = C_{G_1}(z)$ of the $2$-central involution
$z = (xy^2)^7$ constructed in \cite{kim1}. Let
$A_1 = 2Aut(H_1)$. Then the following assertions hold:
\begin{enumerate}
\item[\rm(a)] $H_1 = \langle a,b,c,d,e,f,g,h,i,z \rangle$ where
\begin{eqnarray*}
&&a = (xyx)^7,\quad b = [(qy)^2qy^3q^2y^3qy]^7,\quad c = (y^2xyxy^3)^5,\\
&&d = (qyq^2yqyqyqy^2q^2)^{15},\quad e = (yxy^5x)^5,\quad f = (yqyq^2yq^2y^2qy^4q^2)^5,\\
&&g = (xy^2xy^3x)^7,\quad h = (y^5xyx)^5,\quad i = (q^2y^2qyq^2)^7.\\
\end{eqnarray*}

\item[\rm(b)] $H_1 = \langle a,b,c,d,e,f,g,h,i,z \rangle$
satisfies the set $\mathcal R(H_1)$ of defining relations stated
in Lemma \ref{l. aut2Fi22}(a). Its character table is Table 6.5.2
of \cite{michler2}.

\item[\rm(c)] The character table of $G_1$ is stated in the Atlas
\cite{atlas}, its pp. 178 - 179.

\item[\rm(d)] The amalgam $A_1 \leftarrow H_1 \rightarrow G_1$ has
Goldschmidt index $1$.

\item[\rm(e)] The amalgam $A_1 \leftarrow H_1 \rightarrow G_1$ has
eight compatible pairs $$(\chi,\tau)\in mf
\mbox{char}_{\mathbb{C}}(A_1) \times mf
\mbox{char}_{\mathbb{C}}(G_1)$$ of degree $8671$. All have the
same restriction $$\delta_{2} + \delta_{6} + \delta_{7} +
\delta_{8} + \delta_{9} \in mf \mbox{char}_{\mathbb{C}}(H_1).$$
They are:

\begin{enumerate}
\item[\rm(1)] $(\chi_{3} + \chi_{11}  + \chi_{13} + \chi_{\bf 17},
\quad \tau_{\bf 3} + \tau_{\bf 4})$,

\item[\rm(2)] $(\chi_{3} + \chi_{11}  + \chi_{14} + \chi_{\bf 17},
\quad \tau_{\bf 3} + \tau_{\bf 4})$,

\item[\rm(3)] $(\chi_{3} + \chi_{12}  + \chi_{13} + \chi_{\bf 17},
\quad \tau_{\bf 3} + \tau_{\bf 4})$,

\item[\rm(4)] $(\chi_{3} + \chi_{12}  + \chi_{14} + \chi_{\bf 17},
\quad \tau_{\bf 3} + \tau_{\bf 4})$,

\item[\rm(5)] $(\chi_{4} + \chi_{11}  + \chi_{13} + \chi_{\bf 17},
\quad \tau_{\bf 3} + \tau_{\bf 4})$,

\item[\rm(6)] $(\chi_{4} + \chi_{11}  + \chi_{14} + \chi_{\bf 17},
\quad \tau_{\bf 3} + \tau_{\bf 4})$,

\item[\rm(7)] $(\chi_{4} + \chi_{12} + \chi_{13} + \chi_{\bf 17},
\quad \tau_{\bf 3} + \tau_{\bf 4})$,

\item[\rm(8)] $(\chi_{4} + \chi_{12} + \chi_{14} + \chi_{\bf 17},
\quad \tau_{\bf 3} + \tau_{\bf 4})$.

\end{enumerate}
\end{enumerate}
\end{lemma}

\begin{proof}
By Kim's Theorem 6.3.1 of \cite{michler2} the simple matrix
subgroup $G_1 \cong \Fi_{23}$ of $\GL_{782}(17)$ has a faithful
permutation representation $PG_1$ of degree $31671$ with
stabilizer $H_1$. It is used throughout this proof.

(a) The words of the new generators $a$, $b$, etc. of $H_1$ in
terms of the given generators $x$, $y$ and $q$ of $H_1$ are quoted
from Kim's Proposition 6.2.3 of \cite{michler2}.

(b) Using the faithful permutation representation $PG_1$ and MAGMA
it has been checked that the new generators $a$, $b$ etc. of $H_1$
given in statement (a) satisfy all the relations of $\mathcal
R(H_1)$ of Lemma \ref{l. aut2Fi22}(a).

(c) This assertion is a restatement of Theorem 6.3.1 of
\cite{michler2}.

(d) Kratzer's Algorithm 7.1.10 of \cite{michler} could not be
applied to calculate the Goldschmidt index. When trying to
calculate $Aut(G_1)$ using the Cannon-Holt Algorithm of
\cite{cannon} MAGMA answered: ``Sorry, the top factor of order
4089470473293004800 is not currently stored". However, using the
faithful permutation representation $PA_1$ of $A_1$ stated in
Lemma \ref{l. aut2Fi22}(c) MAGMA established that the outer
automorphism groups $Out(H_1)$ and $Out(A_1)$ of $H_1$ and $A_1$
are both cyclic of order $2$. Hence the Goldschmidt index of the
amalgam $A_1 \leftarrow H_1 \rightarrow G_1$ is $1$ by Step 3 of
Algorithm 7.1.10 of \cite{michler}.

(e)  The eight compatible pairs of degree $8671$ of the amalgam
$A_1 \leftarrow H_1 \rightarrow G_1$ were determined by means of Kratzer's Algorithm 7.3.10
of \cite{michler} and MAGMA.
\end{proof}

\section{A semi-simple representation of $\Fi_{23}$ over $GF(13)$}

In this section H. Kim's results of his senior thesis \cite{kim1}
are used for the construction of two irreducible representations
of degrees $3588$ and $5083$ of $G_1 \cong \Fi_{23}$ over the
prime field $GF(13)$. They correspond to the $2$ irreducible
characters $\tau_{\bf 3}$ and $\tau_{\bf 4}$ of $13$-defect zero
of $G_1$ occurring in the $8$ compatible pairs constructed in
Lemma \ref{l. gensFi24}(e). For the construction of these fairly
large representations we first determine generators of two large
subgroups $mH$ and $mE$ of $G_1$ and their intersection $mD$. We
also calculate the character tables of these $3$ subgroups of
$G_1$.

\begin{lemma}\label{l. reps13Fi23}
Keep the notation of Lemmas \ref{l. aut2Fi22} and \ref{l.
gensFi24}. Let $G_1 = \Fi_{23} = \langle x,y,q,w\rangle$ be the
simple subgroup of $\GL_{782}(17)$ of order $2^{18}\cdot
3^{13}\cdot 5^2\cdot7\cdot11\cdot13\cdot17\cdot23$ with faithful
permutation representation $PG_1$ of degree $31671$ and stabilizer
$H_1 = \langle x,y,q\rangle$ constructed in Kim's Theorem 6.3.1 of
\cite{michler2}. Let $r = s_1 = (yqy^2q)^7$, $s_2 = (yqyqy)^7$,
$s_3 = (yeyeq)^{21}$, $s_4 = (qyqey^2)^7$, $u =
s_4^2(s_1s_2s_4)^3$, $s = qyqey^2$ and $v = (s_1s_3s_1s_4)^2
(s_1s_3s_1s_3s_1s_4)^3s_3s_4s_3s_4^2s_3$. Then the following
assertions hold:

\begin{enumerate}
\item[\rm(a)] $mH = \langle s_1, s_2, s_3, s_4 \rangle = \langle
r,u,v \rangle \cong O_8^+(3):S_3$.

\item[\rm(b)] The character table of $mH$ is given in Table
\ref{Fi_24 ct mH}.

\item[\rm(c)] $mE = \langle u,v,s \rangle \cong O_7(3)\times S_3$.

\item[\rm(d)] The character table of $mE$ is Table \ref{Fi_24 ct
mE}.

\item[\rm(e)] $mD = mH \cap mE = \langle u,v \rangle \cong G_2(3)
\times S_3$.

\item[\rm(f)] $G_1 = \langle u,v,r,s \rangle$ and the original
generators $q$, $y$, $w$ and $x$ of $G_1$ are equal to the
following words in its generators $u$, $v$, $r$ and $s$:
\begin{eqnarray*}
&&q = [(us^2vsr)^9[(svs^3rsr)^{12}(vrsrsvu)^{12}]^2]^{24},\\
&&w = [(w_2w_4w_2w_4w_3w_2w_4)^3(w_1w_2w_3w_1w_4w_3w_1w_4)^7]^3,\\
&&y = (m_1m_2m_1m_2m_3)^7[(m_1m_2)^5(n_3n_1n_2n_3n_1n_3^2n_1)^5]^4,\\
&&x = [(yq^2yqyq^2)^{11}(q^2y^2qyqy)^{11}(qy^2qyqyqyq)^4]^{12},\quad \mbox{where}\\
&&w_1 = (vuvs^2)^9,\quad w_2 = (v^2uruvr)^{18},\quad w_3 = (v^2s^2urv)^{11},\\
&&w_4 = (uvrv^4r)^{10},\quad n_1 = m_1m_2m_3m_2^2m_1,\quad n_2 = m_2m_3m_1m_3m_2^3,\\
&&n_3 = (m_3m_1m_2m_3m_2m_1m_2)^2,\quad m_1 = (t_1t_3t_1t_3t_1^2t_3t_1)^2,\\
&&m_2 = (t_1t_3t_1t_3t_1^2t_2t_3t_2)^2,\quad m_3 = (t_2t_1t_3t_2t_3t_1t_2t_3t_2)^9,\\
&&t_1 = (srsrsrs)^5, \quad t_2 = (s^2rs^4)^{11}, \quad \mbox{and} \quad t_3 = (rs^7)^3.\\
\end{eqnarray*}

\item[\rm(g)] The restrictions of the irreducible character
$\tau_{\bf 3}$ of degree $3588$ of $G_1$ to $mH$ and $mE$ are
$\pi_8 + \pi_{16} \in mf \mbox{char}_{\mathbb{C}}(mH)$ and
$\psi_{4} + \psi_{10} + \psi_{19} + \psi_{36} + \psi_{61} \in mf
\mbox{char}_{\mathbb{C}}(mE)$, respectively.

\item[\rm(h)] The restrictions of the irreducible character
$\tau_{\bf 4}$ of degree $5083$ of $G_1$ to $mH$ and $mE$ are
$\pi_{12} + \pi_{15} \in mf \mbox{char}_{\mathbb{C}}(mH)$ and
$\psi_{14} + \psi_{22} + \psi_{45} + \psi_{74}) \in mf
\mbox{char}_{\mathbb{C}}(mE)$, respectively.

\item[\rm(i)] The irreducible characters $\pi_9$, $\pi_{11}$ and
$\pi_{15}$ of $mH$ are constituents of the permutation characters
$1_{mH_9}^{mH}$, $1_{mH_{11}}^{mH}$ and $1_{mH_{15}}^{mH}$ of the
subgroups
\begin{eqnarray*}
&&mH_9 = \langle (vru^2)^{13}, (rvur)^2, (uvrv^2)^4\rangle,\\
&&mH_{11} = \langle (uvurv)^4, (v^5ur)^9, (ururv^2u^2)^{13}\rangle \quad \mbox{and}\\
&&mH_{15} = (uvuru)^{13}, (uv^2uv)^2, (u^5v)^3\rangle\\
\end{eqnarray*}
of $mH$ with indices $3240$, $72800$ and $2274480$, respectively.

\item[\rm(j)] The linear character $\pi_2$ of $mH$ has values $-1$
and $1$ at $v$ and $u$, $r$, respectively. Furthermore,
$\pi_{8} = \pi_{2} \otimes \pi_{9},\quad \pi_{12} = \pi_{2} \otimes \pi_{11},\quad \mbox{and} \quad \pi_{16} = \pi_{2} \otimes \pi_{15}$.\\

\item[\rm(k)] The irreducible characters $\psi_5$ and $\psi_{10}$
of $mE$ are constituents of the permutation character
$1_{mE_1}^{mE}$ of the subgroup $mE_1 = \langle (s^2vs)^9,
(svsv^2s^2)^2 \rangle$ of index $2106$.

\item[\rm(l)] The irreducible characters $\psi_{14}$, $\psi_{22}$,
$\psi_{36}$, $\psi_{45}$, $\psi_{61}$ and $\psi_{74}$ of $mE$ are
constituents of the permutation characters $1_{mE_{14}}^{mE}$,
$1_{mE_{22}}^{mE}$, $1_{mE_{36}}^{mE}$, $1_{mE_{45}}^{mE}$,
$1_{mE_{61}}^{mE}$ and $1_{mE_{74}}^{mE}$ of the subgroups
\begin{eqnarray*}
&&mE_{14} = \langle (s^2us)^6, (s u^2 s^2)^4, (s u s u^2)^2, (us^3 u s)^7 \rangle,\\
&&mE_{22} = \langle (us^4u^2)^6, u s u^2 s u s u s \rangle,\\
&&mE_{36} = \langle (s^2vs)^9, (v^2 s v s v^2)^5, (vsv^2svs)^2 \rangle,\\
&&mE_{45} = \langle (s^2v^2)^7, (v^5 s)^3, (v s v^2 s v^3 s v)^6 \rangle,\\
&&mE_{61} = \langle (vs^2)^7, (v^4 s^2v^2)^{12},(v^2svsv^4)^{10}, (v s^2 v^3 s^3)^{30} \rangle \quad \mbox {and}\\
&&mE_{74} = \langle (s^2v)^7, (v^2 s^2)^{21}, (s^2 v^2)^{21}, (s v s v^2 s v)^3 \rangle\\
\end{eqnarray*}
of $mE$ with indices $702$, $2160$, $19656$, $7280$, $85293$ and
$29484$, respectively.

\item[\rm(m)] The linear character $\psi_2$ of $mE$ has values
$-1$ and $1$ at $v$ and $u$, $s$, respectively. Furthermore,
$\psi_{4} = \psi_{2} \otimes \psi_{5},\quad \mbox{and} \quad
\psi_{19} = \psi_{2} \otimes \psi_{22}$.

\item[\rm(n)] Both $r$ and $f = (u^3vsv)^9$ are involutions of
$G_1$ such that $(r, f) = 1$, $rf \notin mH$ and $rf \notin mE$.
\end{enumerate}
\end{lemma}

\begin{proof}
(a) The subgroup $mH$ of $G_1$ has been constructed by means of
the faithful  permutation representation $PG_1$ of $G_1$ of degree
$31671$ and the MAGMA command\\
$\verb"LowIndexSubgroups(PG_1, 137632)"$. The four generators
$s_i$ of $mH$, $1 \le i \le 4$, were calculated with Kim's program
$\verb"GetShortGens(PG_1, mH)"$. Another application of MAGMA
determined the composition factors of $mH$.

(b) The character table of $mH$ was calculated by MAGMA using
$PG_1$.

(c) and (e) By Table 6.5.4 of \cite{michler2} $|C_{G_1}((qw)^4)| =
2^9\cdot3^{10}\cdot5\cdot7\cdot13$. Let $mX = N_{G_1}(\langle
(qw)^4)\rangle$. Using $PG_1$ and MAGMA we searched for an element
$x \in mX$ of order $3$ such that $|N_{mH}(\langle x \rangle)| =
2^7\cdot3^7\cdot7\cdot13$. MAGMA found such an element and stated
that $mD = N_{mH}(\langle x \rangle) = \langle u, v\rangle$.
Furthermore, $mE = N_{G_1}(\langle x \rangle) = \langle mD,
s\rangle$ where $u$, $v$ and $s$ are defined in the statement of
this lemma. The composition factors of $mD$ and $mE$ have been
determined by means of MAGMA.

(d) The character table of $mE$ was calculated by means of $PG_1$
and MAGMA.

(f) Using the faithful permutation representation $PG_1$ of $G_1$
and MAGMA one verifies that $G_1 = \langle mH, mE\rangle$. Hence
$G_1 = \langle u,v,r,s \rangle$ by (a) and (c). The words for $q$,
$y$, $w$ and $x$ can easily be checked computationally.

(g) $G_1$ has a unique character $\tau_{\bf 3}$ of degree $3588$
by the character table of $G_1 \cong \Fi_{23}$, see \cite{atlas},
p. 178. Its restrictions to $mH$ and $mE$ given in the statements
have been determined by means of $PG_1$, the character tables of
the subgroups $mH$ and $mE$ of $G_1$, and MAGMA.

(h) This assertion is proved as (g).

(i) Using the MAGMA command $\verb"LowIndexSubgroups(mH, k)"$ we
searched for conjugacy classes of subgroups $H_k$ of index $|mH :
mH_k| = m_k$ such that $\pi_k$ is an irreducible constituent of
the permutation character $1_{mH_k}^{mH}$ for $k \in \{9, 11,
15\}$. Thus we found $3$ subgroups $mH_k$ of respective indices
$m_{9} = 3240$, $m_{11} = 72800$ and $m_{15} = 2274480$. Their
given generators have been obtained by means of Kim's program
$\verb"GetShortGens(mH, mH_k)"$.

(j) $\pi_2$ is the unique non trivial linear character of $mH$ by
its character table. The character equations of the statement are
easily verified by means of Table \ref{Fi_24 ct mH}.

The statements (k) and (l) are proved similarly as (i).

(m) $mE$ has a unique non trivial linear character $\psi_2$, see
Table \ref{Fi_24 ct mE}. The character equations of the statement
are easily verified by means of Table \ref{Fi_24 ct mH}.

(n) Using $PG_1$ and MAGMA we checked that $r$ and $f$ are
commuting involutions such that $rf \notin mH$ and $rf \notin mE$.
\end{proof}

\begin{proposition}\label{prop. rep3588Fi23}
Keep the notation of Lemmas \ref{l. aut2Fi22}, \ref{l. gensFi24}
and \ref{l. reps13Fi23}. Let $PG_1$ be the faithful permutation
representation of the simple group $G_1 = \langle x,y,q,w\rangle =
\langle u,v,r,s \rangle$ of degree $31671$ with stabilizer $H_1 =
\langle x,y,q\rangle$. Let $mH = \langle u,v,r \rangle$, $mD =
\langle u,v \rangle$ and $mE = \langle u,v,s \rangle$. Let $F^{*}$
be the multiplicative group of the prime field $F = \GF(13)$. Let
$Y = \GL_{3588}(13)$.

Let $\mathfrak V$ and $\mathfrak W$ be the up to isomorphism
uniquely determined faithful semi-simple $3588$-dimensional
modules of $mH$ and $mE$ over $F$ corresponding to the
restrictions ${\tau_{\bf 3}}_{|mH}$ and ${\tau_{\bf 3}}_{|mE}$ of
the irreducible character $\tau_3$ of $G_1$, respectively.

Let $\kappa_\mathfrak V : mH \rightarrow \GL_{3588}(13)$ and
$\kappa_\mathfrak W : mE \rightarrow \GL_{3588}(13)$ be the
representations of $mH$ and $mE$ afforded by the modules
$\mathfrak V$ and $\mathfrak W$, respectively.

Let $\mathfrak r = \kappa_\mathfrak V(r)$, $\mathfrak u =
\kappa_\mathfrak V(u)$, $\mathfrak v = \kappa_\mathfrak V(v)$ in $
\kappa_\mathfrak V(mH) \le \GL_{3588}(13)$.

Then $\mathfrak V_{|mD} \cong \mathfrak W_{|mD}$, and there is a
transformation matrix $\mathcal T_1 \in \GL_{3588}(13)$ such that
$$
\mathfrak u = \mathcal T_1^{-1} \kappa_\mathfrak W (u) \mathcal
T_1, \mathfrak v = \mathcal T_1^{-1} \kappa_\mathfrak W(v)
\mathcal T_1.
$$

Let $\mathfrak {mD} = \langle \mathfrak u, \mathfrak v \rangle$,
$\mathfrak {mH} = \langle \mathfrak u, \mathfrak v, \mathfrak r
\rangle$. Let $\mathcal D = C_Y(\mathfrak {mD})$ and $\mathcal H =
C_Y(\mathfrak {mH})$. Let $\mathfrak s_1 = \mathcal T_1^{-1}
\kappa_\mathfrak W (s) \mathcal T_1$. Let $\mathfrak {mE} =
\langle \mathfrak {mD}, \mathfrak s_1 \rangle$ and $\mathcal E =
C_Y(\mathfrak {mE})$. Then the following statements hold:

\begin{enumerate}
\item[\rm(a)] There is an isomorphism $$\alpha: \mathcal D
\rightarrow \mathcal D_1 = \GL_2(13)\times \GL_2(13) \times
{F^{*}}^5 \le \GL_9(13).$$

\item[\rm(b)] $\mathcal H_1 = \alpha(\mathcal H)$ is generated by
the two blocked diagonal matrices\\
\begin{align*}
a_1 = diag( \left( \begin{smallmatrix} 2 & 0 \\  0 & 1
\end{smallmatrix} \right), \left( \begin{smallmatrix} 2 & 0 \\  0
& 1
\end{smallmatrix} \right), 2,1,2,1,1) \quad \mbox{and} \quad
a_2 = diag( \left( \begin{smallmatrix} 1 & 0 \\  0 & 2
\end{smallmatrix} \right), \left( \begin{smallmatrix} 1 & 0 \\  0
& 2
\end{smallmatrix} \right), 1,2,1,2,2),
\end{align*}

\item[\rm(c)] $\mathcal E_1 = \alpha(\mathcal E)$ is generated by
the five blocked diagonal matrices
\begin{align*}
b_1 =diag( \left( \begin{smallmatrix} 2 & 0 \\  0 & 1
\end{smallmatrix} \right), \left( \begin{smallmatrix} 1 & 0 \\  0
& 1
\end{smallmatrix} \right), 1,1,1,1,1), \quad
b_2 = diag( \left(
\begin{smallmatrix} 1 & 0 \\  0 & 1 \end{smallmatrix} \right),
\left( \begin{smallmatrix} 2 & 0 \\  0 & 1 \end{smallmatrix}
\right),
1,1,1,1,1), \\
b_3 = diag( \left( \begin{smallmatrix} 1 & 0 \\  0 & 2
\end{smallmatrix} \right), \left( \begin{smallmatrix} 1 & 0 \\  0
& 1
\end{smallmatrix} \right), 1,2,1,1,1), \quad
b_4 = diag( \left(
\begin{smallmatrix} 1 & 0 \\  0 & 1 \end{smallmatrix} \right),
\left( \begin{smallmatrix} 1 & 0 \\  0 & 1 \end{smallmatrix}
\right),
2,1,1,2,1), \\
\mbox{and} \quad
b_5 = diag( \left( \begin{smallmatrix} 1 & 0 \\
0 & 1
\end{smallmatrix} \right), \left( \begin{smallmatrix} 1 & 0 \\  0
& 2 \end{smallmatrix} \right), 1,1,2,1,2).
\end{align*}

\item[\rm(d)] $\mathcal D$ has $2184^2 \times 12 = 57238272$
$\mathcal H$-$\mathcal E$ double cosets.

\item[\rm(e)] The free product $mH*_{mD}mE$ of $mH$ and $mE$ with
amalgamated subgroup $mD$ has an irreducible $3588$-dimensional
representation over $F$ which induces an irreducible
representation of $G_1$. It corresponds to the $\mathcal
H$-$\mathcal E$ double coset representative

$$\mathcal F = diag(w,z,1^{182},1^{182},1^{364},1^{728},1^{1664}) \in \GL_{3588}(13), \quad
{where}$$

{\renewcommand{\arraystretch}{0.5} \scriptsize
$$
w = \left( \begin{array}{*{2}{c@{\,}}c}
 1^{78} & 4^{78}\\
12^{78} & 2^{78}
\end{array} \right),\quad
z = \left( \begin{array}{*{2}{c@{\,}}c}
1^{156} &  6^{156}\\
5^{156} & 10^{156}
\end{array} \right),
$$
} and $a^n$ denotes a diagonal $n \times n$ matrix with unique
diagonal non zero entry $a \in GF(13)$.

Let $\mathfrak s = \mathcal F^{-1} \mathfrak s_1 \mathcal F$ and
$\mathfrak G_1 = \langle \mathfrak u, \mathfrak v, \mathfrak r,
\mathfrak s \rangle$. Inserting these four generating matrices of
$\mathfrak G_1$ into the formulas of Lemma \ref{l. reps13Fi23}(f)
one obtains the matrices $\mathfrak x_{3588}$, $\mathfrak
y_{3588}$, $\mathfrak q_{3588}$, and $\mathfrak w_{3588}$ of the
original generators $x$, $y$, $q$ and $w$ of $G_1 = \Fi_{23}$ as
words in the generators $u$, $v$, $r$ and $s$. The matrices
$\mathfrak x_{3588}$, $\mathfrak y_{3588}$, $\mathfrak q_{3588}$,
and $\mathfrak w_{3588}$ can be downloaded from the first author's
website\\
$\verb"http://www.math.yale.edu/~hk47/Fi24/index.html"$.
\end{enumerate}
\end{proposition}

\begin{proof}
Let $\mathfrak V$ be the up to isomorphism uniquely determined
faithful semi-simple $3588$-dimensional module of $mH$ over $F =
GF(13)$ corresponding to the restriction ${\tau_{\bf 3}}_{|mH}$.
By Lemma \ref{l. reps13Fi23}(g) ${\tau_{\bf 3}}_{|mH} = \pi_8 +
\pi_{16}$. Lemma \ref{l. reps13Fi23}(j) states that $\pi_{8} =
\pi_{2} \otimes \pi_{9}$.

The irreducible characters $\pi_9$ and $\pi_{15}$ of $mH$ are
constituents of the permutation characters $1_{mH_9}^{mH}$ and
$1_{mH_{15}}^{mH}$ of the respective subgroups $mH_9$ and
$mH_{15}$ of $mH$ determined in Lemma \ref{l. reps13Fi23}(i).
Using a stand-alone program of the first author which is based on
Algorithm 5.7.1 of \cite{michler} we calculated the primitive
idempotents of the endomorphism rings of these permutation
modules. Thus we obtained the corresponding irreducible
representations $M(\pi_9)$ and $M(\pi_{15})$ of the respective
dimensions $780$ and $2808$ over $F$. The irreducible
$FmH$-modules $M(\pi_8)$ and $M(\pi_{16})$ are the tensor products
of $M(\pi_9)$ and $M(\pi_{15})$ with the linear character $\pi_2$
of $mH$ over $F$. Thus $\mathfrak V = M(\pi_8) \oplus
M(\pi_{16})$.

Let $\mathfrak W$ be the up to isomorphism uniquely determined
faithful semi-simple $3588$-dimensional module of $mE$ over $F$
corresponding to the restriction ${\tau_{\bf 3}}_{|mE}$. By Lemma
\ref{l. reps13Fi23}(g) ${\tau_{\bf 3}}_{|mE} = \psi_4 + \psi_{10}
+ \psi_{19} + \psi_{36} + \psi_{61}$. Lemma \ref{l. reps13Fi23}(m)
states that $\psi_{4} = \psi_{2} \otimes \pi_{5}$ and $\psi_{19} =
\psi_{2} \otimes \pi_{22}$.

The irreducible characters $\psi_5$ and $\psi_{10}$ of $mE$ are
constituents of the permutation character $1_{mE_1}^{mE}$ by
\ref{l. reps13Fi23}(k). The irreducible characters $\psi_{22}$,
$\psi_{36}$ and $\psi_{61}$ of $mE$ are constituents of the
permutation characters $1_{mE_{22}}^{mE}$, $1_{mE_{36}}^{mE}$ and
$1_{mE_{61}}^{mE}$, respectively, see Lemma \ref{l.
reps13Fi23}(m). Using a standalone program of the first author
which is based on Algorithm 5.7.1 of \cite{michler} we calculated
the primitive idempotents of the endomorphism rings of these four
permutation modules. Thus we obtained the corresponding
irreducible representations $N(\psi_5)$, $N(\psi_{10})$,
$N(\psi_{22})$, $N(\psi_{36})$ and $N(\psi_{61})$ of the
respective dimensions $78$, $156$, $260$, $910$ and $2184$ over
$F$. The irreducible $FmE$-modules $N(\psi_4)$ and $N(\psi_{19})$
are the tensor products of $N(\psi_5)$ and $N(\psi_{22})$ with the
linear character $\psi_2$ of $mE$ over $F$. Thus $$\mathfrak W =
N(\psi_4) \oplus N(\psi_{10}) \oplus N(\psi_{19}) \oplus
N(\psi_{36}) \oplus  N(\psi_{61}).$$

Fixing a basis in each irreducible constituent $M\pi_k$ of
$\mathfrak V$ we get a basis $\mathcal B_V$ of $\mathfrak V$. It
induces a representation $\kappa_\mathfrak V : mH \rightarrow
\GL_{3588}(13)$ of $mH$. Let $\mathfrak r = \kappa_\mathfrak
V(r)$, $\mathfrak u = \kappa_\mathfrak V(u)$, $\mathfrak v =
\kappa_\mathfrak V(v)$ in $ \kappa_\mathfrak V(mH) \le
\GL_{3588}(13)$.

Fixing a basis in each irreducible constituent $N\psi_j$ of
$\mathfrak W$ we get a basis $\mathcal B_W$ of $\mathfrak W$. It
induces a representation $\kappa_\mathfrak W : mE \rightarrow
\GL_{3588}(13)$ of $mE$. By Lemma \ref{l. reps13Fi23}(g)
$\mathfrak V_{|mD} \cong \mathfrak W_{|mD}$. Let $Y =
\GL_{3588}(13)$. Applying now Parker's isomorphism test of
Proposition 6.1.6 of \cite{michler} by means of the MAGMA command
$$\verb"IsIsomorphic(GModule(sub<Y|V(u),V(v)>),GModule(sub<Y|W(u),W(v)>))"$$
one obtains the transformation matrix $\mathcal T_1$ satisfying
$\mathfrak u = \kappa_{\mathfrak W}(u)^{\mathcal T_1}$ and
$\mathfrak v = \kappa_{\mathfrak W}(v)^{\mathcal T_1}$.

(a) Let $\mathfrak {mD} = \langle \mathfrak u, \mathfrak v
\rangle$, $\mathfrak {mH} = \langle \mathfrak u, \mathfrak v,
\mathfrak r \rangle$. Let $\mathcal D = C_Y(\mathfrak {mD})$ and
$\mathcal H = C_Y(\mathfrak {mH})$. Let $\delta_{12a}$,
$\delta_{12b}$ and $\delta_{23a}$, $\delta_{23b}$ be two distinct
copies of the irreducible characters $\delta_{12}$ and
$\delta_{23}$ of $mD$, respectively. Using $PG_1$ and MAGMA we
checked that the irreducible characters $\pi_8$ and $\pi_{16}$ of
$mH$ have the following restrictions to $mD = \langle u, v
\rangle$:

$${\pi_8}_{|mD} = \delta_{12a} + \delta_{23a} + \delta_{27} + \delta_{39}, \quad {\pi_{16}}_{|mD} = \delta_{12b} + \delta_{23b} + \delta_{29} + \delta_{54} +
\delta_{69},$$ where the irreducible characters $\delta_{12}$,
$\delta_{23}$, $\delta_{27}$, $\delta_{29}$, $\delta_{39}$,
$\delta_{54}$, and $\delta_{69}$ of $mD \cong G_2(3) \times S_3$
have degrees $78$, $156$, $182$, $182$, $364$, $728$ and $1664$,
respectively.

(b) Furthermore, Schur's Lemma asserts that $\mathcal H_1 =
\alpha(\mathcal H)$ is generated by the two blocked diagonal
matrices given in the statement because $2$ is a primitive element
of the multiplicative group $F^{*}$ of $F = GF(13)$.

(c)  Let $\mathfrak s_1 = \mathcal T_1^{-1} \kappa_\mathfrak W (s)
\mathcal T_1$. Let $\mathfrak {mE} = \langle \mathfrak {mD},
\mathfrak s_1 \rangle$ and $\mathcal E = C_Y(\mathfrak {mE})$. Let
$\delta_{12a}$, $\delta_{12b}$ and $\delta_{23a}$, $\delta_{23b}$
are two distinct copies of the irreducible characters
$\delta_{12}$ and $\delta_{23}$ of $mD$, respectively. Using
$PG_1$ and MAGMA we checked that the irreducible characters
$\psi_4$, $\psi_{10}$, $\psi_{19}$, $\psi_{36}$ and $\pi_{61}$ of
$mE$ have the following restrictions to $mD = \langle u, v
\rangle$:
\begin{eqnarray*}
&&{\psi_4}_{|_mD} = \delta_{12a}, \quad {\psi_{10}}_{|_mD} = \delta_{23a}, \quad {\psi_{19}}_{|_mD} = \delta_{12b} + \delta_{29},\\
&&{\psi_{36}}_{|_mD} = \delta_{27} + \delta_{54}, \quad {\psi_{61}}_{|_mD} = \delta_{23b} + \delta_{39} + \delta_{69}.\\
\end{eqnarray*}

Now Schur's Lemma implies that $\mathcal E_1 = \alpha(\mathcal E)$
is generated by the five blocked diagonal matrices $b_j$ given in
the statement.

(d) Every $\mathcal H$-$\mathcal E$ double coset representative is
of the form $diag(A,B,1,1,1,1,v)$ for some $A, B \in Y$ and $v\in
F^{*}$. By multiplying from left and right, we observe that
$diag(A,B,1,1,1,1,v)$ and $diag(A',B',1,1,1,1,v)$ represent the
same double coset if and only if the first columns of $A$ and $B$
are each a scalar multiple of the first columns of $A'$ and $B'$,
respectively. So, we have $12$ choices for $v$, and
$|\GL(2,13)|/12) = 2184$ choices for $A$ and $B$. Thus there are
$2184^2 \cdot 12 = 57238272$ $\mathcal H$-$\mathcal E$ double
cosets.

(e) By Theorem 7.2.2 of \cite{michler} the irreducible
representations of the free product $mH*_{mD}mE$ of the groups
$mH$ and $mE$ with amalgamated subgroup $mD$ are described by the
$\mathcal H$-$\mathcal E$ double coset representatives $T$ of
$\mathcal D$. The elements $r$ and $f = (u^3 v s v)^9$ are two
commuting involutions of $G_1 \cong \Fi_{23}$ by Lemma \ref{l.
reps13Fi23}(o). Let $\mathfrak u$ and $\mathfrak f$ be their
matrices in $\mathfrak {mH}$ and $\mathfrak {mE}$, respectively.
If $T = diag(\left( \begin{smallmatrix} a & c \\  b & d
\end{smallmatrix} \right), \left( \begin{smallmatrix} p & t \\  q
& u \end{smallmatrix} \right),1,1,1,1,v) \quad$ describes a
$3588$-dimensional representation of $G_1$ over $F$ then (*)
$(\mathfrak r, \mathcal T^{-1}\mathfrak f\mathcal T) = 1$ holds,
where $\mathcal T \in \GL_{3588}(13)$ corresponds to $T$.

Since $\mathfrak V_{|mD} \cong \mathfrak W_{|mD}$ is a direct sum
of $9$ irreducible $FmD$-modules both matrices $\mathfrak r$ and
$\mathfrak f$ consist of $81$ blocks $R_{i,j}$ and $F_{i,j}$, $1
\le i, j \le 9$, respectively, such that all diagonal blocks
$R_{i,i}$ and $F_{i,i}$ are non zero. Furthermore a non diagonal
block $R_{i,j}$ of $\mathfrak r$ is non zero if and only if the
$i$-th irreducible and the $j$-th irreducible representations of
$mD$ appear in the restriction of an irreducible representation of
$mH$ to $mD$. A similar description holds for the blocks of
$\mathfrak f$. Hence the system of equations in the proofs of (a)
and (c) imply that

{\renewcommand{\arraystretch}{0.5}
\begin{align*}
\mathfrak{r} = \left( \begin{array}{*{9}{c@{\,}}c}
R_{1,1}&.&R_{1,3}&.&R_{1,5}&.&R_{1,7}&.&.\\
.&R_{2,2}&.&R_{2,4}&.&R_{2,6}&.&R_{2,8}&R_{2,9}\\
R_{3,1}&.&R_{3,3}&.&R_{3,5}&.&R_{3,7}&.&.\\
.&R_{4,2}&.&R_{4,4}&.&R_{4,6}&.&R_{4,8}&R_{4,9}\\
R_{5,1}&.&R_{5,3}&.&R_{5,5}&.&R_{5,7}&.&.\\
.&R_{6,2}&.&R_{6,4}&.&R_{6,6}&.&R_{6,8}&R_{6,9}\\
R_{7,1}&.&R_{7,3}&.&R_{7,5}&.&R_{7,7}&.&.\\
.&R_{8,2}&.&R_{8,4}&.&R_{8,6}&.&R_{8,8}&R_{8,9}\\
.&R_{9,2}&.&R_{9,4}&.&R_{9,6}&.&R_{9,8}&R_{9,9}\\
\end{array} \right), \\
\mathfrak{f} = \left( \begin{array}{*{9}{c@{\,}}c}
F_{1,1}&.&.&.&.&.&.&.&.\\
.&F_{2,2}&.&.&.&F_{2,6}&.&.&.\\
.&.&F_{3,3}&.&.&.&.&.&.\\
.&.&.&F_{4,4}&.&.&F_{4,7}&.&F_{4,9}\\
.&.&.&.&F_{5,5}&.&.&F_{5,8}&.\\
.&F_{6,2}&.&.&.&F_{6,6}&.&.&.\\
.&.&.&F_{7,4}&.&.&F_{7,7}&.&F_{7,9}\\
.&.&.&.&F_{8,5}&.&.&F_{8,8}&.\\
.&.&.&F_{9,4}&.&.&F_{9,7}&.&F_{9,9}\\
\end{array} \right).
\end{align*}
}

Let $e=(ad-bc)^{-1}$ and $g=(pu-tq)^{-1}$. Then $e \neq 0 \neq g$.
For each integer $k$ let $I_k$ denote the $k \times k$ identity
matrix over $F$. Then
\begin{align*}
\mathcal T^{-1} = \left( \begin{array}{*{9}{c@{\,}}c}
ed(I_{78})&-ec(I_{78})&.&.&.&.&.&.&.\\
-eb(I_{78})&ea(I_{78})&.&.&.&.&.&.&.\\
.&.&gu(I_{156})&-gt(I_{156})&.&.&.&.&.\\
.&.&-gq(I_{156})&gp(I_{156})&.&.&.&.&.\\
.&.&.&.&I_{182}&.&.&.&.\\
.&.&.&.&.&I_{182}&.&.&.\\
.&.&.&.&.&.&I_{364}&.&.\\
.&.&.&.&.&.&.&I_{728}&.\\
.&.&.&.&.&.&.&.&v^{-1}(I_{1664})\\
\end{array} \right).
\end{align*}

Hence $\mathfrak f' = \mathcal T^{-1}\mathfrak f\mathcal T$ equals
the matrix

{\renewcommand{\arraystretch}{0.5}
$$
\left(
\begin{array}{*{9}{c@{\,}}c}
G_{1,1}&G_{1,2}&.&.&.&-ecF_{2,6}&.&.&.\\
G_{2,1}&G_{2,2}&.&.&.&eaF_{2,6}&.&.&.\\
.&.&G_{3,3}&G_{3,4}&.&.&-gtF_{4,7}&.&-vgtF_{4,9}\\
.&.&G_{4,3}&G_{4,4}&.&.&gpF_{4,7}&.&vgpF_{4,9}\\
.&.&.&.&F_{5,5}&.&.&F_{5,8}&.\\
G_{6,1}&G_{6,2}&.&.&.&F_{6,6}&.&.&.\\
.&.&G_{7,3}&G_{7,4}&.&.&F_{7,7}&.&vF_{7,9}\\
.&.&.&.&F_{8,5}&.&.&F_{8,8}&.\\
.&.&G_{9,3}&G_{9,4}&.&.&v^{-1}F_{9,7}&.&F_{9,9}\\
\end{array} \right),
$$
}

where

\begin{eqnarray*}
&& G_{1,1} = e(adF_{1,1} - bcF_{2,2}), \quad G_{1,2} = ecd(F_{1,1}-F_{2,2}), \\
&& G_{2,1} = -eab(F_{1,1} - F_{2,2}), \quad G_{2,2} = -e(bcF_{1,1} - adF_{2,2}), \\
&& G_{6,1} = bF_{6,2}, \quad G_{6,2} = dF_{6,2},\\
&&G_{3,3} = g(puF_{3,3} - qtF_{4,4}), \\
&& G_{3,4} = gut(F_{3,3} - F_{4,4}), \\
&& G_{4,3} = -gpq(F_{3,3} - F_{4,4}), \\
&& G_{4,4} = -g(tqF_{3,3} - puF_{4,4}), \\
&& G_{7,3} = qF_{7,4}, \quad G_{7,4} = uF_{7,4}, \\
&& G_{9,3} = qv^{-1}F_{9,4}, \quad G_{9,4} = uv^{-1}F_{9,4}.
\end{eqnarray*}

Now (*) implies the following equations
\begin{align*}
(1,1) & : G_{1,1}R_{1,1} = R_{1,1}G_{1,1}, \\
(1,2) & : G_{1,2}R_{2,2} + (-ec)F_{2,6}R_{6,2} = R_{1,1}G_{1,2}, \\
(6,1) & : G_{6,1}R_{1,1} = R_{6,2}G_{2,1} + R_{6,6}G_{6,1}, \\
(8,1) & : F_{8,5}R_{5,1} = R_{8,2}G_{2,1} + R_{8,6}G_{6,1}.
\end{align*}

Inserting the previous equations into these $4$ equations yields:
\begin{align*}
(1,1) & : e(adF_{1,1} - bcF_{2,2}) R_{1,1} = eR_{1,1}(adF_{1,1} - bcF_{2,2}) , \\
(1,2) & : ec(d(F_{1,1}-F_{2,2}) R_{2,2} -F_{2,6}R_{6,2}) = ecdR_{1,1}(F_{1,1}-F_{2,2}), \\
(6,1) & : bF_{6,2}R_{1,1} = b(-eaR_{6,2}(F_{1,1} - F_{2,2}) + R_{6,6}F_{6,2}), \\
(8,1) & : F_{8,5}R_{5,1} = b(-eaR_{8,2}(F_{1,1} - F_{2,2}) +
R_{8,6}F_{6,2}).
\end{align*}

Since $e^{-1} = ad-bc \neq 0$, at least one of $ad$ or $bc$ is
nonzero. Suppose $ad$ is zero, and $bc$ is nonzero. Then $(1,1)$
yields $F_{1,1}R_{1,1} = R_{1,1}F_{1,1}$. A MAGMA calculation
disproves this equation. Hence $ad \neq 0$. If $bc$ is zero and
$ad$ is nonzero, then $(1,1)$ equation implies $F_{2,2}R_{1,1} =
R_{1,1}F_{2,2}$ which is also wrong by MAGMA. Therefore all
$a,b,c,d$ are nonzero. We modify $T$ so that $a=1$ by multiplying
some power of $diag( \left(\begin{smallmatrix} 2 & 0 \\  0 & 1
\end{smallmatrix} \right), \left( \begin{smallmatrix} 1 & 0 \\  0
& 1 \end{smallmatrix} \right), 1,1,1,1,1)$ from the right.

Since $ec$ is nonzero, it can be cancelled on both sides of
equation $(1,2)$. Hence
\begin{align*}
(1,2) & : d(F_{1,1}-F_{2,2}) R_{2,2} -F_{2,6}R_{6,2} = dR_{1,1}(F_{1,1}-F_{2,2}).\\
\end{align*}
Using MAGMA it can be verified that this equation holds only for
$d = 2$. As $b \neq 0$ and $a = 1$ equation (6,1) implies
\begin{align*}
F_{6,2}R_{1,1} = -eR_{6,2}(F_{1,1} - F_{2,2}) + R_{6,6}F_{6,2}.
\end{align*}
By MAGMA it has the solution $e = 11$. Now equation (8,1) implies
that
\begin{align*}
F_{8,5}R_{5,1} = b(-11R_{8,2}(F_{1,1} - F_{2,2}) + R_{8,6}F_{6,2}).\\
\end{align*}
This equation has the solution $b = 4$ by MAGMA. From the equation
$ad-bc = e^{-1}$ we now deduce that $c = 12$.

In order to determine the remaining coefficients of the matrix $T$
we use the following matrix equations derived from (*).
\begin{align*}
(9,8) & : G_{9,4}R_{4,8} + F_{9,9}R_{9,8} = R_{9,8}F_{8,8}, \\
(9,9) & : G_{9,4}R_{4,9} + F_{9,9}R_{9,9} = vgp R_{9,4}F_{4,9} + R_{9,9}F_{9,9}, \\
(4,5) & : G_{4,3}R_{3,5} + gpF_{4,7}R_{7,5} = R_{4,8}F_{8,5}.
\end{align*}
Inserting the first set of equations yields:
\begin{align*}
(9,8) & : uv^{-1}F_{9,4}R_{4,8} + F_{9,9}R_{9,8} = R_{9,8}F_{8,8}, \\
(9,9) & : uv^{-1}F_{9,4}R_{4,9} + F_{9,9}R_{9,9} = vgp R_{9,4}F_{4,9} + R_{9,9}F_{9,9}, \\
(4,5) & : -gpq(F_{3,3} - F_{4,4})R_{3,5} + gpF_{4,7}R_{7,5} =
R_{4,8}F_{8,5}.
\end{align*}
By MAGMA the equation (9,8) has the solution $uv^{-1} = 10$. Now
MAGMA asserts that the equation (9,9) has the solution $vgp = 11$.
Hence $p \neq 0$. By multiplying some power of $diag( \left(
\begin{smallmatrix} 2 & 0 \\  0 & 1
\end{smallmatrix} \right), \left( \begin{smallmatrix} 1 & 0 \\  0
& 1 \end{smallmatrix} \right), 1,1,1,1,1)$ from the right we can
modify $T$ so that $p=1$. Thus $vg = 11$ and the third equation
(4,5) implies that
\begin{align*}
(4,5) & : g( -q(F_{3,3} - F_{4,4})R_{3,5} + F_{4,7}R_{7,5}) = R_{4,8}F_{8,5}.\\
\end{align*}
This is non linear equation in the unknowns $q$ and $g$ has a
unique solution $q = 6$ and $g = 11$ as has been checked by
running with MAGMA through all $13^2$ cases. From $6 = g^{-1} =
pu-tq = 10 - 6t$ we now deduce that $t = 5$. This completes the
determination of the coefficients of the two matrices $w$ and $z$
of statement (e). The remaining assertions are now clear.
\end{proof}

\begin{proposition}\label{prop. rep5083Fi23}
Keep the notation of Lemmas \ref{l. aut2Fi22}, \ref{l. gensFi24}
and \ref{l. reps13Fi23}. Let $PG_1$ be the faithful permutation
representation of the simple group $G_1 = \langle x,y,q,w\rangle =
\langle u,v,r,s \rangle$ of degree $31671$ with stabilizer $H_1 =
\langle x,y,q\rangle$. Let $mH = \langle u,v,r \rangle$, $mD =
\langle u,v,r \rangle$ and $mE = \langle u,v,s \rangle$. Let
$F^{*}$ be the multiplicative group of the prime field $F =
\GF(13)$. Let $Y = \GL_{5083}(13)$.

Let $\mathfrak V$ and $\mathfrak W$ be the up to isomorphism
uniquely determined faithful semi-simple $5083$-dimensional
modules of $mH$ and $mE$ over $F$ corresponding to the
restrictions ${\tau_{\bf 4}}_{|mH}$ and ${\tau_{\bf 4}}_{|mE}$ of
the irreducible character $\tau_4$ of $G_1$, respectively.

Let $\kappa_\mathfrak V : mH \rightarrow \GL_{5083}(13)$ and
$\kappa_\mathfrak W : mE \rightarrow \GL_{5083}(13)$ be the
representations of $mH$ and $mE$ afforded by the modules
$\mathfrak V$ and $\mathfrak W$, respectively.

Let $\mathfrak r = \kappa_\mathfrak V(r)$, $\mathfrak u =
\kappa_\mathfrak V(u)$, $\mathfrak v = \kappa_\mathfrak V(v)$ in $
\kappa_\mathfrak V(mH) \le \GL_{5083}(13)$.

Then $\mathfrak V_{|mD} \cong \mathfrak W_{|mD}$, and there is a
transformation matrix $\mathcal T \in \GL_{5083}(13)$ such that
$$
\mathfrak u = \mathcal T^{-1} \kappa_\mathfrak W (u) \mathcal T,
\mathfrak v = \mathcal T^{-1} \kappa_\mathfrak W(v) \mathcal T.
$$

Let $\mathfrak {mD} = \langle \mathfrak u, \mathfrak v \rangle$,
$\mathfrak {mH} = \langle \mathfrak u, \mathfrak v, \mathfrak r
\rangle$. Let $\mathcal D = C_Y(\mathfrak {mD})$ and $\mathcal H =
C_Y(\mathfrak {mH})$. Let $\mathfrak s_1 = \mathcal T^{-1}
\kappa_\mathfrak W (s) \mathcal T$. Let $\mathfrak {mE} = \langle
\mathfrak {mD}, \mathfrak s_1 \rangle$ and $\mathcal E =
C_Y(\mathfrak {mE})$. Then the following statements hold:

\begin{enumerate}
\item[\rm(a)] There is an isomorphism $$\alpha: \mathcal D
\rightarrow \mathcal D_1 = \GL_2(13) \times {F^{*}}^8 \le
\GL_{10}(13).$$

\item[\rm(b)] $\mathcal H_1 = \alpha(\mathcal H)$ is generated by
the two blocked diagonal matrices\\
\begin{align*}
a_1 = diag( \left( \begin{smallmatrix} 2 & 0 \\  0 & 1
\end{smallmatrix} \right), 1,2,1,2,2,1,2,1) \quad \mbox{and} \quad
a_2 = diag( \left(
\begin{smallmatrix} 1 & 0 \\  0 & 2 \end{smallmatrix} \right),
2,1,2,1,1,2,1,2),
\end{align*}

\item[\rm(c)] $\mathcal E_1 = \alpha(\mathcal E)$ is generated by
the four blocked diagonal matrices
\begin{align*}
b_1 = diag( \left( \begin{smallmatrix} 2 & 0 \\  0 & 1
\end{smallmatrix} \right), 1,1,1,1,1,1,1,1), \quad
b_2 = diag( \left(
\begin{smallmatrix} 1 & 0 \\  0 & 2 \end{smallmatrix} \right),
2,1,1,1,1,1,1,1), \\
b_3 = diag( \left( \begin{smallmatrix} 1 & 0 \\  0 & 1
\end{smallmatrix} \right), 1,2,1,2,2,2,1,1), \quad
b_4 = diag( \left(
\begin{smallmatrix} 1 & 0 \\  0 & 1 \end{smallmatrix} \right),
1,1,2,1,1,1,2,2).
\end{align*}

\item[\rm(d)] $\mathcal D$ has $2184 \times 12^4 = 4587424$
$\mathcal H$-$\mathcal E$ double cosets.

\item[\rm(e)] The free product $mH*_{mD}mE$ of $mH$ and $mE$ with
amalgamated subgroup $mD$ has an irreducible $5083$-dimensional
representation over $F$ which induces an irreducible
representation of $G_1$. It corresponds to the $\mathcal
H$-$\mathcal E$ double coset representative

$$\mathcal F = diag(w,1^{78},1^{91},11^{156},6^{273},11^{273},1^{728},1^{1456},1^{1664}) \in \GL_{5083}(13), \quad
{where}$$

{\renewcommand{\arraystretch}{0.5} \scriptsize
$$
w = \left( \begin{array}{*{2}{c@{\,}}c}
 1^{182} & 9^{182}\\
11^{182} & 9^{182}
\end{array} \right),
$$
} and $a^n$ denotes a diagonal $n \times n$ matrix with unique
diagonal non zero entry $a \in GF(13)$.

Let $\mathfrak s = \mathcal F^{-1} \mathfrak s_1 \mathcal F$ and
$\mathfrak G_1 = \langle \mathfrak u, \mathfrak v, \mathfrak r,
\mathfrak s \rangle$. Inserting these four generating matrices of
$\mathfrak G_1$ into the formulas of Lemma \ref{l. reps13Fi23}(f)
one obtains the matrices $\mathfrak x_{5083}$, $\mathfrak
y_{5083}$, $\mathfrak q_{5083}$, and $\mathfrak w_{5083}$ of the
original generators $x$, $y$, $q$ and $w$ of $G_1 = \Fi_{23}$ as
words in the generators $u$, $v$, $r$ and $s$. The matrices
$\mathfrak x_{5083}$, $\mathfrak y_{5083}$, $\mathfrak q_{5083}$,
and $\mathfrak w_{5083}$ can be downloaded from the first author's
website\\
$\verb"http://www.math.yale.edu/~hk47/Fi24/index.html"$.
\end{enumerate}
\end{proposition}

\begin{proof}
Let $\mathfrak V$ be the up to isomorphism uniquely determined
faithful semi-simple $5083$-dimensional module of $mH$ over $F =
GF(13)$ corresponding to the restriction ${\tau_{\bf 4}}_{|mH}$.
By Lemma \ref{l. reps13Fi23}(g) ${\tau_{\bf 4}}_{|mH} = \pi_{12} +
\pi_{15}$.

The irreducible characters $\pi_{11}$ and $\pi_{15}$ of $mH$ are
constituents of the permutation characters $1_{mH_{11}}^{mH}$ and
$1_{mH_{15}}^{mH}$ of the respective subgroups $mH_{11}$ and
$mH_{15}$ of $mH$ determined in Lemma \ref{l. reps13Fi23}(i).
Using a stand-alone program of the first author which is based on
Algorithm 5.7.1 of \cite{michler} we calculated the primitive
idempotents of the endomorphism rings of these permutation
modules. Thus we obtained the corresponding irreducible
representations $M(\pi_{11})$ and $M(\pi_{15})$ of the respective
dimensions $2275$ and $2808$ over $F = GF(13)$. Hence the
irreducible $FmH$-module $M(\pi_{12})$ is the tensor product of
$M(\pi_{11})$ with the linear character $\pi_2$ of $mH$ over $F$.
Thus $\mathfrak V = M(\pi_{12}) \oplus M(\pi_{15})$.

Let $\mathfrak W$ be the up to isomorphism uniquely determined
faithful semi-simple $5083$-dimensional module of $mE$ over $F$
corresponding to the restriction ${\tau_{\bf 4}}_{|mE} = \psi_{14}
+ \psi_{22} + \psi_{45} + \psi_{74}$, see Lemma \ref{l.
reps13Fi23}(i).

The irreducible characters $\psi_{14}$, $\psi_{22}$, $\psi_{45}$
and $\psi_{74}$ of $mE$ are constituents of the permutation
characters $1_{mE_{14}}^{mE}$, $1_{mE_{22}}^{mE}$,
$1_{mE_{45}}^{mE}$ and $1_{mE_{74}}^{mE}$, respectively, by Lemma
\ref{l. reps13Fi23}(l). Using a stand alone program of the first
author which is based on Algorithm 5.7.1 of \cite{michler} we
calculated the primitive idempotents of the endomorphism rings of
these four permutation modules. Thus we obtained the corresponding
irreducible representations $N(\psi_{14})$, $N(\psi_{22})$,
$N(\psi_{45})$ and $N(\psi_{74})$ of the respective dimensions
$182$, $260$, $1365$ and $3276$ over $F$. Hence $$\mathfrak W =
N(\psi_{14}) \oplus N(\psi_{22}) \oplus N(\psi_{45}) \oplus
N(\psi_{74}).$$

Fixing a basis in each irreducible constituent $M\pi_k$ of
$\mathfrak V$ we get a basis $\mathcal B_V$ of $\mathfrak V$. It
induces a representation $\kappa_\mathfrak V : mH \rightarrow
\GL_{5083}(13)$ of $mH$. Let $\mathfrak r = \kappa_\mathfrak
V(r)$, $\mathfrak u = \kappa_\mathfrak V(u)$, $\mathfrak v =
\kappa_\mathfrak V(v)$ in $ \kappa_\mathfrak V(mH) \le
\GL_{5083}(13)$.

Fixing a basis in each irreducible constituent $N\psi_j$ of
$\mathfrak W$ we get a basis $\mathcal B_W$ of $\mathfrak W$. It
induces a representation $\kappa_\mathfrak W : mE \rightarrow
\GL_{5083}(13)$ of $mE$. By Lemma \ref{l. reps13Fi23}(h)
$\mathfrak V_{|mD} \cong \mathfrak W_{|mD}$. Let $Y =
\GL_{5083}(13)$. Applying now Parker's isomorphism test of
Proposition 6.1.6 of \cite{michler} by means of the MAGMA command
$$\verb"IsIsomorphic(GModule(sub<Y|V(u),V(v)>),GModule(sub<Y|W(u),W(v)>))"$$
one obtains the transformation matrix $\mathcal T_1$ satisfying
$\mathfrak u = \kappa_{\mathfrak W}(u)^{\mathcal T_1}$ and
$\mathfrak v = \kappa_{\mathfrak W}(v)^{\mathcal T_1}$.

(a) Let $\mathfrak {mD} = \langle \mathfrak u, \mathfrak v
\rangle$, $\mathfrak {mH} = \langle \mathfrak u, \mathfrak v,
\mathfrak r \rangle$. Let $\mathcal D = C_Y(\mathfrak {mD})$ and
$\mathcal H = C_Y(\mathfrak {mH})$. Let $\delta_{32a}$ and
$\delta_{32b}$ be two distinct copies of the irreducible character
$\delta_{32}$ of $mD$. Using $PG_1$ and MAGMA it can be checked
that the irreducible characters $\pi_{12}$ and $\pi_{15}$ of $mH$
have the following restrictions to $mD = \langle u, v \rangle$:
\begin{eqnarray*}
&& {\pi_{12}}_{|mD} = \delta_{16} + \delta_{32a} + \delta_{34} + \delta_{37} + \delta_{66}, \\
&&{\pi_{15}}_{|mD} = \delta_{11} + \delta_{23} + \delta_{32b} + \delta_{53} + \delta_{69},\\
\end{eqnarray*}
where the irreducible characters $\delta_{11}$, $\delta_{16}$,
$\delta_{23}$, $\delta_{32}$, $\delta_{34}$, $\delta_{37}$,
$\delta_{53}$, $\delta_{66}$, and $\delta_{69}$ of $mD \cong
G_2(3) \times S_3$ have degrees $78$, $91$, $156$, $182$, $273$,
$273$, $728$, $1456$ and $1664$, respectively.

Since $\mathfrak V_{|mD}$ is a semi-simple $FmD$-module the
Theorem 2.1.27 of \cite{michler} implies that there is an
isomorphism
$$\alpha: \mathcal D \rightarrow \mathcal D_1 = \GL_2(13) \times {F^{*}}^8 \le \GL_{10}(13).$$

(b) Furthermore, Schur's Lemma asserts that $\mathcal H_1 =
\alpha(\mathcal H)$ is generated by the two blocked diagonal
matrices $a_1$ and $a_2$ given in the statement because $2$ is a
primitive element in the multiplicative group $F^{*}$ of $F =
GF(13)$.

(c)  Let $\mathfrak s_1 = \mathcal T_1^{-1} \kappa_\mathfrak W (s)
\mathcal T_1$. Let $\mathfrak {mE} = \langle \mathfrak {mD},
\mathfrak s_1 \rangle$ and $\mathcal E = C_Y(\mathfrak {mE})$. Let
$\delta_{32a}$ and $\delta_{32b}$ be two distinct copies of the
irreducible character $\delta_{32}$ of $mD$. Using $PG_1$ and
MAGMA it can be checked that the irreducible characters
$\psi_{14}$, $\psi_{22}$, $\psi_{45}$ and $\psi_{74}$ of $mE$ have
the following restrictions to $mD = \langle u, v \rangle$:
\begin{eqnarray*}
&&{\psi_{14}}_{|_mD} = \delta_{32a}, \quad {\psi_{22}}_{|_mD} = \delta_{11} + \delta_{32b}, \\
&&{\psi_{45}}_{|_mD} = \delta_{16} + \delta_{34} + \delta_{37} + \delta_{53}, \quad {\psi_{74}}_{|_mD} = \delta_{23} + \delta_{66} + \delta_{69}.\\
\end{eqnarray*}

Now Schur's Lemma implies that $\mathcal E_1 = \alpha(\mathcal E)$
is generated by the four blocked diagonal matrices $b_j$ given in
the statement.

(d) Every $\mathcal H$-$\mathcal E$ double coset representative is
of the form $$diag(A,1,v_2,v_3,v_4,v_5,1,1,1)$$ for some $A \in Y$
and $v\in F^{*}$. By multiplying from left and right, we can find
out that $diag(A,1,v_2,v_3,v_4,v_5,1,1,1)$ and
$diag(A',1,v_2,v_3,v_4,v_5,1,1,1)$ represent the same double coset
if and only if the first column of $A$ is a scalar multiple of
that of $A'$. So, we have $12^4$ choices for $v_i$, and
$|\GL(2,13)|/12 = 2184$ choices for $A$. Thus there are $2184
\cdot 12^4 = 4587424$ $\mathcal H$-$\mathcal E$ double cosets.

(e) By Theorem 7.2.2 of \cite{michler} the irreducible
representations of the free product $mH*_{mD}mE$ of the groups
$mH$ and $mE$ with amalgamated subgroup $mD$ are described by the
$\mathcal H$-$\mathcal E$ double coset representatives $T$ of
$\mathcal D$. The elements $r$ and $f = (u^3 v s v)^9$ are two
commuting involutions of $G_1 \cong \Fi_{23}$ by Lemma \ref{l.
reps13Fi23}(o). Let $\mathfrak r$ and $\mathfrak f$ be their
matrices in $\mathfrak {mH}$ and $\mathfrak {mE}$, respectively.
If $T_1 = diag(\left( \begin{smallmatrix} a & c \\  b & d
\end{smallmatrix} \right),1,v_2,v_3,v_4,v_5,1,1,1) \quad$
describes a $5083$-dimensional irreducible representation of $G_1$
over $F$ then (**) $(\mathfrak r, \mathcal T_1^{-1}\mathfrak
f\mathcal T_1) = 1$ holds.

Since $\mathfrak V_{|mD} \cong \mathfrak W_{|mD}$ is a direct sum
of $10$ irreducible $FmD$-modules both matrices $\mathfrak r$ and
$\mathfrak f$ consist of $100$ blocks $R_{i,j}$ and $F_{i,j}$, $1
\le i, j \le 10$ such that all diagonal blocks $R_{i,i}$ and
$F_{i,i}$ are nontrivial. Furthermore, a non diagonal block
$R_{i,j}$ of $\mathfrak r$ is non zero if and only if the $i$-th
irreducible and the $j$-th irreducible representations of $mD$
appear in the restriction of an irreducible representation of $mH$
to $mD$. A similar description holds for the blocks of $\mathfrak
f$. Hence the system of equations in the proofs of (a) and (c)
imply that

{\renewcommand{\arraystretch}{0.5}
\begin{align*}
\mathfrak{r} = \left( \begin{array}{*{10}{c@{\,}}c}
R_{1,1}&.&.&R_{1,4}&.&R_{1,6}&R_{1,7}&.&R_{1,9}&.\\
.&R_{2,2}&R_{2,3}&.&R_{2,5}&.&.&R_{2,8}&.&R_{2,10}\\
.&R_{3,2}&R_{3,3}&.&R_{3,5}&.&.&R_{3,8}&.&R_{3,10}\\
R_{4,1}&.&.&R_{4,4}&.&R_{4,6}&R_{4,7}&.&R_{4,9}&.\\
.&R_{5,2}&R_{5,3}&.&R_{5,5}&.&.&R_{5,8}&.&R_{5,10}\\
R_{6,1}&.&.&R_{6,4}&.&R_{6,6}&R_{6,7}&.&R_{6,9}&.\\
R_{7,1}&.&.&R_{7,4}&.&R_{7,6}&R_{7,7}&.&R_{7,9}&.\\
.&R_{8,2}&R_{8,3}&.&R_{8,5}&.&.&R_{8,8}&.&R_{8,10}\\
R_{9,1}&.&.&R_{9,4}&.&R_{9,6}&R_{9,7}&.&R_{9,9}&.\\
.&R_{10,2}&R_{10,3}&.&R_{10,5}&.&.&R_{10,8}&.&R_{10,10}\\
\end{array} \right),\\
\mathfrak{f} = \left( \begin{array}{*{10}{c@{\,}}c}
F_{1,1}&.&.&.&.&.&.&.&.&.\\
.&F_{2,2}&F_{2,3}&.&.&.&.&.&.&.\\
.&F_{3,2}&F_{3,3}&.&.&.&.&.&.&.\\
.&.&.&F_{4,4}&.&F_{4,6}&F_{4,7}&F_{4,8}&.&.\\
.&.&.&.&F_{5,5}&.&.&.&F_{5,9}&F_{5,10}\\
.&.&.&F_{6,4}&.&F_{6,6}&F_{6,7}&F_{6,8}&.&.\\
.&.&.&F_{7,4}&.&F_{7,6}&F_{7,7}&F_{7,8}&.&.\\
.&.&.&F_{8,4}&.&F_{8,6}&F_{8,7}&F_{8,8}&.&.\\
.&.&.&.&F_{9,5}&.&.&.&F_{9,9}&F_{9,10}\\
.&.&.&.&F_{10,5}&.&.&.&F_{10,9}&F_{10,10}\\
\end{array} \right).
\end{align*}
}

Let $e=(ad-bc)^{-1}$ and $g=(pu-tq)^{-1}$. Then $e \neq 0 \neq g$.
For each integer $k$ let $I_k$ denote the $k \times k$ identity
matrix over $F$. Then

{\renewcommand{\arraystretch}{0.5}
\begin{align*}
\mathcal{T}^{-1} = \left( \begin{array}{*{10}{c@{\,}}c}
ed(I_{182})&-ec(I_{182})&.&.&.&.&.&.&.&.\\
-eb(I_{182})&ea(I_{182})&.&.&.&.&.&.&.&.\\
.&.&I_{78}&.&.&.&.&.&.&.\\
.&.&.&v_2^{-1}(I_{91})&.&.&.&.&.&.\\
.&.&.&.&v_3^{-1}(I_{156})&.&.&.&.&.\\
.&.&.&.&.&v_4^{-1}(I_{273})&.&.&.&.\\
.&.&.&.&.&.&v_5^{-1}(I_{273})&.&.&.\\
.&.&.&.&.&.&.&I_{728}&.&.\\
.&.&.&.&.&.&.&.&I_{1456}&.\\
.&.&.&.&.&.&.&.&.&I_{1664}\\
\end{array} \right).
\end{align*}
}

Hence $\mathfrak f' = \mathcal T^{-1}\mathfrak f\mathcal T$ equals
the matrix

{\renewcommand{\arraystretch}{0.5}
$$
\left( \begin{array}{*{10}{c@{\,}}c}
G_{1,1}&G_{1,2}&G_{1,3}&.&.&.&.&.&.&.\\
G_{2,1}&G_{2,2}&G_{2,3}&.&.&.&.&.&.&.\\
bF_{3,2}&dF_{3,2}&F_{3,3}&.&.&.&.&.&.&.\\
.&.&.&F_{4,4}&.&G_{4,6}&G_{4,7}&v_2^{-1}F_{4,8}&.&.\\
.&.&.&.&F_{5,5}&.&.&.&v_3^{-1}F_{5,9}&v_3^{-1}F_{5,10}\\
.&.&.&G_{6,4}&.&F_{6,6}&G_{6,7}&v_4^{-1}F_{6,8}&.&.\\
.&.&.&G_{7,4}&.&G_{7,6}&F_{7,7}&v_5^{-1}F_{7,8}&.&.\\
.&.&.&v_2F_{8,4}&.&v_4F_{8,6}&v_5   F_{8,7}&F_{8,8}&.&.\\
.&.&.&.&v_3F_{9,5}&.&.&.&F_{9,9}&F_{9,10}\\
.&.&.&.&v_3F_{10,5}&.&.&.&F_{10,9}&F_{10,10}\\
\end{array} \right),
$$
}

where
\begin{eqnarray*}
&&G_{1,1} = e(ad F_{1,1} -bc F_{2,2}), \quad G_{1,2} = ecd(F_{1,1} - F_{2,2}), \quad G_{1,3} = -ec F_{2,3}, \\
&&G_{2,1} = -eab(F_{1,1} - F_{2,2}), \quad G_{2,2} = -e(bcF_{1,1} - adF_{2,2}), \quad G_{2,3} = eaF_{2,3},\\
&&G_{4,6} = v_2^{-1} v_4 F_{4,6},\quad G_{4,7} = v_2^{-1} v_5F_{4,7},\quad G_{6,4} = v_4^{-1} v_2 F_{6,4}, \\
&&G_{6,7} = v_4^{-1} v_5 F_{6,7}, \quad G_{7,4} = v_5^{-1} v_2F_{7,4}, \quad G_{7,6} = v_5^{-1} v_4 F_{7,6}.\\
\end{eqnarray*}

Now (**) implies the following equations
\begin{align*}
(10,1) & : F_{10,9} R_{9,1} = R_{10,2} G_{2,1} + R_{10,3} (bF_{3,2}), \\
(9,1) & : F_{9,9} R_{9,1} = R_{9,1} G_{1,1}, \\
(9,2) & : v_3 F_{9,5} R_{5,2} + F_{9,10} R_{10,2} = R_{9,1} G_{1,2}, \\
(8,1) & : v_2F_{8,4} R_{4,1} + v_4F_{8,6}R_{6,1} + v_5F_{8,7}R_{7,1} = R_{8,2}G_{2,1} + R_{8,3}(bF_{3,2}).\\
\end{align*}

Inserting the first set of equations yields:
\begin{align*}
(10,1) & : F_{10,9} R_{9,1} = b[ -ea R_{10,2} (F_{1,1} - F_{2,2})  + R_{10,3} F_{3,2}], \\
(9,1) & : F_{9,9} R_{9,1} = e R_{9,1} (ad F_{1,1} -bc F_{2,2}), \\
(9,2) & : v_3 F_{9,5} R_{5,2} + F_{9,10} R_{10,2} = ecd R_{9,1} (F_{1,1} - F_{2,2}). \\
(8,1) & : v_2F_{8,4} R_{4,1} + v_4F_{8,6}R_{6,1} + v_5F_{8,7}R_{7,1} = b[ -ea R_{8,2}(F_{1,1} - F_{2,2}) + R_{8,3}F_{3,2}].\\
\end{align*}

Equation $(10,1)$ has only one solution $b = 9$, $ea = 1$ in $F =
GF(13)$ as has been checked by means of MAGMA. Hence $a \neq 0$.
By multiplying some power of $diag( \left( \begin{smallmatrix} 2 &
0 \\  0 & 1 \end{smallmatrix} \right), 1,1,1,1,1,1,1,1)$ from the
right we can modify the matrix $T$ so that $a = 1$. Thus also $e =
1$. Now equation $(9,1)$ becomes
\begin{align*}
(9,1) & : F_{9,9} R_{9,1} = R_{9,1}(d F_{1,1} - 9c F_{2,2}).
\end{align*}
Another MAGMA calculation running through $13^2$ pairs $(c,d)$ of
elements of $F$ shows that $(9,1)$ has the unique solution $c =
11$ and $d = 9$.

Inserting the known values for $c$, $d$ and $e$ into equation
$(9,2)$ yields that $v_3 = 11$ as has been checked by
means of MAGMA. Finally, inserting the known values for $a$, $b$
and $e$ into equation $(8,1)$ yields
\begin{align*}
(8,1) & : v_2F_{8,4} R_{4,1} + v_4F_{8,6}R_{6,1} + v_5F_{8,7}R_{7,1} = 9[ - R_{8,2}(F_{1,1} - F_{2,2}) + R_{8,3}F_{3,2}].\\
\end{align*}
A MAGMA calculation shows that this equation has the unique
solution $v_2 = 1$, $v_4 = 6$, $v_5= 11$. This completes the proof
because all remaining statements of (e) are straightforward.
\end{proof}

\section{Construction of the irreducible subgroup $\mathfrak G$ of $\GL_{8671}(13)$}

In this section we construct the $8$ semi-simple representations
of the $2$-fold cover $A_1$ of the automorphism group $\rm
Aut(\Fi{22})$ of $\Fi_{22}$ corresponding to the $8$ compatible
pairs of Lemma \ref{l. gensFi24}(e). Since $H_1 = \langle
q,y\rangle = C_{G_1}(u) \cong A_1'$ for some involution $u$ of
$G_1 = \langle q,y,w\rangle$ and $A_1'$ has index $2$ in $A_1$ it
is not difficult to construct the irreducible constituents of
these representations of $A_1$ from the $3588$-dimensional
irreducible representation of $G_1 \cong \Fi_{23} = \langle
q,y,w\rangle$ by means of Clifford's Theorem. In fact, we
construct $8$ new matrices $t_i$ of order $2$ such that $K_i =
\langle G_1, t_i\rangle$ is an irreducible subgroup of
$\GL_{8671}(13)$. It turns out that only $K_3$ corresponding to
the compatible pair (3) of Lemma \ref{l. gensFi24}(e) may have a
Sylow $2$-subgroup which is isomorphic to a Sylow $2$-subgroup of
the extension group $E$ of Lemma \ref{l. M24-extensions}.

\begin{proposition}\label{prop. repsA1}
Let $\mathfrak G_1 = \langle \mathfrak y_{3588}, \mathfrak
q_{3588},\mathfrak w_{3588}\rangle$ be the simple subgroup of $Y =
\GL_{3588}(13)$ of order $2^{18}\cdot 3^{13}\cdot
5^2\cdot7\cdot11\cdot13\cdot17\cdot23$ constructed in Proposition
\ref{prop. rep3588Fi23}. Let $y_1 = \mathfrak y_{3588}$, $q_1 =
\mathfrak q_{3588}$. Let $\mathfrak H_1 = \langle y_1,
q_1\rangle$.

Let $A_1 = 2Aut(\Fi_{22}) = \langle a,b,c,d,e,f,g,h,i,z,t\rangle$
be the finitely presented group defined in Lemma \ref{l.
aut2Fi22}. Then the following assertions hold:
\begin{enumerate}
\item[\rm(a)] There is an isomorphism $\phi$ from the subgroup
$H_1 = \langle a,b,c,d,e,f,g,h,i,z\rangle$ of $A_1$ to $\mathfrak
H_1$.

\item[\rm(b)] $A_1 = \langle y, q, t \rangle$, where $y =
\phi^{-1}(y_1)$ and $q = \phi^{-1}(q_1)$.

\item[\rm(c)] There is a transformation matrix $\mathcal T \in Y$
such that
\begin{eqnarray*}
&&\mathcal T^{-1}\mathfrak y \mathcal T = diag(\mathfrak y_{78}, \mathfrak y_{1430}, \mathfrak y_{2080}) \in Y,\\
&&\mathcal T^{-1}\mathfrak q \mathcal T = diag(\mathfrak q_{78}, \mathfrak q_{1430}, \mathfrak q_{2080}) \in Y, \\
\end{eqnarray*}
where $y_k, q_k \in \GL_k(13)$ for $k \in \{78, 1430, 2080\}$.

\item[\rm(d)] $A_1$ has a faithful irreducible representation
$\lambda : A_1 \rightarrow \GL_{4160}(13)$ such that
\begin{eqnarray*}
&&\lambda(y) = diag(\mathfrak y_{2080},\phi(y^t)_{2080}) \in \GL_{4160}(13)),\\
&&\lambda(q) = diag(\mathfrak h_{2080},\phi(q^t)_{2080}) \in \GL_{4160}(13)),\\
&&\lambda(t) =  \left(
\begin{smallmatrix}
0 & I_{2080} \\ I_{2080} & 0
\end{smallmatrix}
\right).
\end{eqnarray*}
where $I_{2080}$ denotes the identity matrix of $\GL_{2080}(13)$.

\item[\rm(e)] The irreducible characters $\chi_3$, $\chi_4$,
$\chi_{11}$ and $\chi_{12}$ of Table \ref{Fi_24 ct A_1} of
respective degrees $78$, $78$, $1430$ and $1430$ are constituents
of the permutation character $1_U^{A_1}$ of the subgroup

$$U = \langle (q_1^2y_1^3q_1y_1^3)^4,(y_1^2q_1y_1^3q_1^2y_1q_1)^6,(y_1^4q_1y_1q_1y_1q_1y_1^2)^2\rangle$$

of $A_1$ of index $2358720$.

\item[\rm(f)] The tensor product $\chi_3 \otimes \chi_3$ contains
a copy of the irreducible character $\chi_{13}$. Furthermore,
$\chi_{14} = \chi_{13} \otimes \chi_2$, where $\chi_2$ is the non
trivial linear character of $A_1$.
\end{enumerate}
\end{proposition}

\begin{proof}
(a) In the simple subgroup $\mathfrak G_1$ of $Y = \GL_{3588}(13)$
let
\begin{eqnarray*}
&&x_1 = [(y_1q_1^2y_1q_1y_1q_1^2)^{11}(q_1^2y_1^2q_1y_1q_1y_1)^{11}(q_1y_1^2q_1y_1q_1y_1q_1y_1q_1)^4]^{12},\\
&&a_1 = (x_1y_1x_1)^7,\quad b_1 = [(q_1y_1)^2q_1y_1^3q_1^2y_1^3q_1y_1]^7,\quad c_1 = (y_1^2x_1y_1x_1y_1^3)^5,\\
&&d_1 = (q_1y_1q_1^2y_1q_1y_1q_1y_1q_1y_1^2q_1^2)^{15},\quad e_1 = (y_1x_1y_1^5x_1)^5,\\
&&f_1 = (y_1q_1y_1q_1^2y_1q_1^2y_1^2q_1y_1^4q_1^2)^5,\quad g_1 = (x_1y_1^2x_1y_1^3x_1)^7,\\
&&h_1 = (y_1^5x_1y_1x_1)^5,\quad i_1 = (q_1^2y_1^2q_1y_1q_1^2)^7.\\
\end{eqnarray*}

By Lemma \ref{l. gensFi24} and Proposition \ref{prop. rep3588Fi23}
the matrix subgroup $G_1 = \langle x,y,q,w\rangle$ of
$\GL_{782}(17)$ and the matrix subgroup $\mathfrak G_1 = \langle
\mathfrak y_{3588}, \mathfrak q_{3588},\mathfrak w_{3588}\rangle$
are isomorphic under the map $\theta$ sending $x$, $y$, $q$ and
$w$ to $x_1$, $y_1$, $q_1$ and $w_{3588}$, respectively. Thus
Lemma \ref{l. gensFi24} (b) implies that the normal subgroup $H_1
= \langle a,b,c,d,e,f,g,h,i,z\rangle$ of $A_1$ is isomorphic to
$\mathfrak H_1 = \langle a_1,b_1,c_1,d_1,e_1,f_1,g_1,h_1,i_1,z_1
\rangle$, where $z_1 = (x_1y_1^2)^7$. In particular, the map
$\phi: H_1 \rightarrow \mathfrak H_1$ sending $a$, $b$, \ldots,
$z$ to $a_1$, $b_1$, \ldots, $z_1$ is an isomorphism such that $x
= \phi^{-1}(x_1)$, $y = \phi^{-1}(y_1)$, $q = \phi^{-1}(q_1)$.

Statement (b) is an immediate consequence of (a) and Lemma \ref{l.
aut2Fi22}(b).

(c) The natural $F$-vector space $V = F^{3588}$ is an $F\mathfrak
H_1$-module because $\mathfrak H_1 = \langle y_1, q_1 \rangle$.
Applying the Meataxe algorithm it follows that $V$ has three
composition factors $V_{78}$, $V_{1430}$ and $V_{2080}$ of
dimensions $78$, $1430$ and $2080$, respectively. Since all three
dimensions are divisible by $13$ and $13$ divides $|\mathfrak H_1|
= |2\Fi_{22}|$ to the first power only all three simple
composition factors of $V$ are projective $F\mathfrak H$-modules
by Theorems 3.12.2 and 3.12.4 of \cite{michler}. Hence $V$ is
isomorphic to their direct sum. Thus (c) holds.

(d) The group $A_1$ of Lemma \ref{l. aut2Fi22} has a unique
irreducible character $\chi_{\bf 17}$ of degree $4160$ Table
\ref{Fi_24 ct A_1}. Clifford's Theorem 2.6.15 of \cite{michler}
asserts that its restriction to $H_1$ is a sum of two inequivalent
irreducible characters $\nu$ and $\nu^t$ of degree $2080$. In
particular, the induced $FA_1$-module $V_{2080}^{A_1}$ is the
reduction modulo $13$ of a lattice which affords the irreducible
character $\chi_{\bf 17}$ of $A_1$. It corresponds to the
irreducible representation $\lambda : A_1 \rightarrow
\GL_{4160}(13)$ of $13$-defect zero defined in statement (d). It
is well defined because Lemma \ref{l. aut2Fi22}(b) implies that
$a^t = g^{-1}$, $b^t = f^{-1}$, $c^t = e^{-1}$,  $d^t = d^{-1}$,
$h^t = h^{-1}$ and $i^t = i^{-1}$. Therefore $\phi(y)_{2080}$ and
$\phi(q)_{2080}$ are well defined by (a) and (c).

(e) Using the MAGMA command $\verb"LowIndexSubgroups(A_1, m)"$ we
searched for conjugacy classes of subgroups $U$ of index $|A_1 :
U| = m$ such that $\chi_k$ is an irreducible constituent of the
permutation character $1_{U}^{A_1}$ for $k \in \{3, 4, 11, 12\}$.
Thus we found a subgroup $U$ of index $m = 2358720$ such that its
permutation character contains all four irreducible characters
$\chi_k$. Its generators have been obtained by means of Kim's
program $\verb"GetShortGens(A_1, U)"$.

Statement (f) can be verified by means of the character table of
$A_1$.
\end{proof}

\begin{proposition}\label{prop. matrixFi_24}
Keep the notation of Lemma \ref{l. aut2Fi22} and Propositions
\ref{prop. rep3588Fi23}, \ref{prop. rep5083Fi23}. Let $A_1
\leftarrow H_1 \rightarrow G_1$ be the amalgam constructed in
Lemma \ref{l. gensFi24}, where $G_1 \cong \Fi_{23}$. Let $\sigma:
H_1 \rightarrow A_1$ denote its monomorphism of $H_1$ into $G_1$.
Let $Y = \GL_{8671}(13)$. Let $\mathfrak \sigma(y) =
diag(\mathfrak y_{3588}, \mathfrak y_{5083})$, $\mathfrak
\sigma(q) = diag(\mathfrak q_{3588}, \mathfrak q_{5083})$,
$\mathfrak w_1 = diag(\mathfrak w_{3588}, \mathfrak w_{5083})$ in
$Y$. Let $\mathfrak H_1 = \langle \mathfrak \sigma(y), \mathfrak
\sigma(q)\rangle$ and $\mathfrak G_1 = \langle \mathfrak H_1,
\mathfrak w_1\rangle$. Keep the notation of Table \ref{Fi_24 ct
A_1}, Table 6.6.3 of \cite{michler2} and of the character table of
$\Fi_{23}$, see \cite{atlas}, its p. 178 - 179.

Then the following statements hold:

\begin{enumerate}
\item[\rm(a)] There is a compatible pair of characters
$$
(\chi, \quad \tau) = (\chi_{3} + \chi_{12} + \chi_{13} + \chi_{\bf
17},\quad \tau_{\bf 3} + \tau_{\bf 4}) \in mf
\mbox{char}_{\mathbb{C}}(A_1) \times mf
\mbox{char}_{\mathbb{C}}(G_1)
$$
of degree $8671$ of the groups $A_1 = \langle H_1, t \rangle$ and
$G_1 = \langle H_1 , e_1 \rangle$ with common restriction
$$
\tau_{|H_1} = \chi_{|H_1} = \delta_{2} + \delta_{6} + \delta_{\bf
7} + \delta_{\bf 8} + \delta_{9} \in mf
\mbox{char}_{\mathbb{C}}(H_1) ,
$$
where irreducible characters with bold face indices denote
faithful irreducible characters.

\item[\rm(b)] Let $\mathfrak V$ and $\mathfrak W$ be the up to
isomorphism uniquely determined faithful semi-simple
multiplicity-free $8671$-dimensional modules of $A_1$ and $G_1$
over $F = \GF(13)$ corresponding to the compatible pair $\chi,
\tau $, respectively.

Let $\kappa_\mathfrak V : H \rightarrow \GL_{8671}(13)$ and
$\kappa_\mathfrak W : E \rightarrow \GL_{8671}(13)$ be the
representations of $A_1$ and $G_1$ afforded by the modules
$\mathfrak V$ and $\mathfrak W$, respectively.

Let $\mathfrak q = \kappa_\mathfrak V(q)$, $\mathfrak y =
\kappa_\mathfrak V(y), \mathfrak t = \kappa_\mathfrak V(t)$ in $
\kappa_\mathfrak V(A_1) \le \GL_{8671}(13)$. Then the following
assertions hold:

\begin{enumerate}
\item[\rm(1)] $\mathfrak V_{|\mathfrak H_1} \cong \mathfrak
W_{|\mathfrak H_1}$, and there is a transformation matrix
$\mathcal T \in \GL_{8671}(13)$ such that
$$
\mathfrak q = \mathcal T^{-1} \kappa_\mathfrak W (\sigma(q))
\mathcal T, \mathfrak y = \mathcal T^{-1} \kappa_\mathfrak
W(\sigma(y))) \mathcal T.
$$
Let $\mathfrak w = \mathcal T^{-1} \kappa_\mathfrak W(w_1)
\mathcal T \in \GL_{8671}(13)$.

\item[\rm(2)] The subgroup $\mathfrak G = \langle \mathfrak q,
\mathfrak y,\mathfrak t, \mathfrak w \rangle$ of $Y$ is the
uniquely determined irreducible representation of the free product
$P = G_1*_{H_1}A_1$ of $G_1$ and $A_1$ with amalgamated subgroup
$H_1$ corresponding to the compatible pair (3) of Lemma \ref{l.
gensFi24}(e). Its four generating matrices can be
downloaded from the first author's website\\
$\verb"http://www.math.yale.edu/~hk47/Fi24/index.html"$.
\end{enumerate}
\end{enumerate}
\end{proposition}

\begin{proof}
(a) By Lemma \ref{l. gensFi24}(e) the amalgam $A_1 \leftarrow H_1
\rightarrow G_1$ has $8$ compatible pairs of degree $8671$. We
constructed the corresponding semi-simple representations of $P =
G_1*_{H_1}A_1$ for each of them. But we give a proof only for that
pair (3) of Lemma \ref{l. gensFi24}(e). It belongs to a group of a
suitable order.

(b) By Propositions \ref{prop. rep3588Fi23} and \ref{prop.
rep5083Fi23} the semi-simple $FG_1$-module $\mathfrak W =
\mathfrak W_3 \oplus \mathfrak W_4$ of dimension $8671$
corresponding to the multiplicity free character $\tau_{\bf 3} +
\tau_{\bf4}$ is described by the three matrices $\mathfrak y =
\sigma(y)$, $\mathfrak q = \sigma(q)$, $\mathfrak w$ of $Y$. The
semi-simple $FA_1$-module $\mathfrak V$ of the same dimension
corresponding to the multiplicity free character $\chi_{3} +
\chi_{12} + \chi_{14} + \chi_{\bf 17}$ of $A_1$ is defined by
three blocked diagonal matrices $\mathfrak q = diag(\mathfrak q_3,
\mathfrak q_{12}, \mathfrak q_{13}, \mathfrak q_{17})$, $\mathfrak
y = diag(\mathfrak y_3, \mathfrak y_{12}, \mathfrak y_{13},
\mathfrak y_{17})$ and $\mathfrak t = diag(\mathfrak t_3,
\mathfrak t_{12}, \mathfrak t_{13}, \mathfrak t_{17})$ whose
entries can be calculated by means of Proposition \ref{prop.
repsA1} as follows.

Assertion (d) of Proposition \ref{prop. repsA1} states that
$\mathfrak q_{17} = \lambda(q)$, $\mathfrak y_{17} = \lambda(y)$,
$\mathfrak t_{17} = \lambda(t)$ in $\GL_{4160}(13)$. By
Proposition \ref{prop. repsA1}(e) the irreducible characters
$\chi_3$ and $\chi_{12}$ are constituents of the permutation
character $1_U^{A_1}$ of a well determined subgroup $U$ of $A_1$
of index $2358720$. Let $P_U$ be the corresponding permutation
module over $F = GF(13)$. Using a stand alone program of the first
author which is based on Algorithm 5.7.1 of \cite{michler} we
calculated the primitive idempotents $E_{\chi_3}$ and
$E_{\chi_{12}}$ of the endomorphism ring $End_{FA_1}(P_U)$ of
$P_U$. Since the irreducible characters $\chi_3$ and $\chi_{12}$
are of $13$-defect zero the $FA_1$-modules $\mathfrak V_3 =
E_{\chi_3}P_U$ and $\mathfrak V_{12} = E_{\chi_{12}}P_U$ are
irreducible by Theorems 3.12.2 and 3.12.4 of \cite{michler}. After
fixing a basis in each of them the actions of the generators $q$,
$y$ and $t$ of $A_1$ on $P_U$ induce the matrices $\mathfrak q_3$,
$\mathfrak q_{12}$, $\mathfrak y_3$, $\mathfrak y_{12}$ and
$\mathfrak t_3$, $\mathfrak t_{12}$, respectively.

The tensor product $\chi_3 \otimes \chi_3$ contains a copy of the
irreducible character $\chi_{13}$ by Proposition \ref{prop.
repsA1}(f). Since $\mathfrak V_3$ has dimension $78$ the Meataxe
algorithm implemented in MAGMA can be applied to the tensor
product $\mathfrak V_3 \otimes \mathfrak V_3$. This application
provides the three $3003 \times 3003$ matrices $\mathfrak q_{13}$,
$\mathfrak y_{13}$ and $\mathfrak t_{13}$ corresponding to the
irreducible $FA_1$-module $\mathfrak V_{13}$. Hence $\mathfrak V =
\mathfrak V_3 \oplus \mathfrak V_{12} \oplus \mathfrak V_{13}
\oplus \mathfrak V_{17}$.

By (a) the restrictions ${\mathfrak V_3}_{|\mathfrak H_1}$,
${\mathfrak V_{12}}_{|\mathfrak H_1}$ and ${\mathfrak
V_{13}}_{|\mathfrak H_1}$ to $\mathfrak H_1$ are irreducible. By
the proof of Proposition \ref{prop. repsA1}(d) we know that
\begin{eqnarray*}
&&{\mathfrak V_{17}}_{|\mathfrak H_1} = \mathfrak V_{2080} \oplus \mathfrak V_{2080}\otimes t,\\
&&{\mathfrak V_3}_{|\mathfrak H_1} \oplus {\mathfrak V_{12}}_{|\mathfrak H_1} \oplus \mathfrak V_{2080} \cong {\mathfrak W_3}_{|\mathfrak H_1},\\
&&{\mathfrak V_{13}}_{|\mathfrak H_1} \oplus \mathfrak V_{2080}\otimes t \cong {\mathfrak W_4}_{|\mathfrak H_1}.\\
\end{eqnarray*}
Let $X_3 = \GL_{3588}(13)$, $X_4 = \GL_{5083}(13)$, $V_3 =
{\mathfrak V_3}_{|\mathfrak H_1} \oplus {\mathfrak
V_{12}}_{|\mathfrak H_1}$, $V_4 = {\mathfrak V_{13}}_{|\mathfrak
H_1} \oplus \mathfrak V_{2080} \otimes t$, $W_3 = {\mathfrak
W_3}_{|\mathfrak H_1}$ and $W_4 = {\mathfrak W_4}_{|\mathfrak
H_1}$. By the proof of Proposition \ref{prop. repsA1}(d) $V_3
\cong W_3$ and $V_4 \cong W_4$ as $F\mathfrak H_1$-modules of
respective dimensions $3588$ and $5083$. Applying now Parker's
isomorphism test of Proposition 6.1.6 of \cite{michler} by means
of the MAGMA command
$$\verb"IsIsomorphic(GModule(sub<X_i|V_i(h),V_i(y)>),GModule(sub<X_i|W_i(h),W_i(y)>))",$$
$i \in \{3,4\}$, one obtains the transformation matrices $\mathcal
T_3 \in X_3$ and $\mathcal T_4 \in X_4$ such that $\mathcal T =
diag(\mathcal T_3, \mathcal T_4) \in Y$ satisfies $\mathfrak q =
\kappa_{\mathfrak W}(\sigma(q))^{\mathcal T}$ and $\mathfrak y =
\kappa_{\mathfrak W}(\sigma(y))^{\mathcal T}$.

Let $\mathfrak w = \mathcal T^{-1} \kappa_\mathfrak W(w_1)
\mathcal T \in \GL_{8671}(13)$. Corollary 7.2.4 of \cite{michler}
asserts that the matrix group $\mathfrak G = \langle \mathfrak q,
\mathfrak y, \mathfrak t, \mathfrak w \rangle$ is uniquely
determined by the compatible given in (a).
\end{proof}

\begin{remark}\label{r. 8pairs}
Using Proposition \ref{prop. repsA1} we constructed for each of the eight compatible
pairs $(k)$ of Lemma \ref{l. gensFi24}(e) a matrix $\mathfrak t_k$
for the new generator $t$ of $A_1 = \langle q, y, w, t\rangle$
using the methods of the proof of Proposition \ref{prop.
matrixFi_24}. Thus we obtained $8$ subgroups $\mathfrak G_k =
\langle \mathfrak q, \mathfrak y, \mathfrak t_k, \mathfrak
w\rangle$ of $\GL_{8671}(13)$. In each of them we tried to
calculate the orders of the following products of the generators:

 $$
 \begin{array}{|c||c|c|c|c|c|
}\hline
 \mbox{ Group Name} &
\mathfrak w\mathfrak t_k & \mathfrak y\mathfrak w\mathfrak t_k &
\mathfrak q\mathfrak w\mathfrak t_k & \mathfrak y\mathfrak w
\mathfrak t_k\mathfrak q & \mathfrak y\mathfrak t_k\mathfrak
w\mathfrak q\mathfrak t_k
\\\hline
\mathfrak G_1 & 12 & \mbox{fail} & - & - & -
\\\hline
\mathfrak G_2 & 24 & \mbox{fail} & - & - & -
\\\hline
\mathfrak G_3 & 4 & 24 & 24 & 21 & 33
\\\hline
\mathfrak G_4 & 8 & 24 & \mbox{fail} & - & -
\\\hline
\mathfrak G_5 & 8 & 24 & \mbox{fail} & - & -
\\\hline
\mathfrak G_6 & 4 & 24 & 24 & 42 & 66
\\\hline
\mathfrak G_7 & 24 & \mbox{fail} & - & - & -
\\\hline
\mathfrak G_8 & 12 & \mbox{fail} & - & - & -
\\\hline
 \end{array}
 $$

where ``fail" means the the product has an order which is greater
than $100$. The group $\mathfrak G$ of Proposition \ref{prop.
matrixFi_24} is $\mathfrak G_3$. Looking at the orders of many
random elements we saw that all such orders were bounded by $60$.
In particular, $\mathfrak p = \mathfrak y^2\cdot\mathfrak t_3\cdot
\mathfrak w \cdot \mathfrak q$ has order $29$.

Therefore we prove in the remainder of the article that $\mathfrak
G$ is isomorphic to Fischer's simple group $\Fi_{24}'$. Most
likely, $\mathfrak G_6$ is isomorphic to Fischer's {\em non
simple} group $\Fi_{24}$. In the other cases we were not able to
calculate the orders of non trivial words of the generators in
reasonable time.
\end{remark}

\section{Isomorphism between $\mathfrak G$ and Fischer's group
$\Fi_{24}'$}

In this section we construct an isomorphism between the matrix
group $\mathfrak G = \langle \mathfrak q, \mathfrak y, \mathfrak
t, \mathfrak w \rangle$ of Proposition \ref{prop. matrixFi_24} and
the commutator subgroup of the finitely presented group $G$ of
Hall and Soicher, see \cite{praeger}, p.111. Hence $\mathfrak G$
is isomorphic to Fischer's simple group $\Fi_{24}'$.

\begin{proposition}\label{prop. E(Fi_24)} Let $\mathfrak G =
\langle \mathfrak q, \mathfrak y, \mathfrak  t, \mathfrak w
\rangle$ be the subgroup of $\GL_{8671}(13)$ constructed in
Proposition \ref{prop. matrixFi_24}. Let $\mathfrak H_1 = \langle
\mathfrak y, \mathfrak q \rangle$, $\mathfrak A_1 = \langle
\mathfrak H_1, \mathfrak t\rangle$ and $\mathfrak G_1 = \langle
\mathfrak H_1, \mathfrak w\rangle$.

Let $E = \langle a,b,c,d,t,g,h,i,j,k,v_i \mid 1 \le i \le 11
\rangle$ be the non-split extension of the Mathieu group $\M_{24}$
by its simple $GF(2)$-module $V_2$ constructed in Lemma \ref{l.
M24-extensions}, and let $E_{23} = \langle
a,b,c,d,t,g,h,i,j\rangle$.

Let $\mathfrak x = [(\mathfrak y\mathfrak q^2\mathfrak y\mathfrak
q\mathfrak y\mathfrak q^2)^{11}(\mathfrak q^2\mathfrak
y^2\mathfrak q\mathfrak y\mathfrak q\mathfrak y)^{11}(\mathfrak
q\mathfrak y^2\mathfrak q\mathfrak y\mathfrak q\mathfrak
y\mathfrak q\mathfrak y\mathfrak q)^4]^{12}$, $\mathfrak u_1 =
(\mathfrak x\mathfrak y\mathfrak x)^7$, $\mathfrak u_2 =
(\mathfrak x\mathfrak y\mathfrak x\mathfrak y\mathfrak x\mathfrak
y\mathfrak x)^4$, $\mathfrak u_3 = (\mathfrak x\mathfrak
y\mathfrak x\mathfrak y^2\mathfrak x\mathfrak y^2\mathfrak
x\mathfrak y\mathfrak x\mathfrak y)^2$, $\mathfrak u_4 =
(\mathfrak x\mathfrak y\mathfrak x\mathfrak y^5\mathfrak
x\mathfrak y^4)^2$, and $\mathfrak u_5 = [\mathfrak y(\mathfrak
s\mathfrak w)^2\mathfrak y\mathfrak w]^7$.

Then the following assertions hold:
\begin{enumerate}
\item[\rm(a)] The subgroup $\mathfrak T_1 = \mathfrak u_i \mid 1
\le i \le 4\rangle$ of $\mathfrak D = \langle \mathfrak x,
\mathfrak y\rangle$ is a Sylow $2$-subgroup of $\mathfrak G_1 =
\langle \mathfrak q, \mathfrak y, \mathfrak w\rangle$ of order
$2^{18}$.

\item[\rm(b)] $\mathfrak T_1$ has a unique maximal elementary
abelian normal subgroup $\mathfrak B$ of order $2^{11}$. It is
generated by the $11$ involutions:
\begin{eqnarray*}
&&\mathfrak u_1,\quad \mathfrak u_2^2,\quad (\mathfrak u_1\mathfrak u_2)^2,\quad (\mathfrak u_1\mathfrak u_3)^2,\quad (\mathfrak u_1\mathfrak u_4)^2,\quad (\mathfrak u_2\mathfrak u_4)^4,\\
&&(\mathfrak u_1\mathfrak u_2\mathfrak u_3)^4,\quad (\mathfrak u_1\mathfrak u_2\mathfrak u_4)^4,\quad (\mathfrak u_1\mathfrak u_3\mathfrak u_4)^4,\quad (\mathfrak u_2^2\mathfrak u_3)^2,\quad (\mathfrak u_1\mathfrak u_2\mathfrak u_4\mathfrak u_2)^2.\\
\end{eqnarray*}

\item[\rm(c)] $\mathfrak s = (\mathfrak y^5\mathfrak t)^7$ is an
involution of $\mathfrak A_1$ such that $\mathfrak A_1 = \langle
\mathfrak H_1, \mathfrak s\rangle$, $\mathfrak T_1^{\mathfrak s} =
\mathfrak T_1$, $\mathfrak B^{\mathfrak s} = \mathfrak B$, and
$\mathfrak S = \langle \mathfrak T_1, \mathfrak s\rangle$ is a
Sylow $2$-subgroup of $\mathfrak A_1$ of order $2^{19}$.

\item[\rm(d)] $\mathfrak N_1 = N_{\mathfrak G_1}(\mathfrak B) =
\langle \mathfrak x,\mathfrak y,\mathfrak w \rangle$ is isomorphic
to a non-split extension of $\M_{23}$ by $V_2{|\M_{23}}$ and
$\mathfrak D_1 = N_{\mathfrak A_1}(\mathfrak B) = \langle
\mathfrak x, \mathfrak y, \mathfrak s\rangle$

\item[\rm(e)] There is an isomorphism $\rho$ between $\mathfrak
N_1$ and the subgroup $E_{23}$ of $E$ such that $\rho(\mathfrak y)
=(y_2^5y_3y_2y_3)^3(w_3w_1w_2w_1w_2w_1w_3w_1^2w_3w_2w_3w_2w_1w_3)^{20}$,
$\rho(\mathfrak x) = (x_1x_2x_4x_5x_4x_2x_5)^3$, $\rho(\mathfrak
w) = (e_2e_3e_2e_3^2)^7$, where
\begin{eqnarray*}
&&x_1 = (ij)^3,\quad x_2 = (gahigai)^2,\quad x_3 = (aghijagh)^4,\\
&&x_4 = (j h i g h a j i)^4,\quad x_5 = (a i g h j i g a i)^4,\\
&&y_1 = i,\quad y_2 = a g,\quad y_3 = (a h j)^3,\\
&&w_1 = (y_2y_3^2)^2,\quad w_2 = (y_1y_2y_1y_2y_3)^3,\quad w_3 = (y_1y_2y_3y_2^2)^3,\\
&&e_1 = (a g i j i h)^4,\quad e_2 = (a g^3 i h j)^7,\quad e_3 = g h g h i a i.\\
\end{eqnarray*}

\item[\rm(f)] There is an isomorphism $\mu$ between $\mathfrak D_1
= N_{\mathfrak A_1}(\mathfrak B)$ and the centralizer $C_{E}(u)$
of the involution $u = (\rho(\mathfrak x)\rho(\mathfrak y)^2)^7$
of $E$ such that $\mu(\mathfrak x) = \rho(\mathfrak x)$,
$\mu(\mathfrak y) = \rho(\mathfrak y)$ and $\mu(\mathfrak s) =
(m_1^4m_2m_1m_2)^2$, where $m_1 = a g a h j$, $m_2 = (i j h k
j)^2$, $m_3 = (a h j a g k)^5$.

\item[\rm(g)] The subgroup $\mathfrak E = \langle \mathfrak
x,\mathfrak y,\mathfrak w,\mathfrak s\rangle$ of $\mathfrak G$ has
a faithful permutation representation $P\mathfrak E$ of degree
$1518$ with stabilizer $\langle (\mathfrak y\mathfrak s)^7,
(\mathfrak w\mathfrak y\mathfrak s\mathfrak y)^3, (\mathfrak
s\mathfrak y^3)^2, (\mathfrak y^2\mathfrak w\mathfrak
y^2)^3\rangle$.

\item[\rm(h)] The groups $\mathfrak E$ and $E$ are isomorphic.

\item[\rm(i)] $\mathfrak z = (\mathfrak x\mathfrak y\mathfrak
w)^8$ is a $2$-central involution of $\mathfrak E$ with
centralizer $C_{\mathfrak E}(\mathfrak z)$ of order $2^{21} \cdot
3^3 \cdot 5$ generated by  the elements $\mathfrak r_1 =
(\mathfrak s\mathfrak y^3)^3$, $\mathfrak r_2 = (\mathfrak
y^2\mathfrak w\mathfrak y\mathfrak s)^6$, $\mathfrak r_3 =
(\mathfrak s\mathfrak y\mathfrak w\mathfrak y\mathfrak s)^2$, and
$\mathfrak r_4 = (\mathfrak s\mathfrak w\mathfrak y\mathfrak
s\mathfrak w)^6$ with respective orders $2$, $4$, $4$, and $2$.

\item[\rm(j)] $C_{\mathfrak G_1}(\mathfrak z)$ has order
$2^{18}\cdot 3^5\cdot 5$. It is generated by $\mathfrak f_1 =
\mathfrak x\mathfrak y\mathfrak w$, $\mathfrak f_2 = (\mathfrak
x\mathfrak y\mathfrak x\mathfrak w)^7$, $\mathfrak f_3 =
(\mathfrak x\mathfrak w\mathfrak y\mathfrak w\mathfrak y^2)^7$,
and $\mathfrak v = (\mathfrak w\mathfrak q\mathfrak w\mathfrak
y\mathfrak q\mathfrak y)^7$.

\item[\rm(k)] The subgroup $\langle \mathfrak v\mathfrak
f_2\mathfrak v, (\mathfrak f_1\mathfrak f_2\mathfrak f_1\mathfrak
v\mathfrak f_1)^4, (\mathfrak f_2\mathfrak f_1\mathfrak v\mathfrak
f_1\mathfrak v)^4\rangle$ of $C_{\mathfrak G_1}(\mathfrak z)$ has
index $512$. Furthermore, it does not contain $\mathfrak z$.

\item[\rm(l)] $\mathfrak B$ is the Fitting subgroup of $\mathfrak
E$. It is also the unique maximal elementary abelian normal
subgroup of the Sylow $2$-subgroup $\mathfrak S = \langle
\mathfrak u_i \mid 1 \le i \le 5\rangle$ of $\mathfrak E$
contained in $C_{\mathfrak E}(\mathfrak z)$.
\end{enumerate}
\end{proposition}

\begin{proof}
In order to simplify the notation of the proof we replace the
German letters by Roman letters. In particular, we let $ME =
\langle x,y,w,s\rangle$ be the subgroup $\mathfrak E$.

(a) Let $PA_1$ be the faithful permutation representation of $A_1$
of degree $56320$ constructed in Lemma \ref{l. aut2Fi22}(c). By
Lemma \ref{l. gensFi24} and Lemma \ref{l. reps13Fi23}(f) we know
that $H_1 = \langle y, q\rangle = \langle x, y, q\rangle$. Now
\cite{kim1} asserts that $D = \langle x,y \rangle$ has odd index
in $H_1$ and therefore in $G_1$. Thus $D$ contains a Sylow
$2$-subgroup of $G_1$. The given Sylow $2$-subgroup $T_1$ of $G_1$
and its generators $t_i$ have been found by using MAGMA, the
permutation representation $PA_1$, and the program
$\verb"GetShortGens(H_1,T_1)"$.

(b) Applying the MAGMA command
$$\verb"Subgroups(T_1: Al:=Normal, IsElementaryAbelian := true)"$$
we observed that $T_1$ has $44$ elementary abelian normal
subgroups. Exactly one of them is maximal and has order $2^{11}$.
It is denoted by $B$. Its given generators have been calculated by
means of the first author's program $\verb"GetShortGens(T_1, B)"$.

(c) Since $|A_1 : H_1| = 2$ a Sylow $2$-subgroup of $A_1$ has
order $2^{19}$. Let $W_1 = N_{A_1}(T_1)$. Applying $PA_1$ and the
MAGMA command

$$\verb"exists(r){x: x \in A_1| T_1^x = T_1 and x^2 eq 1 and x \notin H_1}"$$

we found the involution $s \in A_1$ of the statement satisfying $s
\notin H_1$. It satisfies the equation $T_1^s = T_1$. Hence $B^s =
B$ holds trivially by (b).

(d) By another application of $PA_1$ and MAGMA we verified that
$N_{A_1}(B) = \langle x,y,s\rangle$. Using the faithful
permutation representation $PG_1$ of degree $31671$ with
stabilizer $H_1$ of Kim's Theorem 6.3.1 of \cite{michler2} one
establishes that $N_1 = N_{G_1}(B) = \langle x,y,w\rangle$. Hence
$N_1$ is a non split extension of $\M_{23}$ by $B$, see Lemma
6.1.2 and Theorem 6.3.1 of \cite{michler2}.

(e) By Lemma \ref{l. M24-extensions}(e) $E = \langle
a,b,c,d,t,g,h,i,j,k\rangle$ has a the faithful permutation
representation $PE$ with stabilizer $U_3 = \langle g,h,i,(dg)^5,
(d h j k)^3, (i j k j)^2, (d h j i d g)^3 \rangle$. Lemma 8.2.2 of
\cite{michler} states that its subgroup $E_{23} = \langle
a,b,c,d,t,g,h,i,j\rangle$ has index $|E : E_{23}| = 24$. Applying
the command $\verb"IsIsomorphic(N_1,E_{23})"$ MAGMA establishes an
isomorphism $\rho : N_1 \rightarrow E_{23}$. The words of the
images $\rho(x)$, $\rho(y)$ and $\rho(w)$ of the generators $x$,
$y$ and $w$ of $N_1$ are constructed as follows.  Let $C =
C_{E_{23}}(\rho(x))$. Using $PE$ and MAGMA one sees that $|C| =
2^{14}$. The six generators $x_1 = (ij)^3$, $x_2 = (gahigai)^2$,
$x_3 = (aghijagh)^4$, $x_4 = (jhighaji)^4$, $x_5 = (aighjigai)^4$
of $C$ were obtained computationally by means of the program
$\verb"GetShortGens(E_{23},C)"$. Using the program
$\verb"LookupWord(C, \rho(x))"$ MAGMA returned $\rho(x) =
(x_1x_2x_4x_5x_4x_2x_5)^3$. The expressions for $\rho(y)$ and
$\rho(w)$ are obtained similarly.

(f) By Kim's Theorem 6.3.1 of \cite{michler2} we know that $z_1 =
(xy^2)^7$ is a $2$-central involution of $G_1$. Clearly $u =
(\rho(x)\rho(y)^2)^7$ is an involution of $E_{23}$. Applying MAGMA
and $PE$ the reader can verify that $C_u = C_{E_{23}}(u)$ has
order $2^{19}\cdot3^2\cdot5\cdot7\cdot11$. Similarly one observes
that $D_1 = N_{A_1}(B) = \langle x,y,s\rangle$ has the same order.
Using $PA_1$, $PE$ and the command $\verb"IsIsomorphic(D_1,C_u)"$
MAGMA establishes an isomorphism $\mu : D_1 \rightarrow C_u$.
Applying the MAGMA command
$\verb"IsConjugate(C_u,\mu(y),\rh(y))"$\\ we found an element $c_1
\in C_u$ such that $\mu(y)^{c_1} = \rho(y)$. Using the command
$$\verb"exists(c){q: q in C_{C_u}(\rho(y))| (\mu(x)^{c_1})^q = \rho(x)}"$$
one gets an element $c_2 \in C_u$ such that $\mu(x)^{c_1c_2} =
\rho(x)$. Hence $\mu' : D_1 \rightarrow C_u$ defined by $\mu'(d) =
\mu(d)^{c_1c_2}, d \in D_1,$ is an isomorphism between $D_1$ and
$C_u$ such that $\mu'(x) = \rho(x)$ and $\mu'(y) = \rho(y)$. It
has been checked that $C_{E}(\langle \rho(y),\rho(x)\rangle) =
\langle u \rangle$. Furthermore, $\mu'(s)$ has a centralizer
$C_{M_{23}}(\mu'(s))$ of order $2^{17}$ which is generated by the
three elements $m_1$, $m_2$ and $m_3$ of $M_{23}$ given in the
statement. Another application of the Lookup command yields the
word $\mu'(s) = (m_1^4m_2m_1m_2)^2$. Hence the map $\mu' : D_1
\rightarrow C_{E}(u)$ satisfies all conditions of (f).

(g) Using (e), (f), $PE$ and MAGMA it has been verified that\\
$E = \langle \rho(x), \rho(y), \rho(w), \mu(s)\rangle$. As $U_3 =
\langle g,h,i,(dg)^5, (d h j k)^3, (i j k j)^2, (d h j i d g)^3
\rangle$ is a stabilizer of $PE$ we apply the program
$\verb"GetShortGens(E,U_3)"$ w.r.t. the given generators of
$E$. MAGMA returns\\
$U_3 = \langle (\rho(y)\mu(s))^7, (\rho(w)\rho(y)\mu(s)\rho(y))^3,
(\mu(s)\rho(y)^3)^2, (\rho(y)^2\rho(w)\rho(y)^2)^3\rangle$.

Thus $MU = \langle (ys)^7, (wysy)^3, (sy^3)^2, (y^2wy^2)^3\rangle$
is a subgroup of $ME = \langle x,y,w,s\rangle$ which is isomorphic
to $U_3$. Let $V$ be the $8671$-dimensional vector space over $F =
GF(13)$. Using the Meataxe Algorithm implemented in MAGMA we see
that the restriction $V_{MU}$ of $V$ to the subgroup $MU$ has a
$7$-dimensional $FMU$-submodule $W$ which has a complement of
dimension $8664$. Applying now the algorithm described in Theorem
6.2.1 of \cite{michler} we obtain a faithful permutation
representation $PME$ of the matrix group $ME$ of degree $1518$
with stabilizer $MU$.

(h) Using $PE$, $PME$ and the isomorphism test
$\verb"IsIsomorphic(PE,PME)"$ MAGMA established that $ME \cong E$.

(i) By (d) and Table 6.5.1 of \cite{michler2} we know that $z =
(xyw)^8$ is an involution of $N_1 = N_{G_1}(B) \cong E_{23}$ with
centralizer $C_{N_1}(z)$ of order $2^{18}\cdot3^2\cdot5$.
Therefore we calculate $C_E(z)$ by means of $PE$ and MAGMA. It
follows that $|C_E(z)| = 2^{21}\cdot 3^3\cdot5$. Hence $E$ has a
Sylow $2$-subgroup $S_3$ of order $2^{21}$ with center $Z(S_3) =
\langle z \rangle$ by Table \ref{Fi_24 cc E}. The given generators
$r_i$ of $C_z = C_E(z)$ have been determined by means of MAGMA and
the program $\verb"GetShortGens(E, C_z)"$.

(j) Table 6.5.6 of \cite{michler2} implies that $|C_{G_1}(z)| =
2^{18}\cdot3^5\cdot5$ because $z = (xyw)^8 \in G_1 = \langle
x,y,w,q\rangle$. Using MAGMA and the faithful permutation
representation $PG_1$ of $G_1$ we found the involution $v =
(wqwyqy)^7$ such that $(z,v) = 1$, and $C_{G_1}(z) = \langle f_1,
f_2, f_3, v\rangle$ for the elements $f_i \in G_1$ given in the
statement.

(k) All assertions of the statement are easily checked by means of
MAGMA and the faithful permutation representation $PG_1$ of $G_1$.

(l) By (b) and (c) the elementary abelian subgroup $B$ is normal
in $ME$. Hence it is the Fitting subgroup $ME$ by (h) and Lemma
\ref{l. M24-extensions}. Using the faithful permutation
representation of $ME$ given in (g) the remaining assertions can
be verified by means of MAGMA.
\end{proof}

The following presentation of the $3$-transposition Fischer group
$P = \Fi_{24}$ is taken from \cite{praeger}, its p. 124. It is due
to J. Hall and L. S. Soicher \cite{hall}.

\begin{lemma}\label{l. nicepres.Fi_24}
Let $\mathfrak G = \langle \mathfrak y, \mathfrak q, \mathfrak t,
\mathfrak w\rangle$ be the subgroup of $\GL_{8671}(13)$
constructed in Proposition \ref{prop. matrixFi_24}. Let $\mathfrak
H_1 = \langle \mathfrak y, \mathfrak q \rangle$, $\mathfrak G_1 =
\langle \mathfrak H_1, \mathfrak w\rangle$ and $\mathfrak A_1 =
\langle \mathfrak H_1, \mathfrak t\rangle$. Let $\mathfrak s =
(\mathfrak y^5\mathfrak t)^7$ and $\mathfrak x = [(\mathfrak
y\mathfrak q^2\mathfrak y\mathfrak q\mathfrak y\mathfrak
q^2)^{11}(\mathfrak q^2\mathfrak y^2\mathfrak q\mathfrak
y\mathfrak q\mathfrak y)^{11}(\mathfrak q\mathfrak y^2\mathfrak
q\mathfrak y\mathfrak q\mathfrak y\mathfrak q\mathfrak y\mathfrak
q)^4]^{12}$. Let $\mathfrak E = \langle \mathfrak x, \mathfrak y,
\mathfrak w, \mathfrak s\rangle$

Let $P = \langle a,b,c,d,e,f,g,h,i,j,k,l \rangle $ be the finitely
generated group with the following set $\mathcal R(P)$ of defining
relations:
\begin{eqnarray*}
&&l^2 = k^2 = a^2 = b^2 = c^2 = d^2 = e^2 = f^2 = g^2 = j^2 = h^2 = i^2 = 1, \\
&&(l k)^3 = (k a)^3 = (a b)^3 = (b c)^3 = (c d)^3 = (d e)^3 = (e f)^3 = (f g)^3 = (g j)^3 = 1, \\
&&(l a)^2 = (l b)^2 = (l c)^2 = (l d)^2 = (l e)^2 = (l f)^2 = (l g)^2 = (l j)^2 = (l h)^2 = (l i)^2 = 1,\\
&&(k b)^2 = (k c)^2 = (k d)^2 = (k e)^2 = (k f)^2 = (k g)^2 = (k j)^2 = (k h)^2 = (k i)^2 = 1,\\
&&(a c)^2 = (a d)^2 = (a e)^2 = (a f)^2 = (a g)^2 = (a j)^2 = (a h)^2 = (a i)^2 = 1,\\
&&(b d)^2 = (b e)^2 = (b f)^2 = (b g)^2 = (b j)^2 = (b h)^2 = (b i)^2 = 1,\\
&&(c e)^2 = (c f)^2 = (c g)^2 = (c j)^2 = (c h)^2 = (c i)^2 = (d f)^2 = (d g)^2 = (d j)^2 = 1,\\
&&(e g)^2 = (e j)^2 = (e h)^2 = (e i)^2 = (f j)^2 = (f h)^2 = (f i)^2 = (g h)^2 = (g i)^2 = 1,\\
&&(j h)^2 = (j i)^2 = (d h)^3 = (h i)^3 = (d i)^2 = 1,\\
&&l = (a b c d e f h)^9, (d c b a k l d e f g j d h i)^{17} = 1.
\end{eqnarray*}

Then the following statements hold:
\begin{enumerate}
\item[\rm(a)] $P$ has a faithful permutation representation $PP$
of degree $306936$ with stabilizer $M = \langle
a,b,c,d,e,f,g,h,i,j,l \rangle$.

\item[\rm(b)] The commutator subgroup $G = P'$ is a finite simple
group of order $2^{21} \cdot 3^{16} \cdot 5^2 \cdot 7^3 \cdot 11
\cdot 13 \cdot 17 \cdot 23 \cdot 29$.

\item[\rm(c)] $G = \langle b_1, c_1, d_1, e_1, f_1, g_1, h_1, i_1,
j_1, k_1 \rangle$, where $b_1 = ab$, $c_1 = ac$, $d_1 = ad$, $e_1
= ae$, $f_1 = af$, $g_1 = ag$, $h_1 = ah$, $i_1 = ai$, $j_1 = aj$,
and $k_1 = ak$.

Furthermore, $G$ has a faithful permutation representation $PG$ of
degree $306936$ with stabilizer $M_1 = \langle b_1, c_1, d_1, e_1,
f_1, g_1, h_1, i_1, j_1\rangle$ and $M_1$ is a simple group of
order $2^{18} \cdot 3^{13} \cdot 5^2 \cdot 7 \cdot 11 \cdot 13
\cdot 17 \cdot 23$.

\item[\rm(d)] The centralizer $C_1 = C_G(c_1)$ of the involution
$c_1$ of $G$ has order $2^{19} \cdot 3^9 \cdot 5^2 \cdot 7 \cdot
11 \cdot 13$, and $C_1 = \langle e_1,g_1,i_1,j_1,m_1,n_1\rangle$,
where $m_1 = (c_1d_1e_1h_1)^3$ and $n_1 = (b_1c_1d_1k_1f_1)^3$ are
involutions.

Furthermore, $C_1$ has a faithful permutation representation
$PC_1$ of degree $56320$ with stabilizer $\langle
i_1,j_1,(e_1g_1n_1)^2, n_1m_1n_1\rangle$.

\item[\rm(e)] There is an isomorphism $\phi$ between $\mathfrak
A_1 = \langle \mathfrak y, \mathfrak q, \mathfrak t\rangle$ and
the finitely presented group $A_1 = \langle
a,b,c,d,e,f,g,h,i,z,t\rangle$ constructed in Lemma \ref{l.
aut2Fi22}(b) such that
\begin{eqnarray*}
&&y = \phi(\mathfrak y) = (S_1S_2S_4S_1S_4S_2S_4)^5(T_1T_3^2T_1T_3T_1T_2T_1T_3)^{20},\\
&&q = \phi(\mathfrak q) = [(p\cdot o\cdot j\cdot o\cdot k)^4\cdot(j\cdot k\cdot o\cdot k^2p)^4 ]^4,\\
&&t = \phi(\mathfrak t), \quad \mbox{where}\\
&&S_1 = (b\cdot a\cdot c\cdot b)\cdot (d\cdot e\cdot h\cdot d)\cdot(b\cdot a\cdot c\cdot b),\\
&&S_2 = [(b\cdot a\cdot c\cdot b)\cdot (d\cdot e\cdot h\cdot d)\cdot (f\cdot e\cdot g\cdot f)]^2,\\
&&S_3 = [(b\cdot a\cdot c\cdot b)\cdot(d\cdot e\cdot h\cdot d)\cdot (c\cdot d\cdot e\cdot h\cdot i)^4\cdot (b\cdot a\cdot c\cdot b)]^2,\\
&&S_4 = [(b\cdot a\cdot c\cdot b)\cdot (c\cdot d\cdot e\cdot h\cdot i)^4\cdot (d\cdot e\cdot h\cdot d)\cdot (f\cdot e\cdot g\cdot f)]^4,\\
&&T_1 = S_2S_4S_2^2,\quad T_2 = (S_4S_2S_1S_3S_4)^2,\quad T_3 = (S_4S_2S_4S_2S_1S_3)^2,\\
&&j = (c\cdot d\cdot e\cdot h)^3,\\
&&k = (c\cdot d\cdot e\cdot f\cdot g)^2,\\
&&l = (a\cdot b\cdot c\cdot d\cdot e\cdot h)^5,\\
&&o = (l\cdot b\cdot k\cdot j\cdot b\cdot i)^6,\\
&&p = (j\cdot k\cdot j\cdot l\cdot b\cdot i\cdot j\cdot i)^5.
\end{eqnarray*}

\item[\rm(f)] There is an isomorphism $$\rho: A_1 = \langle
a,b,c,d,e,f,g,h,i,z,t\rangle \rightarrow C_1 = \langle
e_1,g_1,i_1,j_1,m_1,n_1\rangle$$ such that $\rho(t) = (m_1n_1)^3$
and:
\begin{eqnarray*}
&&\rho(a) = [(i_1 m_1 n_1)^6\cdot (j_1 n_1 m_1 n_1)\cdot (g_1 n_1 e_1 g_1 n_1)^2\cdot (i_1 m_1 n_1)^6\\
&&\quad \quad \quad \quad \cdot (e_1 g_1 n_1 m_1 n_1 e_1)^2\cdot (g_1 n_1 e_1 g_1 n_1)^2]^7,\\
&&\rho(b) = [(i_1 n_1 g_1 j_1 m_1)^{12}\cdot (e_1 g_1 j_1 m_1 i_1 n_1)^4\cdot (i_1 n_1 g_1 j_1 m_1)^{12}\cdot (e_1 g_1 j_1 m_1 i_1 n_1)^{12}\\
&&\quad \quad \quad \quad \cdot (i_1 n_1 g_1 j_1 m_1)^{12}\cdot (e_1 g_1 j_1 m_1 i_1 n_1)^8]^{11},\\
&&\rho(c) = [(j_1 n_1 m_1 n_1)\cdot (e_1 n_1 g_1 n_1 e_1)\cdot (e_1 g_1 i_1 j_1 m_1 n_1)^{18}]^{11},\\
&&\rho(d) = [(j_1 m_1)\cdot (n_1 g_1 n_1)\cdot (e_1 i_1 j_1 n_1 g_1 m_1)^{18}\cdot (j_1 m_1)\cdot (e_1 i_1 j_1 n_1 g_1 m_1)^{18}\\
&&\quad \quad \quad \quad \cdot (j_1 m_1)\cdot (e_1 i_1 j_1 n_1 g_1 m_1)^{18}]^{11},\\
&&\rho(e) = (bc)^3\cdot c^t,\quad \rho(f) = (bc)^3\cdot b^t, \quad\rho(g) = a^t,\\
&&\rho(h) = [e_1\cdot i_1\cdot (m_1 i_1 n_1 m_1 n_1)^3\cdot (e_1 g_1 n_1 g_1 j_1 m_1)^6\cdot i_1\cdot (e_1 g_1 n_1 g_1 j_1 m_1)^6]^{11},\\
&&\rho(i) = [(m_1 i_1 m_1)\cdot (e_1 g_1 j_1 n_1 m_1)^{12}]^9.\\
\end{eqnarray*}

In particular, $\psi = \rho  \circ \phi: \mathfrak A_1 \rightarrow
C_1$ is an isomorphism such that\\ $y_1 = \psi(\mathfrak y) =
\rho(y)$, $q_1 = \psi(\mathfrak q) = \rho(q)$, and $t_1 =
\psi(\mathfrak t) = \rho(t)$ generate $C_1$.

\item[\rm(g)] In $C_1$ let $x_1 = \psi(\mathfrak x)$ and $s_1 =
\psi(\mathfrak s) = (y_1^5t_1)^7$.

Then $S_1 = \langle u_1, u_2, u_3,u_4, s_1\rangle$ is a Sylow
$2$-subgroup of $C_1$, where\\ $u_1 = (x_1y_1x_1)^7$, $u_2 =
(x_1y_1 x_1y_1x_1y_1x_1)^4$, $u_3 = (x_1y_1x_1y_1^2x_1y_1^2x_1 y_1
x_1y_1)^2$,\\ $u_4 = (x_1y_1x_1y_1^5x_1y_1^4)^2$.

Moreover, $S_1$ has a unique maximal elementary abelian normal
subgroup $B_1$. It is generated by the eleven involutions $u_1$,
$u_2^2$, $(u_1u_2)^2$, $(u_1u_3)^2$, $(u_1u_4)^2$, $(u_2u_4)^4$,
$(u_1u_2u_3)^4$, $(u_1u_2u_4)^4$, $(u_1u_3u_4)^4$, $(u_2^2u_3)^2$,
$(u_1u_2u_4u_2)^2$.

\item[\rm(h)] $N_1 = N_G(B_1) = \langle x_1,y_1,s_1,o_1\rangle$,
where $o_1 = [d_1 \cdot y_1 \cdot s_1 \cdot (y_1)^2 \cdot
s_1)]^{10}$ and $|N_1| = 2^{21}\cdot 3^3\cdot 5\cdot 7\cdot
11\cdot 23$.

\item[\rm(i)] There is an isomorphism $\mu: \mathfrak E
\rightarrow N_1$ such that $\mu(\mathfrak x) = \psi(\mathfrak x) =
x_1$, $\mu(\mathfrak y) = \psi(\mathfrak y) = y_1$, $\mu(\mathfrak
s) = \psi(\mathfrak t) = s_1$ and
$$\mu(w) = w_1 = [(s_1y_1x_1y_1)^2 \cdot (x_1y_1o_1y_1s_1y_1^2)^5]^7.$$

\item[\rm(j)] $G = \langle C_1, N_1\rangle = \langle q_1, y_1,
w_1, t_1 \rangle$.

\item[\rm(k)] The map $\kappa: G \rightarrow \mathfrak G$ given by
$\kappa(q_1) = \mathfrak q$, $\kappa(y_1) = \mathfrak y$,
$\kappa(w_1) = \mathfrak w$ and $\kappa(s_1) = \mathfrak s$ is a
group isomorphism.

In particular, $\mathfrak G$ is a simple group of order $2^{21}
\cdot 3^{16} \cdot 5^2 \cdot 7^3 \cdot 11 \cdot 13 \cdot 17 \cdot
23 \cdot 29$.
\end{enumerate}
\end{lemma}

\begin{proof}
(a) This statement has been verified by running the Todd-Coxeter
Algorithm $\verb"CosetAction(P,M)"$ built into MAGMA.

(b) Using (a) and MAGMA it has been checked that $G = P'$ is a
simple group of the stated order.

(c) Using then the program $\verb"GetShortGens(P,G)"$ we found the
given generators of $G$. Using the faithful permutation
representation $PP$ of (a) and MAGMA it has been checked that $M_1
= \langle b_1, c_1, d_1, e_1, f_1, g_1, h_1, i_1, j_1 \rangle$ is
a simple group of order $2^{18} \cdot 3^{13} \cdot 5^2 \cdot 7
\cdot 11 \cdot 13 \cdot 17 \cdot 23$. It is the stabilizer of the
faithful permutation representation $PG$ of $G$ of degree $306936$
as has been checked by means of MAGMA and the command
$\verb"CosetAction(G,M_1)"$.

(d) The centralizer $C_1 = C_G(c_1)$ of the involution $c_1$ and
$|C_1|$ have been determined by means of the permutation
representation $PG$ of $G$ and MAGMA. Its given generators have
been found by MAGMA using the program
$\verb"GetShortGens(G,C_1)"$. Applying the MAGMA command
$\verb"DegreeReduction"$ we get a faithful permutation
representation $PC_1$ of $C_1$ having degree $56320$. Its
stabilizer $U_1$ has been found by means of the MAGMA command
$\verb"BasicStabilizer(~,2)"$. Its given generators were gotten by
another application of $\verb"GetShortGens(C_1,U_1)"$.

(e) This statement follows immediately from Lemma \ref{l.
aut2Fi22}(b), Proposition 6.2.3 of \cite{michler2}, and
Proposition \ref{prop. matrixFi_24}.

(f)  By Lemma \ref{l. aut2Fi22}(c) the finitely presented group
$A_1 = \langle a,b,c,d,e,f,g,h,i,z,t\rangle$ has a faithful
permutation representation $PA_1$ of degree $56320$. Let $PC_1$ be
the faithful permutation representation of $C_1$ constructed in
(d). A successful application of the MAGMA command
$\verb"IsIsomorphic(PA_1,PC_1)"$ provides an isomorphism $\eta:
\mathfrak A_1 \rightarrow C_1$. Now $C_{C_1}(\eta(t))$ and
$C_{C_1}((m_1 n_1)^3)$ have the same order. Using the MAGMA
command $\verb"exists(w){x : x in C_1| \eta(t)^x eq (m_1n_1)^3}"$
we found an element $a \in C_1$ such that $\eta(t)^a =
(m_1n_1)^3$. Let $\alpha$ be the inner automorphism of $C_1$
induced by $a$. Then the map $\rho = \alpha \circ \eta$ is an
isomorphism from $A_1$ onto $C_1$ such that $\rho(t) =
(m_1n_1)^3$. The centralizers of the images of $\rho(a), \rho(b),
\ldots, \rho(i)$ in $C_1$ all have the same order $2^{17}\cdot
3^6\cdot 5\cdot 7\cdot 11$. Their generators and their words in
them are obtained computationally using the programs
$\verb"GetShortGens"$ and $\verb"LookupWord"$, respectively, as
follows:
\begin{enumerate}
\item $C_{C_1}(\rho(a))$ is generated by $a_1 = (i_1 m_1 n_1)^6$,
$a_2 = j_1 n_1 m_1 n_1$,\\
$a_3 = (g_1 n_1 e_1 g_1 n_1)^2$, $a_4 = (e_1g_1n_1m_1n_1e_1)^2$,
and $\rho(a) = (a_1 a_2 a_3 a_1 a_4 a_3)^7$.

\item $C_{C_1}(\rho(b))$ is generated by $b_1 = (i_1 n_1 g_1 j_1
m_1)^{12}$, $b_2 = (e_1 g_1 j_1 m_1 i_1n_1)^4$,\\ and $\rho(b) =
(b_1 b_2 b_1 b_2^3 b_1 b_2^2)^{11}$.

\item $C_{C_1}(\rho(c))$ is generated by $c_1 = j_1 n_1 m_1n_1$,
$c_2 = e_1 n_1 g_1 n_1 e_1$,\\ $c_3 = (e_1g_1 i_1 j_1 m_1
n_1)^{18}$, and $\rho(c) = (c_1 c_2 c_3)^{11}$.

\item $C_{C_1}(\rho(d))$ is generated by $d_1 = j_1 m_1$, $d_2 =
n_1 g_1 n_1$, $d_3 = (e_1 i_1 j_1n_1 g_1 m_1)^{18}$,\\
and $\rho(d) = (d_1 d_2 d_3 d_1 d_3 d_1 d_3)^{11}$.

\item $C_{C_1}(\rho(h))$ is generated by $h_1 = e_1$,  $h_2 =
i_1$, $h_3 = (m_1 i_1 n_1 m_1 n_1)^3$,\\ $h_4 = (e_1 g_1 n_1 g_1
j_1 m_1)^6$, and $\rho(h) = (h_1 h_2 h_3 h_4 h_2 h_4)^{11}$.

\item $C_{C_1}(\rho(i))$ is generated by $i_1 = m_1 i_1 m_1$, $i_2
= j_1 n_1 m_1 n_1$, $i_3 = (e_1 g_1 j_1n_1m_1)^6$,\\ and $\rho(i)
= (i_1i_3^2)^9$.
\end{enumerate}

The given words for the images $\rho(e)$, $\rho(f)$ and $\rho(g)$
can now be calculated from these images and the relations of Lemma
\ref{l. aut2Fi22}(b).

Statement (e) implies that the composition of the isomorphisms
$\phi: \mathfrak A_1 \rightarrow A_1$ and $\rho: A_1 \rightarrow
C_1$ is an isomorphism $\psi: \mathfrak A_1 \rightarrow C_1$.
Furthermore, the given images $\psi(\mathfrak q)$, $\psi(\mathfrak
q)$ and $\psi(\mathfrak q)$ in $C_1$ of the $3$ generators of
$\mathfrak A_1$ are well defined.

(g) Since $\mathfrak x \in \mathfrak A_1$ its image
$$\psi(\mathfrak x) = [(y_1q_1^2y_1q_1y_1q_1^2)^{11}(q_1^2y_1^2q_1y_1q_1y_1)^{11}(q_1y_1^2q_1y_1q_1y_1q_1y_1q_1)^4]^{12} = x_1$$\\
is a well defined element of $C_1$.

Let $\mathfrak u_i$ be the generators of the Sylow $2$-subgroup
$\mathfrak T_1$ of the subgroup $\mathfrak G_1 = \langle \mathfrak
q, \mathfrak y, \mathfrak w\rangle$ of $\mathfrak G$ given in
Proposition \ref{prop. E(Fi_24)}(a). Hence $\mathfrak T_1$ is a
subgroup of $\mathfrak A_1$ because its generators $\mathfrak u_i$
are words in $\mathfrak x$ and $\mathfrak y$ by Proposition
\ref{prop. E(Fi_24)}. Therefore their images $u_i = \psi(\mathfrak
u_i)$, $1 \le i \le 4$, are well defined. So is $s_1 =
\psi(\mathfrak s) = (y_1^5t_1)^7$. Hence $S_1 = \psi(\mathfrak S)
= \langle u_1, u_2, u_3, u_4, s_1\rangle$ is a Sylow $2$-subgroup
of $C_1$ by Proposition \ref{prop. E(Fi_24)}(c)and (f).

Let $B_1 \psi(\mathfrak B$ where $\mathfrak B$ is the unique
maximal elementary abelian normal subgroup of $\mathfrak S$
defined in Proposition \ref{prop. E(Fi_24)}(b and (d)). Then $B_1$
is generated by the $11$ involutions given in the statement.

(h) Let $N_1 = N_{G}(B_1)$. Using the faithful permutation
representation $PG$ of $G$ and MAGMA the reader can easily check
that $|N_1| = 2^{21} \cdot 3^3 \cdot 5 \cdot 7 \cdot 11 \cdot 23$
and that $N_1$ is generated by $x_1$, $y_1$, $s_1$ and the element
$o_1 \in G$ given in the statement.

(i) By Proposition \ref{prop. E(Fi_24)}(g) the subgroup $\mathfrak
E$ of $\mathfrak G$ has a faithful permutation representation
$P\mathfrak E$ of degree $1518$ with stabilizer $\langle
(\mathfrak y\mathfrak s)^7, (\mathfrak w\mathfrak y\mathfrak
s\mathfrak y)^3, (\mathfrak s\mathfrak y^3)^2, (\mathfrak
y^2\mathfrak w\mathfrak y^2)^3\rangle$. Let $PN_1$ be the
reduction of $PG$ to $N_1 = N_G(B_1)$. A successful application of
the MAGMA command $\verb"IsIsomorphic(PE,PN_1)"$ provides an
isomorphism $\tau: \mathfrak E \rightarrow N_1$. As in the proof
of (f) we find an element $a \in N_1$ such that $(\tau(\mathfrak
y))^a = y_1 = \psi(\mathfrak y)$. Using MAGMA again we verified
that $C_{N_1}(y_1)$ has order $56$. Searching through its elements
we find an element $b \in C_{N_1}(y_1)$ such that $(\tau(\mathfrak
x))^{ab} = \psi(\mathfrak x) = x_1$ and $(\tau(\mathfrak s))^{ab}
= \psi(\mathfrak s) = s_1$. Let $\beta$ denote the inner
automorphism of $N_1$ induced by conjugation with $ab$. Then the
map $\mu = \beta \circ \nu$ is an isomorphism from $\mathfrak E$
onto $N_1 = N_G(B_1)$ such that $\mu(\mathfrak x) = x_1$,
$\mu(\mathfrak y) = y_1$ and $\mu(\mathfrak s) = s_1$. Let $w_1 =
\mu(\mathfrak w)$. Another application of MAGMA and $PN_1$ yields
that $N_1 = \langle x_1, y_1, s_1, w_1\rangle$.

The word for $w_1$ in the generators of $G$ is obtained as
follows. Let $C_2 = C_{N_1}(w_1)$. Using the generators of $N_1$
given in (h), MAGMA, and the program\\
$\verb"GetShortGens(N_1,C_2)"$ we see that $C_2 = \langle v_i \mid
1 \le i \le 3\rangle$, where $v_1 = (s_1y_1x_1y_1)^2$, $v_2 =
(y_1o_1x_1s_1o_1)^4$, $v_3 = (x_1y_1o_1y_1s_1y_1^2)^5$. Applying
then the command\\
$\verb"LookupWord(C_2, w_1)"$ MAGMA provides the solution $w_1 =
(v_1v_3)^7$.

(j) Using the faithful permutation representation $PG$ of $G$ and
MAGMA the reader easily verifies that $G = \langle C_1,
N_1\rangle$. Hence (f) and (i) imply that $G = \langle
q_1,y_1,t_1,x_1,s_1,w_1\rangle = \langle q_1,y_1,t_1,w_1\rangle$.

(k) In $\mathfrak G$ let $\mathfrak E_{23} = \langle \mathfrak x,
\mathfrak y, \mathfrak w\rangle$ and $\mathfrak H_1 = \langle
\mathfrak q, \mathfrak y\rangle$. Then Proposition \ref{prop.
matrixFi_24} and Kim's Theorem 6.3.1 of \cite{michler2} imply that
$\mathfrak G_1 = \langle \mathfrak H_1, \mathfrak E_{23}\rangle$
is a simple subgroup of $\mathfrak G$ such that $\mathfrak D_1 =
\mathfrak E_{23} \cap \mathfrak H_1 = \langle \mathfrak x,
\mathfrak y\rangle$ and $\mathfrak H_1 = C_{\mathfrak
G_1}(\mathfrak z_1)$, where $\mathfrak z_1 = (\mathfrak x\mathfrak
y^2)^7$ is a $2$-central involution of $\mathfrak G_1$.

Let $E_{23} = \mu(\mathfrak E_{23})$, $H_1 = \psi(\mathfrak H_1)$
and $G_1 =\langle E_{23}, H_1\rangle$ in $G$. As $\mu$ and $\psi$
agree on $\mathfrak H_1$ by (j) we have $D_1 = \langle x_1, y_1
\rangle =  E_{23} \cap H_1$. Using $PG$ and MAGMA it has been
checked that $G_1$ is a simple group of order $2^{18}\cdot
3^{13}\cdot 5^2\cdot7\cdot11\cdot13\cdot17\cdot23$. Furthermore,
$C_{G_1}(z_1) = H_1$, where $z_1 = (x_1y_1^2)^7 = \psi(\mathfrak
z_1)$ is a $2$-central involution of $G_1$. Thus Kim's Theorem
6.3.1 of \cite{michler2} implies that $G_1 \cong \Fi_{23}$. In
particular, the map $\kappa_1: G_1 \rightarrow \mathfrak G_1$
defined by $\kappa_1(q_1) = \mathfrak q$, $\kappa_1(y_1) =
\mathfrak x$, and $\kappa_1(w_1) = \mathfrak w$ is an isomorphism.

By Proposition \ref{prop. matrixFi_24} and Lemma \ref{l.
gensFi24}(e) the amalgam $\mathfrak A_1 \leftarrow \mathfrak H_1
\rightarrow \mathfrak G_1$ has Goldschmidt index $1$. Therefore by
Corollary 7.1.9 of \cite{michler}, (f) and (i) imply that the free
products ${\mathfrak A_1}*_{\mathfrak H_1}{\mathfrak G_1}$ and
$C_1*_{H_1}G_1$ with respective amalgamated subgroups $\mathfrak
H_1$ and $H_1$ are isomorphic. Thus (j) and Proposition \ref{prop.
matrixFi_24} imply that the map $\kappa: G \rightarrow \mathfrak
G$ defined by $\kappa(q_1) = \mathfrak q$, $\kappa(y_1) =
\mathfrak y$, $\kappa(w_1) = \mathfrak w$ and $\kappa(t_1) =
\mathfrak t$ is an irreducible $8671$-dimensional representation
of the group $G$ over $GF(13)$. Hence  $\kappa: G \rightarrow
\mathfrak G$ is an isomorphism because $\kappa$ is surjective and
$G$ is simple by (b). This completes the proof.
\end{proof}

\begin{theorem}\label{thm. existenceFi_24}
Let $\mathfrak G = \langle \mathfrak q, \mathfrak y,\mathfrak t,
\mathfrak w \rangle$ be the subgroup of $\GL_{8671}(13)$
constructed in Proposition \ref{prop. matrixFi_24}. Then the
following statements hold:

\begin{enumerate}
\item[\rm(a)] $\mathfrak G = \langle \mathfrak q, \mathfrak
y,\mathfrak t, \mathfrak w \rangle$ has a faithful permutation
representation of degree $306936$ with stabilizer $\mathfrak G_1 =
\langle \mathfrak q, \mathfrak y, \mathfrak w \rangle$.

\item[\rm(b)] $\mathfrak G$ is a finite simple group of order
$$|\mathfrak G| = 2^{21}\cdot 3^{16}\cdot 5^2\cdot 7^3\cdot 11
\cdot 13 \cdot 17 \cdot 23 \cdot 29.$$

\item[\rm(c)] $\mathfrak G$ is isomorphic to the commutator
subgroup $P'$ of the finitely presented group $P = \langle
a,b,c,d,e,f,g,h,i,j,k,l \rangle $ with set of defining relations
$\mathcal R(P)$ stated in Lemma \ref{l. nicepres.Fi_24}.

\item[\rm(d)] The character table of $\mathfrak G$ is equivalent to
that of $\Fi_{24}'$ in the Atlas \cite{atlas}, its pp. 200 --202.
\end{enumerate}
\end{theorem}

\begin{proof}
The first three statements hold by Lemma \ref{l.
nicepres.Fi_24}(j) and (k).

(d) Let $G = P'$ be the commutator subgroup of the finitely
presented group $P$. Then $\mathfrak G \cong G$ by (c). The
character table of $G$ has been calculated by means of MAGMA and
the faithful permutation representation $PG$ of $G$ with
stabilizer\\
$M = \langle b_1,c_1,d_1,e_1,f_1,g_1,h_1,i_1,j_1\rangle$ given in
Lemma \ref{l. nicepres.Fi_24}(c).
\end{proof}

\section{Presentation of $2$-central involution centralizer}

In this section we determine generators and a presentation of the
centralizer $\mathfrak H = C_{\mathfrak G}(\mathfrak z)$ of a
$2$-central involution $\mathfrak z$ of the simple subgroup
$\mathfrak G = \langle \mathfrak q, \mathfrak y, \mathfrak w,
\mathfrak t\rangle$ of $\GL_{8671}(13)$. Thus we are able to show
that $\mathfrak G$ satisfies all conditions of Algorithm 2.5 of
\cite{michler1}. In particular, $\mathfrak G = \langle \mathfrak
H, \mathfrak E\rangle$ and $\mathfrak D = N_{\mathfrak
H}(\mathfrak B) = C_{\mathfrak E}(\mathfrak z)$ where $\mathfrak
B$ is the unique maximal elementary abelian normal subgroup of a
well defined Sylow $2$-subgroup $\mathfrak S$ of $\mathfrak G$ and
$\mathfrak E = N_{\mathfrak G}(\mathfrak B)$. However, the free
product $P = H*_DE$ of the groups $H = \mathfrak H$ and $E =
\mathfrak E$ with amalgamated subgroup $D = \mathfrak D$ has
$939,080,024,064$ irreducible representations of dimension $8671$
over $GF(13)$. Therefore we have not tried to find one satisfying
the Sylow $2$-subgroup test of Step 5 c) of Algorithm 7.4.8 of
\cite{michler}.

\begin{proposition}\label{prop. 2-centralFi_24}
Let $\mathfrak G = \langle \mathfrak q, \mathfrak y,\mathfrak t,
\mathfrak w \rangle$ be the subgroup of $Y = \GL_{8671}(13)$
constructed in Proposition \ref{prop. matrixFi_24}. Let $\mathfrak
x = [(\mathfrak y\mathfrak q^2\mathfrak y\mathfrak q\mathfrak
y\mathfrak q^2)^{11}(\mathfrak q^2\mathfrak y^2\mathfrak
q\mathfrak y\mathfrak q\mathfrak y)^{11}(\mathfrak q\mathfrak
y^2\mathfrak q\mathfrak y\mathfrak q\mathfrak y\mathfrak
q\mathfrak y\mathfrak q)^4]^{12}$ and $\mathfrak s = (\mathfrak
y\mathfrak t)^7$. Let $\mathfrak r_1 = (\mathfrak s\mathfrak
y^3)^3$, $\mathfrak r_2 = (\mathfrak y^2\mathfrak w\mathfrak
y\mathfrak s)^6$, $\mathfrak r_3 = (\mathfrak s\mathfrak
y\mathfrak w\mathfrak y\mathfrak s)^2$, $\mathfrak r_4 =
(\mathfrak s\mathfrak w\mathfrak y\mathfrak s\mathfrak w)^6$, and
$\mathfrak v = (\mathfrak w\mathfrak q\mathfrak w\mathfrak
y\mathfrak q\mathfrak y)^7$. Let $\mathfrak H = \langle \mathfrak
r_1, \mathfrak r_2, \mathfrak r_3, \mathfrak r_4, \mathfrak v
\rangle$.

Then the following statements hold:
\begin{enumerate}
\item[\rm(a)] The element $\mathfrak z = (\mathfrak x\mathfrak
y\mathfrak w)^8$ is a $2$-central involution $\mathfrak z =
(\mathfrak x\mathfrak y\mathfrak w)^8$ of $\mathfrak G$ such that
$C_{\mathfrak G}(\mathfrak z) = \mathfrak H$ has order
$2^{21}\cdot 3^7\cdot 5\cdot 7$.

\item[\rm(b)] $\mathfrak H $ has a faithful permutation
representation $P\mathfrak H$ of degree $258048$ with stabilizer
$\mathfrak U_1 = \langle \mathfrak v \mathfrak f_2 \mathfrak v,
(\mathfrak f_1 \mathfrak f_2 \mathfrak f_1 \mathfrak v \mathfrak
f_1)^4, (\mathfrak f_2 \mathfrak f_1 \mathfrak v \mathfrak f_1
\mathfrak v)^4\rangle$, where $\mathfrak f_1 = \mathfrak x
\mathfrak y \mathfrak w$, $\mathfrak f_2 = (\mathfrak x \mathfrak
y \mathfrak x \mathfrak w)^7$, $\mathfrak f_3 = (\mathfrak x
\mathfrak w \mathfrak y \mathfrak w \mathfrak y^2)^7$.

\item[\rm(c)] $\mathfrak S = \langle \mathfrak s, \mathfrak u_i
\mid 1 \le i \le 5 \rangle$ is a Sylow $2$-subgroup of $\mathfrak
E$ contained in $\mathfrak H$ with center $Z(\mathfrak H) =
\langle \mathfrak z \rangle$, where $\mathfrak u_1 = (\mathfrak
x\mathfrak y\mathfrak x)^7$, $\mathfrak u_2 = (\mathfrak
x\mathfrak y\mathfrak x\mathfrak y\mathfrak x\mathfrak y\mathfrak
x)^4$, $\mathfrak u_3 = (\mathfrak x\mathfrak y\mathfrak
x\mathfrak y^2\mathfrak x\mathfrak y^2\mathfrak x\mathfrak
y\mathfrak x\mathfrak y)^2$, $\mathfrak u_4 = (\mathfrak
x\mathfrak y\mathfrak x\mathfrak y^5\mathfrak x\mathfrak y^4)^2$
and $\mathfrak u_5 = (\mathfrak y\mathfrak s \mathfrak w\mathfrak
s\mathfrak w\mathfrak y\mathfrak w)^7$.

\item[\rm(d)] $\mathfrak B = \langle \mathfrak b_i\mid 1 \le i \le
11 \rangle$ is the unique maximal elementary abelian normal
subgroup of $\mathfrak S$ where $\mathfrak b_1 = \mathfrak u_1$,
$\mathfrak b_2 = \mathfrak u_2^2$, $\mathfrak b_3 = (\mathfrak
u_1\mathfrak u_2)^2$, $\mathfrak b_4 = (\mathfrak u_1\mathfrak
u_3)^2$, $\mathfrak b_5 = (\mathfrak u_1\mathfrak u_4)^2$,
$\mathfrak b_6 = (\mathfrak u_2\mathfrak u_4)^4$, $\mathfrak b_7 =
(\mathfrak u_1\mathfrak u_2\mathfrak u_3)^4$, $\mathfrak b_8 =
(\mathfrak u_1\mathfrak u_2\mathfrak u_4)^4$, $\mathfrak b_9 =
(\mathfrak u_1\mathfrak u_3\mathfrak u_4)^4$, $\mathfrak b_{10} =
(\mathfrak u_2^2\mathfrak u_3)^2$, $\mathfrak b_{11} = (\mathfrak
u_1\mathfrak u_2\mathfrak u_4\mathfrak u_2)^2$.

\item[\rm(e)] $N_{\mathfrak G}(\mathfrak B) = \langle \mathfrak x,
\mathfrak y, \mathfrak w, \mathfrak s\rangle = \mathfrak E$.

\item[\rm(f)] $\mathfrak D = \langle \mathfrak r_i \mid 1 \le i
\le 4\rangle = N_{\mathfrak H}(\mathfrak B) = C_{\mathfrak
E}(\mathfrak z)$.

\item[\rm(g)] The Fitting subgroup $\mathfrak O$ of $\mathfrak H$
is extra-special of order $2^{13}$ and center $Z(\mathfrak O) =
\langle \mathfrak z \rangle$. It is generated by the twelve
involutions
\begin{eqnarray*}
&&\mathfrak p_1 = (\mathfrak r_2)^2,\quad \mathfrak p_2 = (\mathfrak r_1\mathfrak r_2)^4,\quad \mathfrak p_3 = (\mathfrak r_3\mathfrak r_4)^3,\quad \mathfrak p_4 = (\mathfrak r_1\mathfrak r_2\mathfrak r_1)^2,\\
&&\mathfrak p_5 = (\mathfrak r_1\mathfrak r_3\mathfrak r_4)^6,\quad \mathfrak p_6 = (\mathfrak r_2^2\mathfrak r_4)^2,\quad \mathfrak p_7 = (r_3\mathfrak r_4\mathfrak r_1)^6,\quad \mathfrak p_8 = (\mathfrak r_4\mathfrak r_1\mathfrak r_3)^6,\\
&&\mathfrak p_9 = (\mathfrak r_1\mathfrak r_2\mathfrak r_3^2)^4,\quad \mathfrak p_{10} = (\mathfrak r_1\mathfrak r_2\mathfrak r_3\mathfrak r_4)^4,\quad \mathfrak p_{11} = (r_1\mathfrak r_2\mathfrak r_1\mathfrak r_2\mathfrak r_4)^5, \quad \mathfrak p_{12} = (\mathfrak r_1\mathfrak r_2\mathfrak r_3^2\mathfrak r_4)^4,\\
&&\mbox{and} \quad \mathfrak z = (\mathfrak p_1\mathfrak p_5)^2.
\end{eqnarray*}

\item[\rm(h)] $\mathfrak O/Z(\mathfrak O)$ has a complement
$\mathfrak K \cong 3U_4(3):2$ in $\mathfrak H/Z(\mathfrak O)$.

\item[\rm(i)] $\mathfrak H = \langle \mathfrak a, \mathfrak b,
\mathfrak p_i \mid 1 \le i \le 12\rangle = \langle \mathfrak
b,\mathfrak c, \mathfrak d,\mathfrak f,\mathfrak z,\mathfrak
p_i\mid 1 \le i \le 12\rangle$, where \\ $\mathfrak a = \mathfrak
r_1\mathfrak r_3^3\mathfrak v\mathfrak r_3$, $\mathfrak b =
[\mathfrak r_1\mathfrak r_3\mathfrak r_1\mathfrak r_3^2\mathfrak
v]^6[\mathfrak r_1\mathfrak r_3\mathfrak r_1\mathfrak r_3\mathfrak
v\mathfrak r_3\mathfrak v]^{12}$, $\mathfrak c = (\mathfrak
a\mathfrak b)^2$, $\mathfrak d = (\mathfrak b\mathfrak c)^7$ and
$\mathfrak f =  [(\mathfrak b\mathfrak a^3)^5\cdot (\mathfrak
a^4\mathfrak b^2\mathfrak a)^9\cdot (\mathfrak a^2\mathfrak
b^3\mathfrak a^4\mathfrak b)^3\cdot(\mathfrak b\mathfrak
a^3)^5]^3$.

\item[\rm(j)] $\mathfrak H$ is isomorphic to the finitely
presented group $H = \langle b,c,d,f,z,p_i\mid 1 \le i \le
12\rangle$ having the following set $\mathcal R(H)$of defining
relations:
\begin{eqnarray*}
&&b^6 = c^9 = d^3 = f^2 = z^2 = 1, \quad p_i^2 = 1 \quad \mbox{for}\quad 1 \le i \le 12,\\
&&(z,b) = (z,c) = (z,d) = (z,f) = 1, \quad (z,p_i) = 1 \quad \mbox{for}\quad 1 \le i \le 12,\\
&&(b,d) = (c,d) = 1,\quad b^{f} = b^5, \quad c^{f} = d(b^3cb^2c^6bcbc), \quad d^{f} = d^2,\\
&&(c^{-1}b^{-1})^7 d = 1,\quad (c^{-1}b)^9 = z,\quad b^{-1}c^{-1}b^{-3}c^{-1}b^3c^{-1}b^3c^{-1}b^{-2}d = 1, \\
&&(bc^{-2}b)^4d^2 = (bc^{-3}bc^2b)^2d = bc^3b^{-2}cb^3c^{-1}bc^{-1}bcb^{-1}cbd^2 = 1,\\
&&b^{-1}c^{-3}b^{-1}c^{-1}bc^{-1}b^2cb^{-1}cbc^3b^{-1}d = b^{-1}c^3bc^{-1}bc^{-1}b^2cb^{-1}cb^{-1}c^{-3}b^{-1}d = 1,\\
&&c^{-2}b^{-2}c^{-1}bc^{-1}b^{-3}c^{-1}bcb^{-2}cbc^{-1}bd^2 = z, \\
&&b^{-2}c^{-1}bc^{-1}b^{-2}c^{-1}b^{-2}cb^{-1}cb^{-1}cb^2c^2d^2 = cb^{-1}c^4bc^{-1}b^{-2}c^{-1}bcb^{-1}c^2bcbd^2 = 1, \\
&&c^{-2}b^{-1}c^2b^{-1}c^{-2}b^{-2}c^2bcb^2c^{-2}b^{-1}d = b^{-1}c^{-3}bc^2bc^{-1}b^{-1}cb^{-1}cb^{-1}c^2b^{-1}c^{-2} = z, \\
&&c^2b^{-1}c^{-1}bc^{-1}b^{-1}c^{-1}b^3c^2b^{-1}c^{-1}bc^{-1}b^2c = z,\\
&&bc^{-2}bc^{-2}bc^{-1}b^3c^{-1}b^2c^{-1}b^{-1}cb^2c^{-1} = b^{-1}c^{-2}b^{-1}c^2bc^{-1}b^3cb^{-1}c^{-2}bc^2b^{-2} = 1, \\
&&b^{-2}c^2b^{-2}c^{-2}b^3c^{-1}b^{-1}c^{-2}bc^{-1}b^{-1}c^{-2}d = cbcb^{-1}cb^{-1}c^{-2}b^3c^2b^{-1}cb^{-1}c^2bcd = 1,\\
&&(p_1,p_2)= (p_1,p_3)= (p_1,p_4)=1, \quad (p_1,p_5)=z, \\
&&(p_1,p_6)= (p_1,p_7)= (p_1,p_8)= (p_1,p_9)= (p_1,p_{10})=1,\\
&&(p_1,p_{11})=(p_1,p_{12})= (p_2,p_3)=z,\\
&&(p_2,p_4)= (p_2,p_5)=1, \quad (p_2,p_6)=z, \quad (p_2,p_7)=1, \quad (p_2,p_8)=z, \\
&&(p_2,p_9)=1,\quad (p_2,p_{10})=z, \quad (p_2,p_{11})=1, \quad (p_2,p_{12})=z, \\
&&(p_3,p_4)=1, \quad (p_3,p_5)=z, \quad (p_3,p_6)=1, \quad (p_3,p_7)= (p_3,p_8) = (p_3,p_9)=z, \\
&&(p_3,p_{10})=1, \quad (p_3,p_{11}) = z, \quad (p_3,p_{12}) = 1,\\
&&(p_4,p_5)= (p_4,p_6)=1, \quad (p_4,p_7) = (p_4,p_8)=z, \\
&&(p_4,p_9)= (p_4,p_{10})=1, \quad (p_4,p_{11})= (p_4,p_{12})=z,\\
&&(p_5,p_6)=z, \quad (p_5,p_7) = (p_5,p_8)=1, \\
&&(p_5,p_9)= (p_5,p_{10})=1, \quad (p_5,p_{11})=z,\quad (p_5,p_{12})= 1,\\
&&(p_6,p_7)=1, \quad (p_6,p_8) = (p_6,p_9)=z, \quad (p_6,p_{10})= (p_6,p_{11})= (p_6,p_{12})=1, \\
&&(p_7,p_8)= (p_7,p_9) = (p_7,p_{10})=z, \quad (p_7,p_{11})= (p_7,p_{12})=1, \\
&&(p_8,p_9)=z, \quad (p_8,p_{10}) = 1, \quad (p_8,p_{11}) = (p_8,p_{12}) = (p_9,p_{10}) = (p_9,p_{11})=z, \\
&&(p_9,p_{12}) = 1, \quad (p_{10},p_{11})= (p_{10},p_{12})= (p_{11},p_{12})=z,\\
&&p_1^b  p_1 p_5 p_6 p_7 p_8 p_{12}= p_2^b  p_5 p_8 p_9 p_{12} = p_3^b  p_2 p_3 p_4 p_5 p_6 p_7 p_8 p_9 p_{10} p_{11}=1,\\
&&p_4^b p_1 p_2 p_4 p_5 p_8 p_9 p_{10} p_{11} p_{12} = p_5^b  p_7 p_8 p_9 p_{12} = z,\quad p_6^b  p_1 p_3 p_6 p_7 p_8 p_9= 1,\\
&&p_7^b  p_2 p_5 p_6 p_9 p_{10} p_{11} = p_8^b  p_1 p_2 p_5 p_7 p_{10} p_{11} p_{12} = 1,\\
&&p_9^b  p_3 p_5 p_7 p_{10} p_{11} p_{12}= z,\quad p_{10}^b  p_1 p_4 p_5 p_7 p_8 p_9 p_{10} = 1,\\
&&p_{11}^b  p_2 p_4 p_6 p_7 p_8 p_9 p_{11} = z,\quad p_{12}^b  p_1 p_2 p_4 p_5 p_9 p_{10} p_{11}=1,\\
&&p_1^c  p_1 p_2 p_4 p_6 p_7 p_8 = 1,\quad p_2^c  p_1 p_3 p_4 p_7 p_8 p_9 = z,\quad p_3^c  p_1 p_3 p_6 p_7 p_8 p_9 p_{12}= 1,\\
&&p_4^c  p_3 p_4 p_5 p_6 p_8 p_{10} p_{11} = z,\quad p_5^c  p_1 p_2 p_5 p_7 p_9 p_{12} = z,\quad p_6^c  p_3 p_6 p_{11} p_{12}= 1,\\
\end{eqnarray*}
\begin{eqnarray*}
&&p_7^c  p_1 p_2 p_3 p_6 p_7 p_8 p_9 p_{10} = z,\quad p_8^c  p_2 p_3 p_4 p_5 p_6 p_9 p_{12} = 1,\quad p_9^c  p_1 p_2 p_3 p_4 p_7 p_{11}= z,\\
&&p_{10}^c  p_2 p_3 p_5 p_7 p_9 p_{10} p_{12} = z,\quad p_{11}^c  p_1 p_4 p_7 p_8 p_{11} p_{12} = z,\quad p_{12}^c  p_1 p_2 p_8 p_9= 1,\\
&&p_1^d  p_2 p_6 p_9 p_{10} = p_2^d p_1 p_4 p_5 p_8 p_9 p_{12} = p_3^d  p_1 p_2 p_3 p_4 p_6 p_9 = z,\\
&&p_4^d  p_1 p_2 p_3 p_4 p_9 p_{10} = p_5^d  p_1 p_3 p_5 p_7 p_{10} p_{12} = z,\quad p_6^d  p_2 p_3 p_6 p_9 = z,\\
&&p_7^d  p_3 p_4 p_5 p_7 p_8 p_{10} p_{11} = z,\quad p_8^d  p_2 p_3 p_4 p_5 p_6 p_8 p_{11} p_{12} = z,\quad p_9^d  p_2 p_3 p_5 p_8 p_{12} = z, \\
&&p_{10}^d  p_2 p_4 p_9 p_{10} = p_{11}^d  p_1 p_8 p_9 p_{10} p_{11} = 1,\quad p_{12}^d  p_1 p_2 p_4 p_7 p_8 p_9 p_{10} p_{11} = z, \\
&&p_1^f  p_3 p_9 p_{11} p_{12} = z,\quad p_2^f  p_1 p_5 p_6 p_9 p_{12} = z,\quad p_3^f  p_1 p_2 p_3 p_7 p_8 p_9 =1, \\
&&p_4^f  p_1 p_2 p_4 p_6 p_7 p_8 p_9 p_{10} = p_5^f  p_1 p_5 p_6 p_8 p_{10} p_{12} = 1,\quad p_6^f  p_4 p_5 p_7 p_9 = z, \\
&&p_7^f  p_3 p_4 p_5 p_8 p_{10} p_{11} = p_8^f  p_3 p_4 p_8 p_{11} p_{12} = 1,\quad p_9^f  p_2 p_3 p_6 p_7 p_8 p_9 p_{10} p_{11} p_{12}=z, \\
&&p_{10}^f  p_4 p_5 p_6 p_7 p_9 p_{10} = 1,\quad p_{11}^f  p_1 p_4 p_6 p_9 p_{10} p_{11} = z,\quad p_{12}^f  p_1 p_4 p_9 p_{12}=1. \\
\end{eqnarray*}

\item[\rm(k)] $H = \langle b,c,d,f,z,p_i\mid 1 \le i \le
12\rangle$ has a faithful permutation representation of degree
$258048$ with stabilizer $L_1 = \langle b, c^3, f,
p_5p_6p_7p_9p_{10}p_{11}\rangle$.

\item[\rm(l)] In $H = \langle b,c,d,f,z,p_i\mid 1 \le i \le
12\rangle$ let $r = [(fp_1)^2(cfb)^3(fp_1)^2(fc^4)]^6$ and $p =
[(fp_1)^2(cfb)^3(fp_1)^2(b^2fc)]^5$. Then $H = \langle b, p,
r\rangle$, and it has $167$ conjugacy classes whose
representatives are given in Table \ref{Fi_24 cc H}.

\item[\rm(m)] Table \ref{Fi_24 ct H} is the character table of
$H$.
\end{enumerate}
\end{proposition}

\begin{proof}
(a) By Lemma \ref{l. nicepres.Fi_24} there is an isomorphism
$\sigma: \mathfrak G \rightarrow G$ such that $\sigma(\mathfrak q)
= q_1$, $\sigma(\mathfrak y) = y_1$, $\sigma(\mathfrak w) = w_1$,
and $\sigma(\mathfrak t) = t_1$. In particular, $G = \langle
q_1,y_1,w_1,t_1\rangle$. Let $x_1$ be as in Lemma \ref{l.
nicepres.Fi_24}(i). Let $z_1 = \sigma (\mathfrak z) =
(x_1y_1w_1)^8)$. Using the faithful permutation representation
$PG$ of $G$ given in Lemma \ref{l. nicepres.Fi_24}(a) and MAGMA
one sees that $C_G(z_1)$ has order $2^{21}\cdot 3^7\cdot 5\cdot
7$. Hence $\mathfrak z$ is a $2$-central involution of $\mathfrak
G$ by Lemma \ref{l. nicepres.Fi_24}(k).

In $G$ let $r_1 = s_1y_1^3$, $r_2 = y_1^2w_1y_1s_1$, $r_3 =
(s_1y_1w_1s_1)^2$, $r_4 = (s_1w_1y_1s_1w_1)^6$, and $v = (w_1
q_1w_1y_1q_1y_1)^7$. Let $H = \langle r_1, r_2, r_3, r_4, v\rangle
= \sigma(\mathfrak H)$. Then Proposition \ref{prop. E(Fi_24)}(i),
(j) and Lemma \ref{l. nicepres.Fi_24}(l) imply that $H \le
C_G(z_1)$. Another application of $PG$ and MAGMA yields that
$C_G(z_1) = H$ and that $Z(H) = \langle z_1 \rangle$.

(b) From Proposition \ref{prop. E(Fi_24)}(j) we deduce that $U =
C_{G_1}(z_1) = \langle f_1, f_2,f_3,v\rangle$ has order
$2^{18}\cdot 3^5\cdot 5$, where $f_1 = x_1y_1w_1$, $f_2 =
(x_1y_1x_1w_1)^7$, $f_3 = (x_1w_1y_1w_1y_1^2)^7$.  By Proposition
\ref{prop. E(Fi_24)}(k) $U$ has a subgroup $U_1$ of order
$2^9\cdot 3^5\cdot 5$ which does not contain $z$ generated by the
elements of the statement. Hence $H_1 = C_G(z_1)$ has a faithful
permutation representation $PH_1$ of degree $258048$ by (a).

(c) By Proposition \ref{prop. E(Fi_24)}(l) the subgroup $\mathfrak
S$ is a Sylow $2$-subgroup of $\mathfrak E = \langle \mathfrak x,
\mathfrak y, \mathfrak w, \mathfrak s\rangle$. Lemmas \ref{l.
nicepres.Fi_24}(h) states that $\mathfrak E \cong E$ where $E$ is
the finitely presented group of Lemma \ref{l. M24-extensions}.
Thus $|\mathfrak S| = 2^{21}$ by Table \ref{Fi_24 ct_E}. Hence
$\mathfrak S$ is a Sylow $2$-subgroup of $\mathfrak G$ by (a). The
equality $Z(\mathfrak H) = \langle \mathfrak z \rangle$ has been
checked computationally in $PG$.

(d) This statement follows now immediately from Proposition
\ref{prop. E(Fi_24)}(l).

(e) This assertion is true by (d) and Lemma \ref{l.
nicepres.Fi_24}(h), (i) and (k).

(f) By Proposition \ref{prop. E(Fi_24)}(i) and Lemma \ref{l.
nicepres.Fi_24}(l) we know that $\mathfrak D = C_{\mathfrak
E}(\mathfrak z) = \langle \mathfrak r_i \mid 1 \le i \le
4\rangle$, where $\mathfrak z = \kappa(z_1)$. Thus it suffices to
check that $N_H(B_1) = \langle r_1,r_2,r_3,r_4\rangle$. This has
been done using the faithful permutation representation $PG$ and
MAGMA.

(g) We verified computationally that the Fitting subgroup $O$ of
$H$ is extra-special of order $2^{13}$ and has center $Z(O) =
\langle z \rangle$. The twelve involutions $p_i$ generating the
subgroup $O$ have been calculated with MAGMA by means of $PG$ and
the program $\verb"GetShortGens(H,O)"$. We also verified that $z =
(p_1p_5)^2$.

(h) Let $\alpha: H \rightarrow H/Z(H) = H_1$ be the canonical
epimorphism of $H = C_G(z)$ with kernel $Z(H) = \langle z
\rangle$. By (c) $H$ has a faithful permutation representation
$PH$ with stabilizer $U_1 = \langle vf_2v, (f_1f_2f_1vf_1)^4,
(f_2f_1vf_1v)^4\rangle$. As $z$ does not belong to $U_1$ its
cosets provide a faithful permutation representation $PH$ of $H$
having degree $258048$ by (b). Using MAGMA we checked that the
subgroup $\alpha(U_1)$ of $H_1$ is a stabilizer of a faithful
permutation representation $PH_1$ of $H_1$ of degree $129024$.
Applying $PH_1$ and the MAGMA command
$\verb"DegreeReduction(H_1)"$ MAGMA calculated a faithful
permutation representation $PH_{11}$ of $H_1 = \langle
\alpha(r_i),\alpha(v)\rangle$ of degree $504$ with stabilizer\\
$U_{11} = \langle \alpha(r_1)\alpha(r_3)\alpha(r_1),
[\alpha(r_1)\alpha(v)\alpha(r_3)\alpha(v)]^3,
[\alpha(r_1)\alpha(v)\alpha(r_1)\alpha(r_3)\alpha(v)\alpha(r_3)]^3\rangle$.

Clearly, $V = O/Z(H)$ is an elementary abelian normal subgroup of
$H_1$ of order $2^{12}$. Using $PH_{11}$ and the command
$\verb"HasComplement(H_1,V)"$ MAGMA established a complement $K_1$
of $V$ in $H_1$. By means of the command
$\verb"CompositionFactors(K_1)"$ we saw that $|K_1 : K'_1| = 2$,
$|Z(K'_1)| = 3$, and $K'/Z(K') \cong U_4(3)$.

(i) Using $PH_{11}$ a MAGMA calculation confirmed that $H_1$ is generated by\\ $a_1 = \alpha(r_1)(\alpha(r_3)^3\alpha(v)\alpha(r_3))$ and\\
$b_1 = [\alpha(r_1)\alpha(r_3)\alpha(r_1)\alpha(r_3)^2\alpha(v)]^6[\alpha(r_1)\alpha(r_3)\alpha(r_1)\alpha(r_3)\alpha(v)\alpha(r_3)\alpha(v)]^{12}$.\\
Both generators have order $6$. Furthermore, $K_1' = \langle
b_1,c_1\rangle$, where $c_1 = (a_1b_1)^2$. Using the command
$\verb"GetShortGens(K_1',Z(K_1'))"$ we observed that the center
$Z(K_1')$ of $K_1'$ is generated by the element $d_1 = (b_1c_1)^7$
of order $3$.

Let $\beta: K_1' \rightarrow K_1'/Z(K_1') = K_2$ be the canonical
epimorphism of $K_1$ with kernel $Z(K_1') = \langle d_1 \rangle$.
Let $U_{12} = U_{11} \cap K_1'$ and $U_{13} = \langle U_{12},
d_1\rangle$. Then $K_2$ has a faithful permutation representation
$PK_2$ of degree $126$ with stabilizer $\beta(U_{13})$. Let $a_1 =
\beta(b_1)$ and $a_2 = \beta(c_1)$. Then $K_2 = \langle a_1,
a_2\rangle$. Using the command $\verb"FPGroup(K_2)"$ MAGMA
calculates the following set $\mathcal R(K_2)$ of defining
relations of $K_2$:
\begin{eqnarray*}
&&a_1^6 = 1, \quad a_2^9 = 1,\\
&&(a_2^{-1}a_1^{-1})^7 = 1, \quad (a_1a_2^{-2}a_1)^4 = 1,\quad(a_2^{-1} a_1)^9 = 1,\\
&&a_1^{-1} a_2^{-1} a_1^{-3} a_2^{-1} a_1^3 a_2^{-1} a_1^3a_2^{-1} a_1^{-2} = a_1 a_2^3 a_1^{-2} a_2 a_1^3 a_2^{-1} a_1 a_2^{-1} a_1 a_2a_1^{-1} a_2 a_1 = 1,\\
&&a_1^{-1}a_2^{-3}a_1^{-1}a_2^{-1}a_1a_2^{-1}a_1^2a_2a_1^{-1}a_2a_1a_2^3a_1^{-1} = 1,\\
&&a_1^{-1} a_2^3 a_1 a_2^{-1} a_1 a_2^{-1} a_1^2 a_2 a_1^{-1}a_2 a_1^{-1} a_2^{-3} a_1^{-1} = (a_1a_2^{-3}a_1a_2^2a_1)^2 = 1,\\
&&a_2^{-2}a_1^{-2}a_2^{-1}a_1a_2^{-1}a_1^{-3}a_2^{-1}a_1a_2a_1^{-2}a_2a_1a_2^{-1}a_1 = 1,\\
&&a_1^{-2}a_2^{-1}a_1a_2^{-1}a_1^{-2}a_2^{-1}a_1^{-2}a_2a_1^{-1}a_2a_1^{-1}a_2a_1^2a_2^2 = 1,\\
&&a_2 a_1^{-1} a_2^4 a_1 a_2^{-1} a_1^{-2} a_2^{-1} a_1 a_2a_1^{-1} a_2^2 a_1 a_2 a_1 = 1,\\
&&a_2^{-2} a_1^{-1} a_2^2 a_1^{-1} a_2^{-2} a_1^{-2} a_2^2 a_1a_2 a_1^2 a_2^{-2} a_1^{-1} = 1,\\
&&a_1^{-1}a_2^{-3}a_1a_2^2a_1a_2^{-1}a_1^{-1}a_2a_1^{-1}a_2a_1^{-1}a_2^2a_1^{-1}a_2^{-2} = 1,\\
&&a_2^2a_1^{-1}a_2^{-1}a_1a_2^{-1}a_1^{-1}a_2^{-1}a_1^3a_2^2a_1^{-1}a_2^{-1}a_1a_2^{-1}a_1^2a_2 = 1,\\
&&a_1a_2^{-2}a_1a_2^{-2}a_1a_2^{-1}a_1^3a_2^{-1}a_1^2a_2^{-1}a_1^{-1}a_2a_1^2a_2^{-1} = 1,\\
&&a_1^{-1}a_2^{-2}a_1^{-1}a_2^2a_1a_2^{-1}a_1^3a_2a_1^{-1}a_2^{-2}a_1a_2^2a_1^{-2} = 1,\\
&&a_1^{-2}a_2^2a_1^{-2}a_2^{-2}a_1^3a_2^{-1}a_1^{-1}a_2^{-2}a_1a_2^{-1}a_1^{-1}a_2^{-2} = 1,\\
&&a_2 a_1 a_2 a_1^{-1} a_2 a_1^{-1} a_2^{-2} a_1^3 a_2^2a_1^{-1} a_2 a_1^{-1} a_2^2 a_1 a_2 = 1.\\
\end{eqnarray*}

Since $K_1' = \langle b_1,c\rangle$, $Z(K_1') = \langle
d_1\rangle$, $d_1 = (b_1c_1)^7 \neq 1$, $d_1^3 = 1$ we lift the
relations of $\mathcal R(K_2)$ to $K_1'$ by evaluating them in the
permutation representation of $K_1' = \langle b_1,c_1,d_1\rangle$.
Thus we obtain the following set $\mathcal R(K_1')$ of defining
relations of $K_1'$:
\begin{eqnarray*}
&&b_1^6 = c_1^9 = d_1^3 = 1,\quad (b_1,d_1) = (c_1,d_1) = 1,\\
&&(c_1^{-1}b_1^{-1})^7 d_1 = 1,\quad (c_1^{-1}b_1)^9 = 1,\quad b_1^{-1}c_1^{-1}b_1^{-3}c_1^{-1}b_1^3c_1^{-1}b_1^3c_1^{-1}b_1^{-2}d_1 = 1, \\
&&(b_1c_1^{-2}b_1)^4d_1^2 = (b_1c_1^{-3}b_1c_1^2b_1)^2d_1 = b_1c_1^3b_1^{-2}c_1b_1^3c_1^{-1}b_1c_1^{-1}b_1c_1b_1^{-1}c_1b_1d_1^2 = 1,\\
&&b_1^{-1}c_1^{-3}b_1^{-1}c_1^{-1}b_1c_1^{-1}b_1^2c_1b_1^{-1}c_1b_1c_1^3b_1^{-1}d_1 = 1,\\
&&b_1^{-1}c_1^3b_1c_1^{-1}b_1c_1^{-1}b_1^2c_1b_1^{-1}c_1b_1^{-1}c_1^{-3}b_1^{-1}d_1 = 1,\\
&&c_1^{-2}b_1^{-2}c_1^{-1}b_1c_1^{-1}b_1^{-3}c_1^{-1}b_1c_1b_1^{-2}c_1b_1c_1^{-1}b_1d_1^2 = 1, \\
&&b_1^{-2}c_1^{-1}b_1c_1^{-1}b_1^{-2}c_1^{-1}b_1^{-2}c_1b_1^{-1}c_1b_1^{-1}c_1b_1^2c_1^2d_1^2 = 1,\\
&&c_1b_1^{-1}c_1^4b_1c_1^{-1}b_1^{-2}c_1^{-1}b_1c_1b_1^{-1}c_1^2b_1c_1b_1d_1^2 = 1,\\
&&c_1^{-2}b_1^{-1}c_1^2b_1^{-1}c_1^{-2}b_1^{-2}c_1^2b_1c_1b_1^2c_1^{-2}b_1^{-1}d_1 = 1,\\
&&b_1^{-1}c_1^{-3}b_1c_1^2b_1c_1^{-1}b_1^{-1}c_1b_1^{-1}c_1b_1^{-1}c_1^2b_1^{-1}c_1^{-2} = 1, \\
&&c_1^2b_1^{-1}c_1^{-1}b_1c_1^{-1}b_1^{-1}c_1^{-1}b_1^3c_1^2b_1^{-1}c_1^{-1}b_1c_1^{-1}b_1^2c_1 = 1,\\
&&b_1c_1^{-2}b_1c_1^{-2}b_1c_1^{-1}b_1^3c_1^{-1}b_1^2c_1^{-1}b_1^{-1}c_1b_1^2c_1^{-1} = 1,\\
&&b_1^{-1}c_1^{-2}b_1^{-1}c_1^2b_1c_1^{-1}b_1^3c_1b_1^{-1}c_1^{-2}b_1c_1^2b_1^{-2} = 1, \\
&&b_1^{-2}c_1^2b_1^{-2}c_1^{-2}b_1^3c_1^{-1}b_1^{-1}c_1^{-2}b_1c_1^{-1}b_1^{-1}c_1^{-2}d_1 = 1,\\
&&c_1b_1c_1b_1^{-1}c_1b_1^{-1}c_1^{-2}b_1^3c_1^2b_1^{-1}c_1b_1^{-1}c_1^2b_1c_1d_1 = 1.\\
\end{eqnarray*}

Using the faithful permutation representation $PH_{11}$ of $H_1$
and the MAGMA command $\verb"HasComplement(K_1, K_1')")$ we see
that $K_1'$ has a complement $\langle f_1 \rangle$ of order $2$
such that $b_1^{f_1} \in \{b_1, b_1^{-1}\}$. In order to find a
generator $f_1$ of a suitable complement we searched for an
involution $f_1$ of $K_1$ so that at least one of the elements
$c_1^{f_1}, d_1\cdot c_1^{f_1}$, or $d_1^2 \cdot c_1^{f_1}$ is in
the collection of all words in $b_1$ and $c_1$ of length at most
$16$. This search was successful. The involution\\
$f_1 = [(b_1a_1^3)^5\cdot (a_1^4b_1^2a_1)^9\cdot
(a_1^2b_1^3a_1^4b_1)^3\cdot(b_1a_1^3)^5]^3$ of $K_1 = \langle a_1,
b_1\rangle$ satisfies the following set $\mathcal R(f_1)$ of
relations: $f_1^2 = 1$, $b_1^{f_1} = b_1^5$, $c_1^{f_1} =
d_1(b_1^3c_1b_1^2c_1^6b_1c_1b_1c_1)$, $d_1^{f_1} = d_1^2$. Hence
$K_1 = \langle a_1,b_1 \rangle = \langle b_1,c_1,d_1,f_1\rangle$
has a set $\mathcal R(K_1)$ of defining relations consisting of
$\mathcal R(K_1')$ and $\mathcal R(f_1)$.

The elementary abelian normal subgroup $V = \alpha(O)$ is
generated by the involutions $q_i = \alpha(p_i)$ which are the
images of the twelve generating involutions $r_i$ of the Fitting
subgroup $O$ of $H$. Using the faithful permutation representation
$PH_{11}$ we calculated the images $q_i^{x}$ in $V$ for all 4
generators $x \in \{b_1,c_1,d_1,f_1\}$ of $K_1$. Thus we obtained
the following set of essential relations ${\mathcal R_2}(V\rtimes
K_1)$ of the semi-direct product $(V\rtimes K_1)$:
\begin{eqnarray*}
&&q_1^{b_1}  q_1 q_5 q_6 q_7 q_8 q_{12}= q_2^{b_1}  q_5 q_8 q_9 q_{12} = q_3^{b_1}  q_2 q_3 q_4 q_5 q_6 q_7 q_8 q_9 q_{10} q_{11}=1,\\
&&q_4^{b_1} q_1 q_2 q_4 q_5 q_8 q_9 q_{10} q_{11} q_{12} = q_5^{b_1}  q_7 q_8 q_9 q_{12} = q_6^{b_1}  q_1 q_3 q_6 q_7 q_8 q_9= 1,\\
&&q_7^{b_1}  q_2 q_5 q_6 q_9 q_{10} q_{11} = q_8^{b_1}  q_1 q_2 q_5 q_7 q_{10} q_{11} q_{12} = q_9^{b_1}  q_3 q_5 q_7 q_{10} q_{11} q_{12}= 1,\\
&&q_{10}^{b_1}  q_1 q_4 q_5 q_7 q_8 q_9 q_{10} = q_{11}^{b_1}  q_2 q_4 q_6 q_7 q_8 q_9 q_{11} = q_{12}^{b_1}  q_1 q_2 q_4 q_5 q_9 q_{10} q_{11}=1,\\
&&q_1^{c_1}  q_1 q_2 q_4 q_6 q_7 q_8 = q_2^{c_1}  q_1 q_3 q_4 q_7 q_8 q_9 = q_3^{c_1}  q_1 q_3 q_6 q_7 q_8 q_9 q_{12}= 1,\\
&&q_4^{c_1}  q_3 q_4 q_5 q_6 q_8 q_{10} q_{11} = q_5^{c_1}  q_1 q_2 q_5 q_7 q_9 q_{12} = q_6^{c_1}  q_3 q_6 q_{11} q_{12}= 1,\\
&&q_7^{c_1}  q_1 q_2 q_3 q_6 q_7 q_8 q_9 q_{10} = q_8^{c_1}  q_2 q_3 q_4 q_5 q_6 q_9 q_{12} = q_9^{c_1}  q_1 q_2 q_3 q_4 q_7 q_{11}= 1,\\
&&q_{10}^{c_1}  q_2 q_3 q_5 q_7 q_9 q_{10} q_{12} = q_{11}^{c_1}  q_1 q_4 q_7 q_8 q_{11} q_{12} = q_{12}^{c_1}  q_1 q_2 q_8 q_9= 1,\\
&&q_1^{d_1}  q_2 q_6 q_9 q_{10} = q_2^{d_1} q_1 q_4 q_5 q_8 q_9 q_{12} = q_3^{d_1}  q_1 q_2 q_3 q_4 q_6 q_9 = 1,\\
\end{eqnarray*}
\begin{eqnarray*}
&&q_4^{d_1}  q_1 q_2 q_3 q_4 q_9 q_{10} = q_5^{d_1}  q_1 q_3 q_5 q_7 q_{10} q_{12} = q_6^{d_1}  q_2 q_3 q_6 q_9 = 1,\\
&&q_7^{d_1}  q_3 q_4 q_5 q_7 q_8 q_{10} q_{11} = q_8^{d_1}  q_2 q_3 q_4 q_5 q_6 q_8 q_{11} q_{12} = q_9^{d_1}  q_2 q_3 q_5 q_8 q_{12} =1, \\
&&q_{10}^{d_1}  q_2 q_4 q_9 q_{10} = q_{11}^{d_1}  q_1 q_8 q_9 q_{10} q_{11} = q_{12}^{d_1}  q_1 q_2 q_4 q_7 q_8 q_9 q_{10} q_{11} =1, \\
&&q_1^{f_1}  q_3 q_9 q_{11} q_{12} = q_2^{f_1}  q_1 q_5 q_6 q_9 q_{12} = q_3^{f_1}  q_1 q_2 q_3 q_7 q_8 q_9 =1, \\
&&q_4^{f_1}  q_1 q_2 q_4 q_6 q_7 q_8 q_9 q_{10} = q_5^{f_1}  q_1 q_5 q_6 q_8 q_{10} q_{12} = q_6^{f_1}  q_4 q_5 q_7 q_9=1, \\
&&q_7^{f_1}  q_3 q_4 q_5 q_8 q_{10} q_{11} = q_8^{f_1}  q_3 q_4 q_8 q_{11} q_{12} = q_9^{f_1}  q_2 q_3 q_6 q_7 q_8 q_9 q_{10} q_{11} q_{12}=1, \\
&&q_{10}^{f_1}  q_4 q_5 q_6 q_7 q_9 q_{10} = q_{11}^{f_1}  q_1 q_4 q_6 q_9 q_{10} q_{11} = q_{12}^{f_1}  q_1 q_4 q_9 q_{12}=1. \\
\end{eqnarray*}

Hence the set $\mathcal R(H_1)$ of defining relations of the
semi-direct product $H_1 = (V\rtimes K_1)$ consists $\mathcal
R(K_1)$, ${\mathcal R_2}(V\rtimes K_1)$ and the following
relations:
\begin{eqnarray*}
&&q_j^2 = 1 \quad \mbox{for all}\quad 1\leq j\leq 12,\\
&&q_k\cdot q_j = q_j\cdot q_k \quad \mbox{for all}\quad 1\leq j,k\leq 12.\\
\end{eqnarray*}

In order to get a presentation of $H$ we lift the generators $a_1$
and $b_1$ of $K_1$ to $H$. Clearly, $a = r_1 r_3^3 v r_3$ and $b =
(r_1 r_3 r_1 r_3^2 v)^6 (r_1 r_3 r_1 r_3 v r_3 v)^{12}$ of $H$ map
onto $a_1$ and $b_1$ of $H_1$, respectively. Let $c = (ab)^2$, $d
= (bc)^7$ and $f = [(ba^3)^5\cdot (a^4b^2a)^9\cdot
(a^2b^3a^4b)^3\cdot(ba^3)^5]^3$. Then $c_1$, $d_1$, $f_1$ in $H_1$
are images of $c$, $d$, $f$ in $H$ under $\alpha$, respectively.
Since $z = (p_1p_5)^2$ generates the center $Z(O)$ of the Frattini
subgroup $O = \langle p_i \mid 1 \le i \le 12 \rangle$ it follows
that
$$H = \langle a, b, p_i\mid 1\le i \le 12\rangle = \langle b,c,d,f, p_i\mid 1\le i \le 12\rangle.$$
The set of defining relations $\mathcal R(H)$ of $H = \langle
b,c,d,f, p_i\mid 1\le i \le 12\rangle$ has been obtained by
evaluating the lifted equations of the presentation $\mathcal
R(H_1)$ of $H_1$ in the permutation representation $PH$ with
stabilizer $U_1$ defined in the proof of (c). The resulting
equations of $\mathcal R(H)$ are stated in assertion (k).

The map $\mathfrak H \rightarrow H$ sending each generator
$\mathfrak x$ of $\mathfrak H$ in (i) to the corresponding
generator $x \in H$ in (j) is an isomorphism by (a) and the order
of $H$.

(k) The given stabilizer of the group $H = \langle b,c,d,f,
p_i\mid 1\le i \le 12\rangle$ has been found as follows. In the
original permutation representation of the finitely presented
group $G$ of degree $306936$ we checked that the subgroup $L =
\langle b, c^3, f\rangle$ has index $1032192$ in $H$ and that $z
\notin L$. Using then MAGMA and the command
$\verb"MyCosetAction(H,L: maxsize:=10000000)"$ we verified that
$L$ has the same index in the finitely presented group $H =
\langle b,c,d,f, p_i\mid 1\le i \le 12\rangle$. In the
corresponding permutation representation $pH$ of this group we
searched then for an element $p \in O = \langle p_i \rangle$ such
that $z \notin L_1 = \langle L, p\rangle$. MAGMA found an
involution $p \in O$ with these properties. Since $O$ is
extra-special of order $2^{13}$ the command
$\verb"LookupWord(O,p)"$ worked well. The word of $p$ is stated in
the assertion. Using the command
$\verb"MyCosetAction(H,L_1:maxsize:=10000000)"$ MAGMA established
in $70$ seconds the index $|H : L_1| = 258048$.

(l) Both elements $r$ and $p$ of $H$ have order $6$. Using the
faithful permutation representation $PH$ of $H$ with stabilizer
$L_1$ and MAGMA it has been verified that $H = \langle
r,p,b\rangle$. Since $H_1 = H/\langle z\rangle$ has a faithful
permutation representation of degree $504$ we used it and
Kratzer's Algorithm 5.3.18 of \cite{michler} to calculate a system
of representatives of the classes of $H_1 = \langle
\alpha(r),\alpha(p),\alpha(b)\rangle$. $H_1$ has $123$ conjugacy
classes. Their representatives have been lifted to $H$. Using $PH$
we have checked the conjugacy of the lifted representatives and
the products with the central involution $z$ of $H$.

(m) The character table of $H$ was calculated automatically by
MAGMA using $PH$.
\end{proof}

\section{Group order}

In this section we check the group order of $\mathfrak G$ by means
of Thompson's group order formula and Theorem 6.1.4 of
\cite{michler2}.

\begin{proposition}\label{prop. Fi24'order}
Let $\mathfrak G = \langle \mathfrak q, \mathfrak y,\mathfrak t,
\mathfrak w \rangle$ be the subgroup of $\GL_{8671}(13)$
constructed in Proposition \ref{prop. matrixFi_24}. Let $\mathfrak
x = [(\mathfrak y\mathfrak q^2\mathfrak y\mathfrak q\mathfrak
y\mathfrak q^2)^{11}(\mathfrak q^2\mathfrak y^2\mathfrak
q\mathfrak y\mathfrak q\mathfrak y)^{11}(\mathfrak q\mathfrak
y^2\mathfrak q\mathfrak y\mathfrak q\mathfrak y\mathfrak
q\mathfrak y\mathfrak q)^4]^{12}$ and $\mathfrak s = (\mathfrak
y^5\mathfrak t)^7$.

Let $\mathfrak r_1 = (\mathfrak s\mathfrak y^3)^3$, $\mathfrak r_2
= (\mathfrak y^2\mathfrak w\mathfrak y\mathfrak s)^6$, $\mathfrak
r_3 = (\mathfrak s\mathfrak y\mathfrak w\mathfrak y\mathfrak
s)^2$, $\mathfrak r_4 = (\mathfrak s\mathfrak w\mathfrak
y\mathfrak s\mathfrak w)^6$, and $\mathfrak v = (\mathfrak
w\mathfrak q\mathfrak w\mathfrak y\mathfrak q\mathfrak y)^7$. Let
$\mathfrak a = \mathfrak r_1\mathfrak r_3^3\mathfrak v\mathfrak
r_3$, $\mathfrak b = [\mathfrak r_1\mathfrak r_3\mathfrak
r_1\mathfrak r_3^2\mathfrak v]^6[\mathfrak r_1\mathfrak
r_3\mathfrak r_1\mathfrak r_3\mathfrak v\mathfrak r_3\mathfrak
v]^{12}$,\\ $\mathfrak f = [(\mathfrak b\mathfrak a^3)^5
(\mathfrak a^4\mathfrak b^2\mathfrak a)^9(\mathfrak a^2\mathfrak
b^3\mathfrak a^4\mathfrak b)^3(\mathfrak b\mathfrak a^3)^5]^3$,
$\mathfrak r = [(\mathfrak f \mathfrak r_2^2)^2[(\mathfrak
a\mathfrak b)^2\mathfrak f\mathfrak b]^3(f\mathfrak
r_2^2)^2(\mathfrak f(\mathfrak a\mathfrak b)^8)]^6$,\\ and
$\mathfrak p = [(\mathfrak f\mathfrak r_2^2)^2[(\mathfrak
a\mathfrak b)^2\mathfrak f\mathfrak b]^3(\mathfrak f\mathfrak
r_2^2)^2(\mathfrak b^2\mathfrak f(\mathfrak a\mathfrak b)^2)]^5$.

Then the following assertions hold:
\begin{enumerate}
\item[\rm(a)] $\mathfrak z = (\mathfrak x\mathfrak y\mathfrak
w)^8$ is a $2$-central involution of $\mathfrak G$ with
centralizer $\mathfrak H = C_{\mathfrak G}(\mathfrak z) = \langle
\mathfrak b, \mathfrak p, \mathfrak r\rangle$. Furthermore,
$\mathfrak H$ has $9$, conjugacy classes of involutions $2_i$, $1
\le i \le 9$ with representatives given in Table \ref{Fi_24 cc H},
and $|\mathfrak H| = 2^{21}\cdot 3^7\cdot 5\cdot 7$.

\item[\rm(b)] $\mathfrak u = (\mathfrak y \mathfrak s)^{14}$ is an
involution of $\mathfrak G$ with centralizer $\mathfrak U =
C_{\mathfrak G}(\mathfrak u) = \mathfrak A_1 = \langle \mathfrak
q, \mathfrak y, \mathfrak s\rangle$. Furthermore, $\mathfrak U$
has $7$ conjugacy classes of involutions $2_i$, $1 \le i \le 7 $,
with representatives given in Table \ref{Fi_24 cc A_1}, and
$|\mathfrak U| = 2^{19}\cdot 3^9\cdot 5^2\cdot 7\cdot 11\cdot 13$.

\item[\rm(c)] The Fitting subgroup $\mathfrak B$ of $\mathfrak E =
\langle \mathfrak x, \mathfrak y, \mathfrak w, \mathfrak s\rangle$
is an elementary abelian group of order $2^{11}$ such that
$\mathfrak E/\mathfrak B \cong \M_{24}$.

\item[\rm(d)] $\mathfrak D = \langle \mathfrak r_j\mid 1 \le j \le
4\rangle = N_{\mathfrak H}(\mathfrak B) = C_{\mathfrak
E}(\mathfrak z_1)$ where $\mathfrak z_1 = \mathfrak x^2 \in
\mathfrak E$.

\item[\rm(e)] $\mathfrak D_1 = \langle \mathfrak x, \mathfrak y,
\mathfrak s\rangle = N_{\mathfrak A_1}(\mathfrak B) = C_{\mathfrak
E}(\mathfrak z_2)$ where $\mathfrak z_2 = (\mathfrak x\mathfrak
y^2)^7 \in \mathfrak E$.

\item[\rm(f)] $\mathfrak G$ has two conjugacy classes of
involutions represented by $\mathfrak z$ and $\mathfrak u$.

\item[\rm(g)]  The conjugacy classes of involutions $2_1$, $2_3$,
$2_5$, $2_6$, $2_8$ and $2_9$ of $\mathfrak H$ fuse with
$\mathfrak z$ in $\mathfrak G$. Its classes $2_2$, $2_4$ and $2_7$
fuse with $\mathfrak u$ in $\mathfrak G$.

\item[\rm(h)] The conjugacy classes of involutions $2_1$, $2_2$,
$2_3$ and $2_4$ of $\mathfrak U$ fuse with $\mathfrak u$ in
$\mathfrak G$. Its classes $2_5$, $2_6$ and $2_7$ fuse with
$\mathfrak z$ in $\mathfrak G$.

\item[\rm(i)] $\mathfrak G = 2132400816\cdot |\mathfrak U| +
4388805476055\cdot |\mathfrak H| = 2^{21} \cdot 3^{13} \cdot
5^2\cdot 7\cdot 11\cdot 13\cdot 17 \cdot 23\cdot 29$.
\end{enumerate}
\end{proposition}

\begin{proof}
(a) holds by Proposition \ref{prop. 2-centralFi_24}(a), (j), (l)
and Table \ref{Fi_24 cc H}.

(b) By Proposition \ref{prop. E(Fi_24)}(c) we know that $\mathfrak
A_1 = \langle \mathfrak q, \mathfrak y, \mathfrak s\rangle$. Its
conjugacy classes are classified in Table \ref{Fi_24 cc A_1}. It
asserts that the involution $\mathfrak u = (\mathfrak y \mathfrak
s)^{14}$ generates the center of $\mathfrak A_1$. The equation
$C_{\mathfrak G}(\mathfrak u) = \mathfrak A_1$ is a consequence of
Lemma \ref{l. nicepres.Fi_24}(e) and (f). All other assertions
hold by Table \ref{Fi_24 cc A_1}.

(c) and (d) hold by Proposition \ref{prop. 2-centralFi_24}(d) and
(g), respectively.

(e) This is true by Proposition \ref{prop. E(Fi_24)}(d), (f)
and Table \ref{Fi_24 cc E}.

(f) Table \ref{Fi_24 cc E} shows that $E$ has $5$ conjugacy
classes of involutions and that $z_1 = x^2$ is the representative
of the unique $2$-central conjugacy class of $E = \langle
x,y,e\rangle$. By (d) the subgroup $\mathfrak D = \langle
\mathfrak x, \mathfrak y, \mathfrak s\rangle = \mathfrak H \cap
\mathfrak E$. It has $18$ conjugacy classes $2_k$, $1 \le k \le
18$ of involutions. Using MAGMA and the faithful permutation
representations $P\mathfrak H$ and $P\mathfrak E$ of Propositions
\ref{prop. 2-centralFi_24}(c) and \ref{prop. E(Fi_24)}(h),
respectively, we calculated the fusion of the classes $2_k$ of
$\mathfrak D$ in $\mathfrak H$ and also in $\mathfrak E$. Thus we
obtained a fusion graph $\mathcal G(H)$ of the fusion of the
$\mathfrak H$-classes and $\mathfrak E$-classes of order $2$ in
the matrix group $\mathfrak G$. It follows that $\mathfrak G$ has
$2$ conjugacy classes of involutions represented by $\mathfrak z$
and $\mathfrak u$ belonging to the classes $2_1$ and $2_2$ of
$\mathfrak H$, respectively.

(g) This statement follows also from the fusion graph $\mathcal
G$.

(h) In view of (e) we now study the fusion of the $11$ classes of
involutions of $\mathfrak R = N_{\mathfrak A_1}(\mathfrak B) =
\mathfrak A_1 \cap \mathfrak E$ in the two over groups $\mathfrak
A_1$ and $\mathfrak E$. Using MAGMA and the faithful permutation
representation $PA_1$ of Lemma \ref{l. aut2Fi22}(c) it follows
that the classes $2_1$, $2_2$, $2_3$ and $2_4$ of $\mathfrak U$
fuse with $\mathfrak u$, and that the remaining three classes of
$\mathfrak U$ fuse in $\mathfrak G$ with $\mathfrak z$.

(i) In order to simplify notation we replace the Gothic letters by
Roman ones. Let $r(z,u,z)=\left|\left\{(x,y)\in(z^G\cap H) \times
(u^G\cap
H)\big| z\in \langle xy\rangle\right\}\right|$ and\\
$r(z,u,u)=\left|\left\{(x,y)\in(z^G\cap U) \times (u^G\cap U)\big|
u\in \langle xy\rangle\right\}\right|$.

By Table \ref{Fi_24 ct H} $H$ has $45$ real $z$-special conjugacy
classes. Here their representatives are denoted by $\{t_j \mid 1
\le j \le 45\}$. Let $z_i$ and $u_k$ be representatives of the
$H$-classes of involutions fusing to $z$ and $u$ in $G$,
respectively. For each triple $(z_i,u_k,t_j)$ let
$$d(z_i,u_k,t_j) = \frac{|H|^2}{|C_H(z_i)|\cdot
|C_H(u_k)|\cdot|C_H(t_j)|}\cdot\sum_{\psi\in\Irr_{\mathbb{C}}(H)}\psi(z_i)\psi(u_k)\psi(t_j)\psi(1)^{-1}.$$

Then Theorem 1.6.4 of \cite{michler2} and (g) imply that

$$r(z, u, z) = \sum\limits_{i=1}^{6}\sum\limits_{k=1}^{3}\sum\limits_{j=1}^{45}d(z_i,u_k,t_j).$$

Using (g) and the values of the character Table \ref{Fi_24 ct H}
of $H$ these formulas yield that $r(z,u,z) = 2132400816$.

By Table \ref{Fi_24 ct A_1} $U = C_G(u) = A_1$ has $22$ real
$u$-special conjugacy classes. Denote their representatives by
$\{s_n \mid 1 \le n \le 22\}$. Let $u_i$ and $z_k$ be
representatives of the $U$-classes of involutions fusing to $u$
and $z$ in $G$, respectively. For each triple $(u_i,z_k,s_n)$ let
$$d(u_i,z_k,s_n) = \frac{|U|^2}{|C_U(u_i)|\cdot
|C_U(z_k)|\cdot|C_U(s_n)|}\cdot\sum_{\psi\in\Irr_{\mathbb{C}}(U)}\psi(u_i)\psi(z_k)\psi(s_n)\psi(1)^{-1}.$$

Then Theorem 1.6.4 of \cite{michler2} and (h) imply that

$$r(z, u, u) = \sum\limits_{i=1}^{4}\sum\limits_{k=1}^{3}\sum\limits_{n=1}^{22}d(u_i,z_k,s_n).$$

Using (h) and the values of the character Table
\ref{Fi_24 ct A_1} of $U$ these formulas yields that $r(z,u,u) =
4388805476055$. Now Theorem 4.2.1 of \cite{michler} due to J. G.
Thompson implies
$$
|G|=r(z,u,z)\cdot |C_G(u)| + r(z,u,u)\cdot |C_G(z)|\nonumber =
2^{21} \cdot 3^{13} \cdot 5^2\cdot 7\cdot 11\cdot 13\cdot 17 \cdot
23\cdot 29.\nonumber
$$
\end{proof}

\begin{appendix}
\newpage

\section{Representatives of conjugacy classes}

\setlength\textheight{48\baselineskip}
\addtolength\textheight{\topskip}

\begin{cclass}\label{Fi_24 cc E} Conjugacy classes of $E(\Fi_{24}') = \langle x,y,e \rangle$
{ \setlength{\arraycolsep}{1mm}
\renewcommand{\baselinestretch}{0.5}
\scriptsize
$$
\begin{array}{|c|r|r|r|c|c|c|c|c|c|
}\hline
 \mbox{Class} &
 \multicolumn{1}{c|}{\mbox{Representative}} &
 \multicolumn{1}{c|}{|\mbox{Class}|} &
 \multicolumn{1}{c|}{|\mbox{Centralizer}|}
 &
 2{\rm P}
 & 3{\rm P}
 & 5{\rm P}
 & 7{\rm P}
 & 11{\rm P}
 & 23{\rm P}
 \\\hline
 1 & 1 & 1 & 2^{21}\cdot3^{3}\cdot5\cdot7\cdot11\cdot23
 & 1 & 1 & 1 & 1 & 1 & 1 \\\hline

 2_{1} & (e)^{2} & 276 & 2^{19}\cdot3^{2}\cdot5\cdot7\cdot11
 & 1 & 2_{1} & 2_{1} & 2_{1} & 2_{1} & 2_{1} \\\hline

 2_{2} & (x)^{2} & 1771 & 2^{21}\cdot3^{3}\cdot5
 & 1 & 2_{2} & 2_{2} & 2_{2} & 2_{2} & 2_{2} \\\hline

 2_{3} & (y)^{3} & 182160 & 2^{17}\cdot3\cdot7
 & 1 & 2_{3} & 2_{3} & 2_{3} & 2_{3} & 2_{3} \\\hline

 2_{4} & (ye^{2})^{3} & 1275120 & 2^{17}\cdot3
 & 1 & 2_{4} & 2_{4} & 2_{4} & 2_{4} & 2_{4} \\\hline

 2_{5} & (xey)^{3} & 2040192 & 2^{14}\cdot3\cdot5
 & 1 & 2_{5} & 2_{5} & 2_{5} & 2_{5} & 2_{5} \\\hline

 3_{1} & (y)^{2} & 14508032 & 2^{8}\cdot3^{3}\cdot5
 & 3_{1} & 1 & 3_{1} & 3_{1} & 3_{1} & 3_{1} \\
\hline
\end{array}
$$

 \newpage
 \normalsize
 Conjugacy classes of $E(\Fi_{24}') = \langle x,y,e \rangle$
 (continued)
 \scriptsize
 $$
 \begin{array}{|c|r|r|r|c|c|c|c|c|c|
}\hline
 \mbox{Class} &
 \multicolumn{1}{c|}{\mbox{Representative}} &
 \multicolumn{1}{c|}{|\mbox{Class}|} &
 \multicolumn{1}{c|}{|\mbox{Centralizer}|}
 &
 2{\rm P}
 & 3{\rm P}
 & 5{\rm P}
 & 7{\rm P}
 & 11{\rm P}
 & 23{\rm P}
 \\\hline

 21_{1} & xeye & 23876075520 & 3\cdot7
 & 21_{1} & 7_{1} & 21_{2} & 3_{2} & 21_{1} & 21_{1} \\\hline

 21_{2} & xy^{2}eye & 23876075520 & 3\cdot7
 & 21_{2} & 7_{2} & 21_{1} & 3_{2} & 21_{2} & 21_{2} \\\hline

 22 & ye & 22790799360 & 2\cdot11
 & 11 & 22 & 22 & 22 & 2_{1} & 22 \\\hline

 23_{1} & xy^{2}e & 21799895040 & 23
 & 23_{1} & 23_{1} & 23_{2} & 23_{2} & 23_{2} & 1 \\\hline

 23_{2} & xey^{3} & 21799895040 & 23
 & 23_{2} & 23_{2} & 23_{1} & 23_{1} & 23_{1} & 1 \\\hline

 24_{1} & xye & 10445783040 & 2^{4}\cdot3
 & 12_{2} & 8_{1} & 24_{2} & 24_{2} & 24_{1} & 24_{1} \\\hline

 24_{2} & xy^{3} & 10445783040 & 2^{4}\cdot3
 & 12_{2} & 8_{1} & 24_{1} & 24_{1} & 24_{2} & 24_{2} \\\hline

 24_{3} & xey^{2} & 20891566080 & 2^{3}\cdot3
 & 12_{3} & 8_{2} & 24_{3} & 24_{3} & 24_{3} & 24_{3} \\\hline

 24_{4} & xyxy^{2} & 20891566080 & 2^{3}\cdot3
 & 12_{4} & 8_{3} & 24_{4} & 24_{4} & 24_{4} & 24_{4} \\\hline

 28_{1} & xyey^{2} & 17907056640 & 2^{2}\cdot7
 & 14_{1} & 28_{2} & 28_{2} & 4_{5} & 28_{1} & 28_{1} \\\hline

 28_{2} & xye^{2}xe & 17907056640 & 2^{2}\cdot7
 & 14_{2} & 28_{1} & 28_{1} & 4_{5} & 28_{2} & 28_{2} \\\hline

 30_{1} & xye^{2}y & 16713252864 & 2\cdot3\cdot5
 & 15_{1} & 10_{1} & 6_{1} & 30_{2} & 30_{2} & 30_{1} \\\hline

 30_{2} & x^{2}y^{2}xy^{2} & 16713252864 & 2\cdot3\cdot5
 & 15_{2} & 10_{1} & 6_{1} & 30_{1} & 30_{1} & 30_{2} \\\hline

 \end{array}
 $$
}
\end{cclass}

\normalsize
\newpage

\begin{cclass}\label{Fi_24 cc A_1} Conjugacy classes of $A_1 = Aut(2\Fi_{22}) = \langle y,q,s\rangle$
{ \setlength{\arraycolsep}{1mm}
\renewcommand{\baselinestretch}{0.5}
\scriptsize
$$
 \begin{array}{|c|r|r|c|c|c|c|c|c|
}\hline
 \mbox{Class} &
 \multicolumn{1}{c|}{\mbox{Representative}} &
 \multicolumn{1}{c|}{|\mbox{Centralizer}|}
 &
 2{\rm P}
 & 3{\rm P}
 & 5{\rm P}
 & 7{\rm P}
 & 11{\rm P}
 & 13{\rm P}
 \\\hline
 1 & 1 & 2^{19}\cdot3^{9}\cdot5^{2}\cdot7\cdot11\cdot13
 & 1 & 1 & 1 & 1 & 1 & 1 \\\hline

 2_{1} & (ys)^{14} & 2^{19}\cdot3^{9}\cdot5^{2}\cdot7\cdot11\cdot13
 & 1 & 2_{1} & 2_{1} & 2_{1} & 2_{1} & 2_{1} \\\hline

 2_{2} & (y)^{7} & 2^{17}\cdot3^{6}\cdot5\cdot7\cdot11
 & 1 & 2_{2} & 2_{2} & 2_{2} & 2_{2} & 2_{2} \\\hline

 2_{3} & s & 2^{14}\cdot3^{6}\cdot5^{2}\cdot7
 & 1 & 2_{3} & 2_{3} & 2_{3} & 2_{3} & 2_{3} \\\hline

 2_{4} & (yq)^{10} & 2^{19}\cdot3^{4}\cdot5
 & 1 & 2_{4} & 2_{4} & 2_{4} & 2_{4} & 2_{4} \\\hline

 2_{5} & (y^{2}q)^{8} & 2^{19}\cdot3^{4}\cdot5
 & 1 & 2_{5} & 2_{5} & 2_{5} & 2_{5} & 2_{5} \\\hline

 2_{6} & (y^{3}qs)^{9} & 2^{14}\cdot3^{4}\cdot5
 & 1 & 2_{6} & 2_{6} & 2_{6} & 2_{6} & 2_{6} \\
\hline

\end{array}
$$

\newpage
 \normalsize
 Conjugacy classes of $A_1 = Aut(2\Fi_{22}) = \langle y,q,s \rangle$
 (continued)
 \scriptsize

 $$
 \begin{array}{|c|r|r|c|c|c|c|c|c|
}\hline
 \mbox{Class} &
 \multicolumn{1}{c|}{\mbox{Representative}} &
 \multicolumn{1}{c|}{|\mbox{Centralizer}|}
 &
 2{\rm P}
 & 3{\rm P}
 & 5{\rm P}
 & 7{\rm P}
 & 11{\rm P}
 & 13{\rm P}
 \\\hline
 6_{29} & y^{3}qsqs & 2^{5}\cdot3^{3}
 & 3_{4} & 2_{7} & 6_{29} & 6_{29} & 6_{29} & 6_{29} \\\hline

 6_{30} & y^{4}qy^{3}q & 2^{5}\cdot3^{3}
 & 3_{3} & 2_{7} & 6_{30} & 6_{30} & 6_{30} & 6_{30} \\\hline

 6_{31} & y^{3}sysqys & 2^{5}\cdot3^{3}
 & 3_{3} & 2_{6} & 6_{31} & 6_{31} & 6_{31} & 6_{31} \\\hline

 7 & (y)^{2} & 2^{3}\cdot3\cdot7
 & 7 & 7 & 7 & 1 & 7 & 7 \\\hline

 8_{1} & y^{2}s & 2^{8}\cdot3
 & 4_{6} & 8_{1} & 8_{1} & 8_{1} & 8_{1} & 8_{1} \\\hline

 8_{2} & (yqyq^{2})^{3} & 2^{8}\cdot3
 & 4_{4} & 8_{2} & 8_{2} & 8_{2} & 8_{2} & 8_{2} \\\hline

 8_{3} & (y^{4}qs)^{3} & 2^{8}\cdot3
 & 4_{4} & 8_{3} & 8_{3} & 8_{3} & 8_{3} & 8_{3} \\
\hline

\end{array}
$$

\newpage
 \normalsize
 Conjugacy classes of $A_1 = Aut(2\Fi_{22}) = \langle y,q,s \rangle$
 (continued)
 \scriptsize

 $$
 \begin{array}{|c|r|r|c|c|c|c|c|c|
}\hline
 \mbox{Class} &
 \multicolumn{1}{c|}{\mbox{Representative}} &
 \multicolumn{1}{c|}{|\mbox{Centralizer}|}
 &
 2{\rm P}
 & 3{\rm P}
 & 5{\rm P}
 & 7{\rm P}
 & 11{\rm P}
 & 13{\rm P}
 \\\hline
 14_{3} & y & 2^{2}\cdot7
 & 7 & 14_{3} & 14_{3} & 2_{2} & 14_{3} & 14_{3} \\\hline

 15 & (y^{3}q^{2})^{2} & 2^{3}\cdot3\cdot5
 & 15 & 5 & 3_{1} & 15 & 15 & 15 \\\hline

 16_{1} & y^{2}q & 2^{5}
 & 8_{6} & 16_{1} & 16_{1} & 16_{1} & 16_{1} & 16_{1} \\\hline

 16_{2} & y^{2}qs & 2^{5}
 & 8_{6} & 16_{2} & 16_{2} & 16_{2} & 16_{2} & 16_{2} \\\hline

 18_{1} & y^{4}syqyqys & 2^{4}\cdot3^{4}
 & 9_{1} & 6_{2} & 18_{1} & 18_{1} & 18_{1} & 18_{1} \\\hline

 18_{2} & (y^{4}qys)^{2} & 2^{3}\cdot3^{4}
 & 9_{2} & 6_{2} & 18_{2} & 18_{2} & 18_{2} & 18_{2} \\\hline

 18_{3} & y^{4}qsq & 2^{3}\cdot3^{3}
 & 9_{1} & 6_{14} & 18_{6} & 18_{3} & 18_{6} & 18_{3} \\
\hline

 \end{array}
 $$
 }
 \end{cclass}


\newpage

\begin{cclass}\label{Fi_24 cc H} Conjugacy classes of $H(\Fi_{24}') = \langle r, p, b \rangle$
{ \setlength{\arraycolsep}{1mm}
\renewcommand{\baselinestretch}{0.5}
\scriptsize
$$
 \begin{array}{|c|c|r|r|c|c|c|c|}
 \hline
  \mbox{No.} &
 \mbox{Class} &
 \multicolumn{1}{c|}{\mbox{Representative}} &
 \multicolumn{1}{c|}{|\mbox{Centralizer}|}
 & 2{\rm P}
 & 3{\rm P}
 & 5{\rm P}
 & 7{\rm P}
 \\\hline
 1 & 1 &  1 & 2^{21} \cdot 3^{7} \cdot 5^{1} \cdot 7^{1}
 & 1 & 1 & 1 & 1\\\hline
 2 & 2_{1} & z=(p^4 r^2)^6  & 2^{21} \cdot 3^{7} \cdot 5^{1} \cdot 7^{1}
 & 1 & 2_{1} & 2_{1} & 2_{1}\\\hline
 3 & 2_{2} & (p r)^6  & 2^{19} \cdot 3^{4} \cdot 5^{1}
 & 1 & 2_{2} & 2_{2} & 2_{2}\\\hline
 4 & 2_{3} & r^6 & 2^{20} \cdot 3^{2} \cdot 5^{1}
 & 1 & 2_{3} & 2_{3} & 2_{3}\\\hline
 5 & 2_{4} & (p^2 b r)^9 & 2^{14} \cdot 3^{4} \cdot 5^{1}
 & 1 & 2_{4} & 2_{4} & 2_{4}\\\hline
 6 & 2_{5} & z (p^2 b r)^9 & 2^{14} \cdot 3^{4} \cdot 5^{1}
 & 1 & 2_{5} & 2_{5} & 2_{5}\\\hline
 7 & 2_{6} & p^6 & 2^{17} \cdot 3^{3}
 & 1 & 2_{6} & 2_{6} & 2_{6}\\
\hline
\end{array}
$$
 \newpage
 \normalsize
 Conjugacy classes of $H(\Fi_{24}') = \langle r, p, b \rangle$
 (continued)
 \scriptsize
$$
 \begin{array}{|c|c|r|r|c|c|c|c|}
 \hline
 \mbox{No} &
 \mbox{Class} &
 \multicolumn{1}{c|}{\mbox{Representative}} &
 \multicolumn{1}{c|}{|\mbox{Centralizer}|}
 & 2{\rm P}
 & 3{\rm P}
 & 5{\rm P}
 & 7{\rm P}
 \\\hline
 56 & 6_{18} & z b & 2^{8} \cdot 3^{3}
 & 3_{5} & 2_{6} & 6_{18} & 6_{18}\\\hline
 57 & 6_{19} & z (p^2 r b p b r^2 b p r) & 2^{8} \cdot 3^{3}
 & 3_{2} & 2_{7} & 6_{19} & 6_{19}\\\hline
 58 & 6_{20} & (p^3 b r)^4 & 2^{6} \cdot 3^{4}
 & 3_{7} & 2_{1} & 6_{20} & 6_{20}\\\hline
 59 & 6_{21} & (p r)^2 & 2^{7} \cdot 3^{3}
 & 3_{5} & 2_{2} & 6_{21} & 6_{21}\\\hline
 60 & 6_{22} & (p^2 b r)^3 & 2^{5} \cdot 3^{4}
 & 3_{4} & 2_{4} & 6_{23} & 6_{22}\\\hline
 61 & 6_{23} & (p^2 b r)^{15} & 2^{5} \cdot 3^{4}
 & 3_{4} & 2_{4} & 6_{22} & 6_{23}\\\hline
 62 & 6_{24} & z (p^2 b r)^{15} & 2^{5} \cdot 3^{4}
 & 3_{4} & 2_{5} & 6_{25} & 6_{24}\\
\hline
 \end{array}
$$
 \newpage
 \normalsize
 Conjugacy classes of $H(\Fi_{24}') = \langle r, p, b \rangle$
 (continued)
 \scriptsize
  $$
 \begin{array}{|c|c|r|r|c|c|c|c|}
 \hline
 \mbox{No} &
 \mbox{Class} &
 \multicolumn{1}{c|}{\mbox{Representative}} &
 \multicolumn{1}{c|}{|\mbox{Centralizer}|}
 & 2{\rm P}
 & 3{\rm P}
 & 5{\rm P}
 & 7{\rm P}
 \\\hline
 112 & 12_{13} & p^4 r p r b & 2^{6} \cdot 3^{2}
 & 6_{7} & 4_{8} & 12_{13} & 12_{13}\\\hline
 113 & 12_{14} & (p r b p b^2)^2 & 2^{6} \cdot 3^{2}
 & 6_{5} & 4_{2} & 12_{14} & 12_{14}\\\hline
 114 & 12_{15} & z p & 2^{6} \cdot 3^{2}
 & 6_{14} & 4_{9} & 12_{15} & 12_{15}\\\hline
 115 & 12_{16} & p^2 r^2 b p r & 2^{7} \cdot 3^{1}
 & 6_{9} & 4_{3} & 12_{16} & 12_{16}\\\hline
 116 & 12_{17} & p^2 r p^2 b r & 2^{7} \cdot 3^{1}
 & 6_{9} & 4_{6} & 12_{17} & 12_{17}\\\hline
 117 & 12_{18} & z (p^2 b^2 p b^2) & 2^{5} \cdot 3^{2}
 & 6_{15} & 4_{7} & 12_{18} & 12_{18}\\\hline
 118 & 12_{19} & p^2 r p^2 b^2 r & 2^{5} \cdot 3^{2}
 & 6_{6} & 4_{8} & 12_{19} & 12_{19}\\
\hline
\end{array}
$$
}
\end{cclass}

\setlength\textheight{50\baselineskip}
\addtolength\textheight{\topskip}

\newpage

\section{Character tables}

\begin{ctab}\label{Fi_24 ct_E}Character table of $E(\Fi'_{24}) = \langle x,y,e \rangle$
{ \setlength{\arraycolsep}{1mm}
\renewcommand{\baselinestretch}{0.5}

{\tiny \setlength{\arraycolsep}{0,3mm}
$$
\begin{array}{r|rrrrrrrrrrrrrrrrrrrrr}
2&21&19&21&17&17&14&8&6&15&15&15&15&11&13&11&11&11&11&10&10\\
3&3&2&3&1&1&1&3&2&1&1&1&1&1&.&1&1&.&.&.&.\\
5&1&1&1&.&.&1&1&.&1&1&.&.&.&.&.&.&.&.&.&.\\
7&1&1&.&1&.&.&.&1&.&.&.&.&1&.&.&.&.&.&.&.\\
11&1&1&.&.&.&.&.&.&.&.&.&.&.&.&.&.&.&.&.&.\\
23&1&.&.&.&.&.&.&.&.&.&.&.&.&.&.&.&.&.&.&.\\
\hline
&1a&2a&2b&2c&2d&2e&3a&3b&4a&4b&4c&4d&4e&4f&4g&4h&4i&4j&4k&4l\\
\hline
2P&1a&1a&1a&1a&1a&1a&3a&3b&2b&2b&2b&2b&2a&2b&2c&2c&2d&2d&2c&2d\\
\end{array}
$$
} }

\newpage
Character table of $E(Fi'_{24})$ (continued){
\renewcommand{\baselinestretch}{0.5}

{\tiny \setlength{\arraycolsep}{0,3mm}
$$
\begin{array}{r|rrrrrrrrrrrrrrrrrrrrrrrrrrrr}
2&5&8&8&8&6&6&6&4&3&3&8&8&8&9&8&7&6&6&5&4&4&4&1&6&6&5&5\\
3&1&3&2&2&2&1&1&1&1&1&1&1&1&.&.&.&.&.&1&.&.&.&.&1&1&1&1\\
5&1&1&.&.&.&.&.&.&.&.&.&.&.&.&.&.&.&.&1&1&1&1&.&.&.&.&.\\
7&.&.&.&.&.&.&.&.&1&1&.&.&.&.&.&.&.&.&.&.&.&.&.&.&.&.&.\\
11&.&.&.&.&.&.&.&.&.&.&.&.&.&.&.&.&.&.&.&.&.&.&1&.&.&.&.\\
23&.&.&.&.&.&.&.&.&.&.&.&.&.&.&.&.&.&.&.&.&.&.&.&.&.&.&.\\
\hline
&5a&6a&6b&6c&6d&6e&6f&6g&7a&7b&8a&8b&8c&8d&8e&8f&8g&8h&10a&10b
&10c&10d&11a&12a&12b&12c&12d\\
\hline
2P&5a&3a&3a&3a&3b&3a&3a&3b&7a&7b&4d&4a&4b&4d&4c&4f&4i&4j&5a&5a&5a
&5a&11a&6b&6b&6d&6d\\
\end{array}
$$
} }

\newpage
Character table of $E(Fi'_{24})$ (continued){
\renewcommand{\baselinestretch}{0.5}

{\tiny \setlength{\arraycolsep}{0,3mm}
$$
\begin{array}{r|rrrrrrrrrrrrrrrrrrrrrrrrrr}
2&4&4&4&3&3&2&2&1&1&5&4&4&.&.&1&.&.&4&4&3&3&2&2&1&1\\
3&1&1&1&.&.&.&.&1&1&.&.&.&1&1&.&.&.&1&1&1&1&.&.&1&1\\
5&.&.&.&.&.&.&.&1&1&.&1&1&.&.&.&.&.&.&.&.&.&.&.&1&1\\
7&.&.&.&1&1&1&1&.&.&.&.&.&1&1&.&.&.&.&.&.&.&1&1&.&.\\
11&.&.&.&.&.&.&.&.&.&.&.&.&.&.&1&.&.&.&.&.&.&.&.&.&.\\
23&.&.&.&.&.&.&.&.&.&.&.&.&.&.&.&1&1&.&.&.&.&.&.&.&.\\
\hline
&12e&12f&12g&14a&14b&14c&14d&15a&15b&16a&20a&20b&21a&21b&22a&23a
&23b&24a&24b&24c&24d&28a&28b&30a&30b\\
\hline
2P&6c&6e&6e&7b&7a&7a&7b&15a&15b&8d&10a&10a&21a&21b&11a&23a&23b&12b
&12b&12c&12d&14a&14b&15a&15b\\
\end{array}
$$
} } where $ A = \zeta(7)^4 + \zeta(7)^2 + \zeta(7)$, $B =
-2\zeta(15)_3\zeta(15)_5^3 - 2\zeta(15)_3\zeta(15)_5^2 -
\zeta(15)_3 - \zeta(15)_5^3 - \zeta(15)_5^2 - 1$, $C =
-\zeta(23)^{18} - \zeta(23)^{16} - \zeta(23)^{13} - \zeta(23)^{12}
- \zeta(23)^9 - \zeta(23)^8 - \zeta(23)^6 - \zeta(23)^4 -
\zeta(23)^3 - \zeta(23)^2 - \zeta(23) - 1$, $D = -2\zeta(7)^4 -
2\zeta(7)^2 - 2\zeta(7) - 2$, $E = -3\zeta(7)^4 - 3\zeta(7)^2 -
3\zeta(7) - 3$, $F = -4\zeta(5)^3 - 4\zeta(5)^2 - 3$, $G =
4\zeta(5)^3 + 4\zeta(5)^2 + 1$, $H = -4\zeta(12)_4\zeta(12)_3 -
2\zeta(12)_4$.
\end{ctab}



\newpage

\begin{ctab}\label{Fi_24 ct mH}Character table of $mH = \langle r, u, v\rangle$
{ \setlength{\arraycolsep}{1mm}
\renewcommand{\baselinestretch}{0.5}

{\tiny \setlength{\arraycolsep}{0,3mm}
$$
\begin{array}{r|rrrrrrrrrrrrrrrrrr}
2&13&10&11&13&10&8&7&7&6&4&4&4&2&3&9&10&8\\
3&13&9&6&5&5&13&10&10&7&11&10&9&9&7&4&4&3\\
5&2&1&1&.&1&.&1&1&.&.&.&.&.&.&1&.&.\\
7&1&1&1&.&.&.&.&.&1&.&.&.&.&.&.&.&1\\
13&1&1&.&.&.&.&.&.&1&.&.&.&.&.&.&.&.\\
\hline
&1a&2a&2b&2c&2d&3a&3b&3c&3d&3e&3f&3g&3h&3i&4a&4b&4c\\
\hline
2P&1a&1a&1a&1a&1a&3a&3b&3c&3d&3e&3f&3g&3h&3i&2b&2c&2b\\
3P&1a&2a&2b&2c&2d&1a&1a&1a&1a&1a&1a&1a&1a&1a&4a&4b&4c\\
\end{array}
$$
} }

\newpage
Character table of $mH$ (continued){
\renewcommand{\baselinestretch}{0.5}

{\tiny \setlength{\arraycolsep}{0,3mm}
$$
\begin{array}{r|rrrrrrrrrrrrrrrrrrrrrr}
2&10&7&9&8&4&6&5&7&7&7&8&4&4&6&6&6&7&7&3&6&4\\
3&2&2&2&2&2&9&7&7&5&6&5&7&7&5&5&5&4&4&6&4&5\\
5&.&1&.&.&2&.&1&.&1&.&.&.&.&.&.&.&.&.&.&.&.\\
7&.&.&.&.&.&.&.&.&.&.&.&.&.&.&.&.&.&.&.&.&.\\
13&.&.&.&.&.&.&.&.&.&.&.&.&.&.&.&.&.&.&.&.&.\\
\hline
&4d&4e&4f&4g&5a&6_{1}&6_{2}&6_{3}&6_{4}&6_{5}&6_{6}&6_{7}&6_{8}
&6_{9}&6_{10}&6_{11}&6_{12}&6_{13}&6_{14}&6_{15}&6_{16}\\
\hline
2P&2c&2b&2b&2c&5a&3a&3b&3c&3c&3a&3a&3e&3f&3a&3b&3c&3b&3c&3g&3c&3e\\
3P&4d&4e&4f&4g&5a&2a&2a&2a&2d&2b&2c&2a&2a&2d&2b&2b&2c&2c&2a&2d&2c\\
\end{array}
$$
} }

\newpage
Character table of $mH$ (continued){
\renewcommand{\baselinestretch}{0.5}

{\tiny \setlength{\arraycolsep}{0,3mm}
$$
\begin{array}{r|rrrrrrrrrrrrrrrrrrrrrrrrrr}
2&4&4&4&2&5&6&4&4&4&2&3&2&3&3&6&6&6&6&3&3&1&3&1&1&.\\
3&5&5&5&6&4&3&4&4&4&5&4&4&3&1&1&1&1&.&6&6&7&5&6&6&6\\
5&.&.&.&.&.&.&.&.&.&.&.&.&.&.&.&.&.&.&.&.&.&.&.&.&.\\
7&.&.&.&.&.&.&.&.&.&.&.&.&.&1&.&.&.&.&.&.&.&.&.&.&.\\
13&.&.&.&.&.&.&.&.&.&.&.&.&.&.&.&.&.&.&.&.&.&.&.&.&.\\
\hline
&6_{17}&6_{18}&6_{19}&6_{20}&6_{21}&6_{22}&6_{23}&6_{24}&6_{25}
&6_{26}&6_{27}&6_{28}&6_{29}&7a&8a&8b&8c&8d&9a&9b&9c&9d&9e&9f&9g\\
\hline
2P&3g&3e&3f&3h&3b&3d&3g&3f&3g&3h&3g&3h&3i&7a&4b&4d&4b&4d&9a&9b&9c
&9d&9e&9f&9g\\
3P&2c&2d&2d&2a&2d&2c&2b&2c&2c&2c&2d&2d&2c&7a&8a&8b&8c&8d&3a&3a&3a
&3e&3a&3e&3e\\
\end{array}
$$
} }

\newpage
Character table of $mH$ (continued){
\renewcommand{\baselinestretch}{0.5}

{\tiny \setlength{\arraycolsep}{0,3mm}
$$
\begin{array}{r|rrrrrrrrrrrrrrrrrrrrrrrrrr}
2&1&1&1&4&4&4&7&5&7&5&5&5&5&3&4&5&5&5&5&5&3&4&4&3&3\\
3&5&5&4&1&1&.&3&4&2&3&3&3&3&4&3&2&2&2&2&2&3&2&2&2&2\\
5&.&.&.&1&1&1&.&.&.&.&.&.&.&.&.&.&.&.&.&.&.&.&.&.&.\\
7&.&.&.&.&.&.&.&.&.&.&.&.&.&.&.&.&.&.&.&.&.&.&.&.&.\\
13&.&.&.&.&.&.&.&.&.&.&.&.&.&.&.&.&.&.&.&.&.&.&.&.&.\\
\hline
&9h&9i&9j&10a&10b&10c&12a&12b&12c&12d&12e&12f&12g&12h&12i&12j
&12k&12l&12m&12n&12o&12p&12q&12r&12s\\
\hline
2P&9h&9i&9j&5a&5a&5a&6_{6}&6_{5}&6_{6}&6_{12}&6_{13}&6_{10}&6_{5}
&6_{16}&6_{11}&6_{22}&6_{22}&6_{10}&6_{6}&6_{13}&6_{24}&6_{10}&6_{11}
&6_{16}&6_{23}\\
3P&3e&3e&3f&10a&10b&10c&4b&4a&4d&4b&4b&4a&4c&4b&4a&4d&4b&4f&4g&4g
&4b&4e&4e&4g&4c\\
\end{array}
$$
} }

\newpage
Character table of $mH$ (continued){
\renewcommand{\baselinestretch}{0.5}

{\tiny \setlength{\arraycolsep}{0,3mm}
$$
\begin{array}{r|rrrrrrrrrrrrrrrrrrrrrrrrrr}
2&3&3&2&1&1&3&2&1&1&3&2&3&3&3&3&1&1&1&1&1&1&3&3&.&3\\
3&2&2&2&1&1&.&.&2&2&4&4&3&3&3&3&3&3&3&3&3&3&.&.&1&1\\
5&.&.&.&.&.&.&.&1&1&.&.&.&.&.&.&.&.&.&.&.&.&1&1&.&.\\
7&.&.&.&.&.&1&1&.&.&.&.&.&.&.&.&.&.&.&.&.&.&.&.&1&.\\
13&.&.&.&1&1&.&.&.&.&.&.&.&.&.&.&.&.&.&.&.&.&.&.&.&.\\
\hline
&12t&12u&12v&13a&13b&14a&14b&15a&15b&18a&18b&18c&18d&18e&18f&18g
&18h&18i&18j&18k&18l&20a&20b&21a&24a\\
\hline
2P&6_{24}&6_{23}&6_{29}&13b&13a&7a&7a&15a&15b&9b&9a&9d&9a&9b&9b&9h
&9e&9f&9i&9j&9c&10c&10c&21a&12d\\
3P&4g&4f&4d&13a&13b&14a&14b&5a&5a&6_{1}&6_{1}&6_{16}&6_{5}&6_{9}
&6_{5}&6_{16}&6_{9}&6_{16}&6_{7}&6_{8}&6_{9}&20a&20b&7a&8a\\
\end{array}
$$
} }

\newpage
Character table of $mH$ (continued){
\renewcommand{\baselinestretch}{0.5}

{\tiny \setlength{\arraycolsep}{0,3mm}
$$
\begin{array}{r|rrrrrrrrrrrrr}
2&3&3&1&1&.&2&1&1&2&2&.&.\\
3&1&1&.&.&3&.&1&1&2&2&1&1\\
5&.&.&.&.&.&.&1&1&.&.&.&.\\
7&.&.&.&.&.&1&.&.&.&.&.&.\\
13&.&.&1&1&.&.&.&.&.&.&1&1\\
\hline
&24b&24c&26a&26b&27a&28a&30a&30b&36a&36b&39a&39b\\
\hline
2P&12j&12k&13a&13b&27a&14a&15a&15b&18c&18d&39b&39a\\
3P&8b&8c&26a&26b&9c&28a&10a&10b&12h&12b&13a&13b\\
\end{array}
$$
} }

\newpage
Character table of $mH$ (continued){
\renewcommand{\baselinestretch}{0.5}

{\tiny \setlength{\arraycolsep}{0,3mm}
$$
\begin{array}{r|rrrrrrrrrrrrrrrrr}
2&13&10&11&13&10&8&7&7&6&4&4&4&2&3&9&10\\
3&13&9&6&5&5&13&10&10&7&11&10&9&9&7&4&4\\
5&2&1&1&.&1&.&1&1&.&.&.&.&.&.&1&.\\
7&1&1&1&.&.&.&.&.&1&.&.&.&.&.&.&.\\
13&1&1&.&.&.&.&.&.&1&.&.&.&.&.&.&.\\
\hline
&1a&2a&2b&2c&2d&3a&3b&3c&3d&3e&3f&3g&3h&3i&4a&4b\\
\hline
2P&1a&1a&1a&1a&1a&3a&3b&3c&3d&3e&3f&3g&3h&3i&2b&2c\\
3P&1a&2a&2b&2c&2d&1a&1a&1a&1a&1a&1a&1a&1a&1a&4a&4b\\
\end{array}
$$
} }

\newpage
Character table of $mH$ (continued){
\renewcommand{\baselinestretch}{0.5}

{\tiny \setlength{\arraycolsep}{0,3mm}
$$
\begin{array}{r|rrrrrrrrrrrrrrrrrrrrrr}
2&8&10&7&9&8&4&6&5&7&7&7&8&4&4&6&6&6&7&7&3&6\\
3&3&2&2&2&2&2&9&7&7&5&6&5&7&7&5&5&5&4&4&6&4\\
5&.&.&1&.&.&2&.&1&.&1&.&.&.&.&.&.&.&.&.&.&.\\
7&1&.&.&.&.&.&.&.&.&.&.&.&.&.&.&.&.&.&.&.&.\\
13&.&.&.&.&.&.&.&.&.&.&.&.&.&.&.&.&.&.&.&.&.\\
\hline
&4c&4d&4e&4f&4g&5a&6_{1}&6_{2}&6_{3}&6_{4}&6_{5}&6_{6}&6_{7}
&6_{8}&6_{9}&6_{10}&6_{11}&6_{12}&6_{13}&6_{14}&6_{15}\\
\hline
2P&2b&2c&2b&2b&2c&5a&3a&3b&3c&3c&3a&3a&3e&3f&3a&3b&3c&3b&3c&3g&3c\\
3P&4c&4d&4e&4f&4g&5a&2a&2a&2a&2d&2b&2c&2a&2a&2d&2b&2b&2c&2c&2a&2d\\
\end{array}
$$
} }

\newpage
Character table of $mH$ (continued){
\renewcommand{\baselinestretch}{0.5}

{\tiny \setlength{\arraycolsep}{0,3mm}
$$
\begin{array}{r|rrrrrrrrrrrrrrrrrrrrrrrrrr}
2&4&4&4&4&2&5&6&4&4&4&2&3&2&3&3&6&6&6&6&3&3&1&3&1&1\\
3&5&5&5&5&6&4&3&4&4&4&5&4&4&3&1&1&1&1&.&6&6&7&5&6&6\\
5&.&.&.&.&.&.&.&.&.&.&.&.&.&.&.&.&.&.&.&.&.&.&.&.&.\\
7&.&.&.&.&.&.&.&.&.&.&.&.&.&.&1&.&.&.&.&.&.&.&.&.&.\\
13&.&.&.&.&.&.&.&.&.&.&.&.&.&.&.&.&.&.&.&.&.&.&.&.&.\\
\hline
&6_{16}&6_{17}&6_{18}&6_{19}&6_{20}&6_{21}&6_{22}&6_{23}&6_{24}
&6_{25}&6_{26}&6_{27}&6_{28}&6_{29}&7a&8a&8b&8c&8d&9a&9b&9c&9d&9e&9f\\
\hline
2P&3e&3g&3e&3f&3h&3b&3d&3g&3f&3g&3h&3g&3h&3i&7a&4b&4d&4b&4d&9a&9b
&9c&9d&9e&9f\\
3P&2c&2c&2d&2d&2a&2d&2c&2b&2c&2c&2c&2d&2d&2c&7a&8a&8b&8c&8d&3a&3a
&3a&3e&3a&3e\\
\end{array}
$$
} }

\newpage
Character table of $mH$ (continued){
\renewcommand{\baselinestretch}{0.5}

{\tiny \setlength{\arraycolsep}{0,3mm}
$$
\begin{array}{r|rrrrrrrrrrrrrrrrrrrrrrrrrr}
2&.&1&1&1&4&4&4&7&5&7&5&5&5&5&3&4&5&5&5&5&5&3&4&4&3\\
3&6&5&5&4&1&1&.&3&4&2&3&3&3&3&4&3&2&2&2&2&2&3&2&2&2\\
5&.&.&.&.&1&1&1&.&.&.&.&.&.&.&.&.&.&.&.&.&.&.&.&.&.\\
7&.&.&.&.&.&.&.&.&.&.&.&.&.&.&.&.&.&.&.&.&.&.&.&.&.\\
13&.&.&.&.&.&.&.&.&.&.&.&.&.&.&.&.&.&.&.&.&.&.&.&.&.\\
\hline
&9g&9h&9i&9j&10a&10b&10c&12a&12b&12c&12d&12e&12f&12g&12h&12i&12j
&12k&12l&12m&12n&12o&12p&12q&12r\\
\hline
2P&9g&9h&9i&9j&5a&5a&5a&6_{6}&6_{5}&6_{6}&6_{12}&6_{13}&6_{10}
&6_{5}&6_{16}&6_{11}&6_{22}&6_{22}&6_{10}&6_{6}&6_{13}&6_{24}&6_{10}
&6_{11}&6_{16}\\
3P&3e&3e&3e&3f&10a&10b&10c&4b&4a&4d&4b&4b&4a&4c&4b&4a&4d&4b&4f&4g
&4g&4b&4e&4e&4g\\
\end{array}
$$
} }

\newpage
Character table of $mH$ (continued){
\renewcommand{\baselinestretch}{0.5}

{\tiny \setlength{\arraycolsep}{0,3mm}
$$
\begin{array}{r|rrrrrrrrrrrrrrrrrrrrrrrrrr}
2&3&3&3&2&1&1&3&2&1&1&3&2&3&3&3&3&1&1&1&1&1&1&3&3&.\\
3&2&2&2&2&1&1&.&.&2&2&4&4&3&3&3&3&3&3&3&3&3&3&.&.&1\\
5&.&.&.&.&.&.&.&.&1&1&.&.&.&.&.&.&.&.&.&.&.&.&1&1&.\\
7&.&.&.&.&.&.&1&1&.&.&.&.&.&.&.&.&.&.&.&.&.&.&.&.&1\\
13&.&.&.&.&1&1&.&.&.&.&.&.&.&.&.&.&.&.&.&.&.&.&.&.&.\\
\hline
&12s&12t&12u&12v&13a&13b&14a&14b&15a&15b&18a&18b&18c&18d&18e&18f
&18g&18h&18i&18j&18k&18l&20a&20b&21a\\
\hline
2P&6_{23}&6_{24}&6_{23}&6_{29}&13b&13a&7a&7a&15a&15b&9b&9a&9d&9a
&9b&9b&9h&9e&9f&9i&9j&9c&10c&10c&21a\\
3P&4c&4g&4f&4d&13a&13b&14a&14b&5a&5a&6_{1}&6_{1}&6_{16}&6_{5}
&6_{9}&6_{5}&6_{16}&6_{9}&6_{16}&6_{7}&6_{8}&6_{9}&20a&20b&7a\\
\end{array}
$$
} }

\newpage
Character table of $mH$ (continued){
\renewcommand{\baselinestretch}{0.5}

{\tiny \setlength{\arraycolsep}{0,3mm}
$$
\begin{array}{r|rrrrrrrrrrrrrr}
2&3&3&3&1&1&.&2&1&1&2&2&.&.\\
3&1&1&1&.&.&3&.&1&1&2&2&1&1\\
5&.&.&.&.&.&.&.&1&1&.&.&.&.\\
7&.&.&.&.&.&.&1&.&.&.&.&.&.\\
13&.&.&.&1&1&.&.&.&.&.&.&1&1\\
\hline
&24a&24b&24c&26a&26b&27a&28a&30a&30b&36a&36b&39a&39b\\
\hline
2P&12d&12j&12k&13a&13b&27a&14a&15a&15b&18c&18d&39b&39a\\
3P&8a&8b&8c&26a&26b&9c&28a&10a&10b&12h&12b&13a&13b\\
5P&24a&24b&24c&26b&26a&27a&28a&6_{4}&6_{2}&36a&36b&39b&39a\\
7P&24a&24b&24c&26b&26a&27a&4c&30a&30b&36a&36b&39b&39a\\
13P&24a&24b&24c&2a&2a&27a&28a&30a&30b&36a&36b&3d&3d\\
\hline
X.87&.&.&.&.&.&.&.&.&.&.&.&.&.\\
X.88&.&.&.&.&.&.&.&.&.&.&.&.&.\\
X.89&.&.&.&.&.&.&.&.&.&.&1&.&.\\
X.90&.&.&.&.&.&.&.&.&.&.&-1&.&.\\
X.91&.&-1&-1&-1&-1&.&1&.&.&.&.&1&1\\
X.92&.&-1&-1&1&1&.&-1&.&.&.&.&1&1\\
X.93&.&.&.&.&.&.&.&.&.&.&.&.&.\\
X.94&.&.&.&.&.&.&.&.&.&.&.&.&.\\
X.95&-1&.&.&.&.&.&1&.&1&.&.&.&.\\
X.96&1&.&.&.&.&.&-1&.&-1&.&.&.&.\\
X.97&.&.&.&.&.&.&1&-1&.&.&.&.&.\\
X.98&.&.&.&.&.&.&-1&1&.&.&.&.&.\\
X.99&.&.&.&.&.&1&.&.&.&.&.&.&.\\
X.100&.&.&.&.&.&1&.&.&.&.&.&.&.\\
X.101&.&.&.&.&.&-1&.&.&.&.&.&.&.\\
X.102&.&.&.&.&.&.&.&1&.&.&1&.&.\\
X.103&.&.&.&.&.&.&.&.&-1&.&-1&.&.\\
X.104&.&.&.&.&.&.&.&-1&.&.&-1&.&.\\
X.105&.&.&.&.&.&.&.&.&1&.&1&.&.\\
X.106&.&.&.&B&A&1&.&.&.&.&.&A&B\\
X.107&.&.&.&-B&-A&1&.&.&.&.&.&A&B\\
X.108&.&.&.&-A&-B&1&.&.&.&.&.&B&A\\
X.109&.&.&.&A&B&1&.&.&.&.&.&B&A\\
X.110&.&.&.&.&.&.&.&-1&1&.&1&.&.\\
X.111&.&.&.&.&.&.&.&-1&1&.&.&.&.\\
X.112&.&.&.&.&.&.&.&1&-1&.&-1&.&.\\
X.113&.&.&.&.&.&.&.&1&-1&.&.&.&.\\
X.114&.&.&.&.&.&.&.&.&.&.&.&.&.\\
X.115&.&.&.&.&.&.&.&.&.&.&.&.&.\\
X.116&.&.&.&.&.&.&.&.&.&.&.&.&.\\
X.117&.&.&.&.&.&.&-1&.&.&.&.&.&.\\
X.118&.&.&.&.&.&.&1&.&.&.&.&.&.\\
X.119&.&1&1&.&.&.&.&.&.&.&.&-1&-1\\
X.120&.&.&.&.&.&.&.&.&.&.&-1&.&.\\
X.121&.&.&.&.&.&.&.&.&.&.&1&.&.\\
X.122&.&.&.&.&.&.&.&1&-1&.&.&.&.\\
X.123&.&.&.&.&.&.&.&-1&1&.&.&.&.\\
X.124&.&.&.&.&.&-1&.&.&.&.&.&-B&-A\\
X.125&.&.&.&.&.&-1&.&.&.&.&.&-A&-B\\
\end{array}
$$
} }, where $ A = -\zeta(13)^{11} - \zeta(13)^8 - \zeta(13)^7 -
\zeta(13)^6 - \zeta(13)^5 - \zeta(13)^2 - 1$, $B = -A - 1$, $C =
2B$, $D = 2A$.
\end{ctab}

\normalsize

\newpage

\begin{ctab}\label{Fi_24 ct mE}Character table of $mE = \langle u,v,s \rangle$
{ \setlength{\arraycolsep}{1mm}
\renewcommand{\baselinestretch}{0.5}

{\tiny \setlength{\arraycolsep}{0,3mm}
$$
\begin{array}{r|rrrrrrrrrrrrrrrrrrrrr}
2&10&10&10&10&10&10&10&10&9&6&5&5&7&4&6&4&4&3&3&3\\
3&10&9&7&6&5&4&4&3&10&10&10&8&8&8&8&8&8&8&8&7\\
5&1&1&1&1&1&1&.&.&1&.&.&1&.&1&.&.&.&.&.&.\\
7&1&1&1&1&.&.&.&.&1&.&.&.&.&.&.&.&.&.&.&.\\
13&1&1&.&.&.&.&.&.&1&.&.&.&.&.&.&.&.&.&.&.\\
\hline
&1a&2a&2b&2c&2d&2e&2f&2g&3a&3b&3c&3d&3e&3f&3g&3h&3i&3j&3k&3l\\
\hline
2P&1a&1a&1a&1a&1a&1a&1a&1a&3a&3b&3c&3d&3e&3f&3g&3h&3i&3j&3k&3l\\
3P&1a&2a&2b&2c&2d&2e&2f&2g&1a&1a&1a&1a&1a&1a&1a&1a&1a&1a&1a&1a\\
\end{array}
$$
} }

\normalsize

\newpage
Character table of $mE$ (continued){
\renewcommand{\baselinestretch}{0.5}
{\tiny \setlength{\arraycolsep}{0,3mm}
$$
\begin{array}{r|rrrrrrrrrrrrrrrrrrrrrrrr}
2&2&2&1&7&8&7&8&8&7&8&7&4&9&6&9&5&7&6&5&6&9&4&4\\
3&7&7&7&3&3&2&2&2&2&1&1&2&7&9&5&7&7&7&7&6&4&7&7\\
5&.&.&.&1&.&1&.&.&.&.&.&1&1&.&1&1&.&.&.&.&.&.&.\\
7&.&.&.&.&.&.&.&.&.&.&.&.&1&.&.&.&.&.&.&.&.&.&.\\
13&.&.&.&.&.&.&.&.&.&.&.&.&.&.&.&.&.&.&.&.&.&.&.\\
\hline
&3m&3n&3o&4a&4b&4c&4d&4e&4f&4g&4h&5a&6_{1}&6_{2}&6_{3}&6_{4}
&6_{5}&6_{6}&6_{7}&6_{8}&6_{9}&6_{10}&6_{11}\\
\hline
2P&3m&3n&3o&2d&2f&2d&2f&2f&2d&2f&2d&5a&3a&3b&3a&3d&3e&3b&3c&3b&3a
&3h&3i\\
3P&1a&1a&1a&4a&4b&4c&4d&4e&4f&4g&4h&5a&2b&2a&2d&2a&2a&2b&2b&2c&2f
&2a&2a\\
\end{array}
$$
} }
\normalsize
\newpage
Character table of $mE$ (continued){
\renewcommand{\baselinestretch}{0.5}

{\tiny \setlength{\arraycolsep}{0,3mm}
$$
\begin{array}{r|rrrrrrrrrrrrrrrrrrrrrrr}
2&5&5&6&4&4&7&5&5&5&3&6&6&6&6&7&2&5&5&5&5&5&3\\
3&6&6&5&6&6&4&5&5&5&6&4&4&4&4&3&6&4&4&4&4&4&5\\
5&.&.&.&.&.&.&.&.&.&.&.&.&.&.&.&.&.&.&.&.&.&.\\
7&.&.&.&.&.&.&.&.&.&.&.&.&.&.&.&.&.&.&.&.&.&.\\
13&.&.&.&.&.&.&.&.&.&.&.&.&.&.&.&.&.&.&.&.&.&.\\
\hline
&6_{12}&6_{13}&6_{14}&6_{15}&6_{16}&6_{17}&6_{18}&6_{19}&6_{20}
&6_{21}&6_{22}&6_{23}&6_{24}&6_{25}&6_{26}&6_{27}&6_{28}&6_{29}&6_{30}
&6_{31}&6_{32}&6_{33}\\
\hline
2P&3d&3e&3b&3g&3f&3e&3e&3c&3d&3l&3e&3b&3b&3g&3e&3m&3c&3g&3d&3d&3e
&3l\\
3P&2b&2b&2d&2b&2b&2f&2c&2d&2c&2a&2d&2e&2f&2f&2g&2a&2f&2d&2d&2f&2f
&2b\\
\end{array}
$$
} }

\newpage
Character table of $mE$ (continued){
\renewcommand{\baselinestretch}{0.5}

{\tiny \setlength{\arraycolsep}{0,3mm}
$$
\begin{array}{r|rrrrrrrrrrrrrrrrrrrrrrrrrr}
2&6&6&4&4&4&4&4&2&5&5&5&3&3&3&3&3&4&4&2&2&2&3&3&1&2\\
3&3&3&4&4&4&4&4&5&3&3&3&4&4&4&4&4&3&3&4&4&4&3&3&4&3\\
5&.&.&.&.&.&.&.&.&.&.&.&.&.&.&.&.&.&.&.&.&.&.&.&.&.\\
7&.&.&.&.&.&.&.&.&.&.&.&.&.&.&.&.&.&.&.&.&.&.&.&.&.\\
13&.&.&.&.&.&.&.&.&.&.&.&.&.&.&.&.&.&.&.&.&.&.&.&.&.\\
\hline
&6_{34}&6_{35}&6_{36}&6_{37}&6_{38}&6_{39}&6_{40}&6_{41}&6_{42}
&6_{43}&6_{44}&6_{45}&6_{46}&6_{47}&6_{48}&6_{49}&6_{50}&6_{51}&6_{52}
&6_{53}&6_{54}&6_{55}&6_{56}&6_{57}&6_{58}\\
\hline
2P&3b&3e&3i&3f&3h&3g&3f&3n&3e&3d&3d&3j&3k&3l&3l&3l&3h&3i&3m&3n&3n
&3l&3l&3o&3m\\
3P&2g&2e&2f&2d&2f&2f&2f&2b&2g&2g&2e&2f&2f&2f&2f&2c&2g&2g&2f&2f&2f
&2g&2g&2f&2g\\
\end{array}
$$
} }

\newpage
Character table of $mE$ (continued){
\renewcommand{\baselinestretch}{0.5}

{\tiny \setlength{\arraycolsep}{0,3mm}
$$
\begin{array}{r|rrrrrrrrrrrrrrrrrrrrrrrrrrrr}
2&2&5&5&5&5&3&2&2&1&1&1&.&.&4&4&3&4&3&6&7&7&5&5&6&4&4&4\\
3&1&1&1&.&.&5&5&5&5&4&4&4&4&1&1&1&.&.&3&3&2&3&3&2&3&3&3\\
5&.&.&.&.&.&.&.&.&.&.&.&.&.&1&1&1&1&1&1&.&.&.&.&.&.&.&.\\
7&1&.&.&.&.&.&.&.&.&.&.&.&.&.&.&.&.&.&.&.&.&.&.&.&.&.&.\\
13&.&.&.&.&.&.&.&.&.&.&.&.&.&.&.&.&.&.&.&.&.&.&.&.&.&.&.\\
\hline &7a&8a&8b&8c&8d&9a&9b&9c&9d&9e&9f&9g&9h&10a&10b&10c&10d&10e
&12_{1}&12_{2}&12_{3}&12_{4}&12_{5}&12_{6}&12_{7}&12_{8}&12_{9}\\
\hline
2P&7a&4e&4b&4b&4e&9a&9b&9c&9d&9e&9f&9g&9h&5a&5a&5a&5a&5a&6_{3}
&6_{9}&6_{9}&6_{24}&6_{17}&6_{3}&6_{28}&6_{25}&6_{22}\\
3P&7a&8a&8b&8c&8d&3b&3b&3b&3b&3h&3i&3i&3h&10a&10b&10c&10d&10e&4a
&4b&4e&4b&4b&4f&4b&4b&4a\\
\end{array}
$$
} }

\newpage
Character table of $mE$ (continued){
\renewcommand{\baselinestretch}{0.5}

{\tiny \setlength{\arraycolsep}{0,3mm}
$$
\begin{array}{r|rrrrrrrrrrrrrrrrrrrrr}
2&4&5&5&5&3&3&3&3&4&4&4&4&2&2&5&3&3&3&4&1\\
3&3&2&2&2&3&3&3&3&2&2&2&2&3&3&1&2&2&2&1&1\\
5&.&.&.&.&.&.&.&.&.&.&.&.&.&.&.&.&.&.&.&.\\
7&.&.&.&.&.&.&.&.&.&.&.&.&.&.&.&.&.&.&.&.\\
13&.&.&.&.&.&.&.&.&.&.&.&.&.&.&.&.&.&.&.&1\\
\hline
&12_{10}&12_{11}&12_{12}&12_{13}&12_{14}&12_{15}&12_{16}&12_{17}
&12_{18}&12_{19}&12_{20}&12_{21}&12_{22}&12_{23}&12_{24}&12_{25}&12_{26}
&12_{27}&12_{28}&13a\\
\hline
2P&6_{30}&6_{24}&6_{24}&6_{17}&6_{29}&6_{37}&6_{36}&6_{38}&6_{28}
&6_{22}&6_{30}&6_{30}&6_{45}&6_{46}&6_{24}&6_{36}&6_{37}&6_{38}&6_{30}
&13b\\
3P&4a&4d&4e&4d&4a&4a&4b&4b&4e&4c&4c&4f&4b&4b&4g&4d&4f&4d&4h&13a\\
\end{array}
$$
} }

\newpage
Character table of $mE$ (continued){
\renewcommand{\baselinestretch}{0.5}

{\tiny \setlength{\arraycolsep}{0,3mm}
$$
\begin{array}{r|rrrrrrrrrrrrrrrrrrrrrrrrrr}
2&1&2&2&2&3&1&.&3&3&3&2&2&2&2&3&3&3&1&2&2&3&1&1&3&3\\
3&1&1&.&.&2&2&2&4&4&4&4&4&4&4&3&3&3&4&3&3&2&3&3&1&.\\
5&.&.&.&.&1&1&1&.&.&.&.&.&.&.&.&.&.&.&.&.&.&.&.&1&1\\
7&.&1&1&1&.&.&.&.&.&.&.&.&.&.&.&.&.&.&.&.&.&.&.&.&.\\
13&1&.&.&.&.&.&.&.&.&.&.&.&.&.&.&.&.&.&.&.&.&.&.&.&.\\
\hline
&13b&14a&14b&14c&15a&15b&15c&18a&18b&18c&18d&18e&18f&18g&18h&18i
&18j&18k&18l&18m&18n&18o&18p&20a&20b\\
\hline
2P&13a&7a&7a&7a&15a&15b&15c&9a&9a&9a&9c&9b&9b&9c&9a&9a&9a&9d&9b&9c
&9a&9e&9f&10a&10a\\
3P&13b&14a&14b&14c&5a&5a&5a&6_{6}&6_{2}&6_{6}&6_{6}&6_{6}&6_{2}
&6_{6}&6_{8}&6_{14}&6_{8}&6_{6}&6_{8}&6_{14}&6_{23}&6_{10}&6_{11}&20a
&20b\\
\end{array}
$$
} }

\newpage
Character table of $mE$ (continued){
\renewcommand{\baselinestretch}{0.5}

{\tiny \setlength{\arraycolsep}{0,3mm}
$$
\begin{array}{r|rrrrrrrrrrrrr}
2&1&4&4&1&1&3&2&1&.&.&1&2\\
3&1&1&1&.&.&1&1&1&1&1&1&1\\
5&.&.&.&.&.&1&1&1&.&.&.&1\\
7&1&.&.&.&.&.&.&.&.&.&1&.\\
13&.&.&.&1&1&.&.&.&1&1&.&.\\
\hline
&21a&24a&24b&26a&26b&30a&30b&30c&39a&39b&42a&60a\\
\hline
2P&21a&12_{3}&12_{2}&13a&13b&15a&15a&15b&39b&39a&21a&30a\\
3P&7a&8a&8b&26a&26b&10a&10c&10b&13a&13b&14a&20a\\
\end{array}
$$
} }

\newpage
Character table of $mE$ (continued){
\renewcommand{\baselinestretch}{0.5}

{\tiny \setlength{\arraycolsep}{0,3mm}
$$
\begin{array}{r|rrrrrrrrrrrrrrrrrr}
2&10&10&10&10&10&10&10&10&9&6&5&5&7&4&6&4&4\\
3&10&9&7&6&5&4&4&3&10&10&10&8&8&8&8&8&8\\
5&1&1&1&1&1&1&.&.&1&.&.&1&.&1&.&.&.\\
7&1&1&1&1&.&.&.&.&1&.&.&.&.&.&.&.&.\\
13&1&1&.&.&.&.&.&.&1&.&.&.&.&.&.&.&.\\
\hline
&1a&2a&2b&2c&2d&2e&2f&2g&3a&3b&3c&3d&3e&3f&3g&3h&3i\\
\hline
2P&1a&1a&1a&1a&1a&1a&1a&1a&3a&3b&3c&3d&3e&3f&3g&3h&3i\\
3P&1a&2a&2b&2c&2d&2e&2f&2g&1a&1a&1a&1a&1a&1a&1a&1a&1a\\
\end{array}
$$
} }

\newpage
Character table of $mE$ (continued){
\renewcommand{\baselinestretch}{0.5}

{\tiny \setlength{\arraycolsep}{0,3mm}
$$
\begin{array}{r|rrrrrrrrrrrrrrrrrrrrrr}
2&3&3&3&2&2&1&7&8&7&8&8&7&8&7&4&9&6&9&5&7&6\\
3&8&8&7&7&7&7&3&3&2&2&2&2&1&1&2&7&9&5&7&7&7\\
5&.&.&.&.&.&.&1&.&1&.&.&.&.&.&1&1&.&1&1&.&.\\
7&.&.&.&.&.&.&.&.&.&.&.&.&.&.&.&1&.&.&.&.&.\\
13&.&.&.&.&.&.&.&.&.&.&.&.&.&.&.&.&.&.&.&.&.\\
\hline
&3j&3k&3l&3m&3n&3o&4a&4b&4c&4d&4e&4f&4g&4h&5a&6_{1}&6_{2}&6_{3}
&6_{4}&6_{5}&6_{6}\\
\hline
2P&3j&3k&3l&3m&3n&3o&2d&2f&2d&2f&2f&2d&2f&2d&5a&3a&3b&3a&3d&3e&3b\\
3P&1a&1a&1a&1a&1a&1a&4a&4b&4c&4d&4e&4f&4g&4h&5a&2b&2a&2d&2a&2a&2b\\
\end{array}
$$
} }

\newpage
Character table of $mE$ (continued){
\renewcommand{\baselinestretch}{0.5}

{\tiny \setlength{\arraycolsep}{0,3mm}
$$
\begin{array}{r|rrrrrrrrrrrrrrrrrrrrrrr}
2&5&6&9&4&4&5&5&6&4&4&7&5&5&5&3&6&6&6&6&7&2&5\\
3&7&6&4&7&7&6&6&5&6&6&4&5&5&5&6&4&4&4&4&3&6&4\\
5&.&.&.&.&.&.&.&.&.&.&.&.&.&.&.&.&.&.&.&.&.&.\\
7&.&.&.&.&.&.&.&.&.&.&.&.&.&.&.&.&.&.&.&.&.&.\\
13&.&.&.&.&.&.&.&.&.&.&.&.&.&.&.&.&.&.&.&.&.&.\\
\hline
&6_{7}&6_{8}&6_{9}&6_{10}&6_{11}&6_{12}&6_{13}&6_{14}&6_{15}
&6_{16}&6_{17}&6_{18}&6_{19}&6_{20}&6_{21}&6_{22}&6_{23}&6_{24}&6_{25}
&6_{26}&6_{27}&6_{28}\\
\hline
2P&3c&3b&3a&3h&3i&3d&3e&3b&3g&3f&3e&3e&3c&3d&3l&3e&3b&3b&3g&3e&3m
&3c\\
3P&2b&2c&2f&2a&2a&2b&2b&2d&2b&2b&2f&2c&2d&2c&2a&2d&2e&2f&2f&2g&2a
&2f\\
\end{array}
$$
} }

\newpage
Character table of $mE$ (continued){
\renewcommand{\baselinestretch}{0.5}

{\tiny \setlength{\arraycolsep}{0,3mm}
$$
\begin{array}{r|rrrrrrrrrrrrrrrrrrrrrrrrrr}
2&5&5&5&5&3&6&6&4&4&4&4&4&2&5&5&5&3&3&3&3&3&4&4&2&2\\
3&4&4&4&4&5&3&3&4&4&4&4&4&5&3&3&3&4&4&4&4&4&3&3&4&4\\
5&.&.&.&.&.&.&.&.&.&.&.&.&.&.&.&.&.&.&.&.&.&.&.&.&.\\
7&.&.&.&.&.&.&.&.&.&.&.&.&.&.&.&.&.&.&.&.&.&.&.&.&.\\
13&.&.&.&.&.&.&.&.&.&.&.&.&.&.&.&.&.&.&.&.&.&.&.&.&.\\
\hline
&6_{29}&6_{30}&6_{31}&6_{32}&6_{33}&6_{34}&6_{35}&6_{36}&6_{37}
&6_{38}&6_{39}&6_{40}&6_{41}&6_{42}&6_{43}&6_{44}&6_{45}&6_{46}&6_{47}
&6_{48}&6_{49}&6_{50}&6_{51}&6_{52}&6_{53}\\
\hline
2P&3g&3d&3d&3e&3l&3b&3e&3i&3f&3h&3g&3f&3n&3e&3d&3d&3j&3k&3l&3l&3l
&3h&3i&3m&3n\\
3P&2d&2d&2f&2f&2b&2g&2e&2f&2d&2f&2f&2f&2b&2g&2g&2e&2f&2f&2f&2f&2c
&2g&2g&2f&2f\\
\end{array}
$$
} }

\newpage
Character table of $mE$ (continued){
\renewcommand{\baselinestretch}{0.5}

{\tiny \setlength{\arraycolsep}{0,3mm}
$$
\begin{array}{r|rrrrrrrrrrrrrrrrrrrrrrrrrrr}
2&2&3&3&1&2&2&5&5&5&5&3&2&2&1&1&1&.&.&4&4&3&4&3&6&7&7\\
3&4&3&3&4&3&1&1&1&.&.&5&5&5&5&4&4&4&4&1&1&1&.&.&3&3&2\\
5&.&.&.&.&.&.&.&.&.&.&.&.&.&.&.&.&.&.&1&1&1&1&1&1&.&.\\
7&.&.&.&.&.&1&.&.&.&.&.&.&.&.&.&.&.&.&.&.&.&.&.&.&.&.\\
13&.&.&.&.&.&.&.&.&.&.&.&.&.&.&.&.&.&.&.&.&.&.&.&.&.&.\\
\hline
&6_{54}&6_{55}&6_{56}&6_{57}&6_{58}&7a&8a&8b&8c&8d&9a&9b&9c&9d
&9e&9f&9g&9h&10a&10b&10c&10d&10e&12_{1}&12_{2}&12_{3}\\
\hline
2P&3n&3l&3l&3o&3m&7a&4e&4b&4b&4e&9a&9b&9c&9d&9e&9f&9g&9h&5a&5a&5a
&5a&5a&6_{3}&6_{9}&6_{9}\\
3P&2f&2g&2g&2f&2g&7a&8a&8b&8c&8d&3b&3b&3b&3b&3h&3i&3i&3h&10a&10b
&10c&10d&10e&4a&4b&4e\\
\end{array}
$$
} }

\newpage
Character table of $mE$ (continued){
\renewcommand{\baselinestretch}{0.5}

{\tiny \setlength{\arraycolsep}{0,3mm}
$$
\begin{array}{r|rrrrrrrrrrrrrrrrrrrrrr}
2&5&5&6&4&4&4&4&5&5&5&3&3&3&3&4&4&4&4&2&2&5\\
3&3&3&2&3&3&3&3&2&2&2&3&3&3&3&2&2&2&2&3&3&1\\
5&.&.&.&.&.&.&.&.&.&.&.&.&.&.&.&.&.&.&.&.&.\\
7&.&.&.&.&.&.&.&.&.&.&.&.&.&.&.&.&.&.&.&.&.\\
13&.&.&.&.&.&.&.&.&.&.&.&.&.&.&.&.&.&.&.&.&.\\
\hline &12_{4}&12_{5}&12_{6}&12_{7}&12_{8}&12_{9}&12_{10}&12_{11}
&12_{12}&12_{13}&12_{14}&12_{15}&12_{16}&12_{17}&12_{18}&12_{19}&12_{20}
&12_{21}&12_{22}&12_{23}&12_{24}\\
\hline
2P&6_{24}&6_{17}&6_{3}&6_{28}&6_{25}&6_{22}&6_{30}&6_{24}&6_{24}
&6_{17}&6_{29}&6_{37}&6_{36}&6_{38}&6_{28}&6_{22}&6_{30}&6_{30}&6_{45}
&6_{46}&6_{24}\\
3P&4b&4b&4f&4b&4b&4a&4a&4d&4e&4d&4a&4a&4b&4b&4e&4c&4c&4f&4b&4b&4g\\
\end{array}
$$
} }

\newpage
Character table of $mE$ (continued){
\renewcommand{\baselinestretch}{0.5}

{\tiny \setlength{\arraycolsep}{0,3mm}
$$
\begin{array}{r|rrrrrrrrrrrrrrrrrrrrrrrrr}
2&3&3&3&4&1&1&2&2&2&3&1&.&3&3&3&2&2&2&2&3&3&3&1&2\\
3&2&2&2&1&1&1&1&.&.&2&2&2&4&4&4&4&4&4&4&3&3&3&4&3\\
5&.&.&.&.&.&.&.&.&.&1&1&1&.&.&.&.&.&.&.&.&.&.&.&.\\
7&.&.&.&.&.&.&1&1&1&.&.&.&.&.&.&.&.&.&.&.&.&.&.&.\\
13&.&.&.&.&1&1&.&.&.&.&.&.&.&.&.&.&.&.&.&.&.&.&.&.\\
\hline
&12_{25}&12_{26}&12_{27}&12_{28}&13a&13b&14a&14b&14c&15a&15b&15c
&18a&18b&18c&18d&18e&18f&18g&18h&18i&18j&18k&18l\\
\hline
2P&6_{36}&6_{37}&6_{38}&6_{30}&13b&13a&7a&7a&7a&15a&15b&15c&9a&9a
&9a&9c&9b&9b&9c&9a&9a&9a&9d&9b\\
3P&4d&4f&4d&4h&13a&13b&14a&14b&14c&5a&5a&5a&6_{6}&6_{2}&6_{6}
&6_{6}&6_{6}&6_{2}&6_{6}&6_{8}&6_{14}&6_{8}&6_{6}&6_{8}\\
\end{array}
$$
} }

\newpage
Character table of $mE$ (continued){
\renewcommand{\baselinestretch}{0.5}

{\tiny \setlength{\arraycolsep}{0,3mm}
$$
\begin{array}{r|rrrrrrrrrrrrrrrrrrr}
2&2&3&1&1&3&3&1&4&4&1&1&3&2&1&.&.&1&2\\
3&3&2&3&3&1&.&1&1&1&.&.&1&1&1&1&1&1&1\\
5&.&.&.&.&1&1&.&.&.&.&.&1&1&1&.&.&.&1\\
7&.&.&.&.&.&.&1&.&.&.&.&.&.&.&.&.&1&.\\
13&.&.&.&.&.&.&.&.&.&1&1&.&.&.&1&1&.&.\\
\hline
&18m&18n&18o&18p&20a&20b&21a&24a&24b&26a&26b&30a&30b&30c&39a&39b
&42a&60a\\
\hline
2P&9c&9a&9e&9f&10a&10a&21a&12_{3}&12_{2}&13a&13b&15a&15a&15b&39b
&39a&21a&30a\\
3P&6_{14}&6_{23}&6_{10}&6_{11}&20a&20b&7a&8a&8b&26a&26b&10a&10c
&10b&13a&13b&14a&20a\\
\end{array}
$$
} } where $ A = -6\zeta(3) - 3$, $B = -12\zeta(3) - 6$, $C =
-\zeta(13)^{11} - \zeta(13)^8 - \zeta(13)^7 - \zeta(13)^6 -
\zeta(13)^5 - \zeta(13)^2 - 1$, $D = -C - 1$, $E = 2D$, $F = 2C$.
\end{ctab}

\normalsize
\newpage

\begin{ctab}\label{Fi_24 ct H}Character table of $H(\Fi'_{24}) = \langle b, p, r \rangle$
{ \setlength{\arraycolsep}{1mm}
\renewcommand{\baselinestretch}{0.5}

{\tiny \setlength{\arraycolsep}{0,3mm}
$$
\begin{array}{r|rrrrrrrrrrrrrrrrrrrr}
2&21&21&19&20&14&14&17&17&16&13&8&12&10&5&8&5&6&15&15\\
3&7&7&4&2&4&4&3&3&1&2&7&6&7&7&5&6&4&5&2\\
5&1&1&1&1&1&1&.&.&.&.&1&.&.&.&.&.&.&1&.\\
7&1&1&.&.&.&.&.&.&.&.&1&.&.&.&.&.&.&.&.\\
\hline
&1a&2a&2b&2c&2d&2e&2f&2g&2h&2i&3a&3b&3c&3d&3e&3f&3g&4a&4b\\
\hline
2P&1a&1a&1a&1a&1a&1a&1a&1a&1a&1a&3a&3b&3c&3d&3e&3f&3g&2a&2a\\
3P&1a&2a&2b&2c&2d&2e&2f&2g&2h&2i&1a&1a&1a&1a&1a&1a&1a&4a&4b\\
5P&1a&2a&2b&2c&2d&2e&2f&2g&2h&2i&3a&3b&3c&3d&3e&3f&3g&4a&4b\\
\end{array}
$$
} }

\newpage
Character table of $H(\Fi'_{24})$ (continued){
\renewcommand{\baselinestretch}{0.5}

{\tiny \setlength{\arraycolsep}{0,3mm}
$$
\begin{array}{r|rrrrrrrrrrrrrrrrrrrrrrrr}
2&13&11&13&14&11&11&11&14&11&12&10&10&10&9&9&10&10&9&6&8&12&10&5\\
3&1&2&2&1&2&2&2&.&1&.&1&1&1&1&1&.&.&.&1&7&6&7&7\\
5&1&1&.&.&.&.&.&.&.&.&.&.&.&.&.&.&.&.&1&1&.&.&.\\
7&.&.&.&.&.&.&.&.&.&.&.&.&.&.&.&.&.&.&.&1&.&.&.\\
\hline
&4c&4d&4e&4f&4g&4h&4i&4j&4k&4l&4m&4n&4o&4p&4q&4r&4s&4t&5a&6_{1}
&6_{2}&6_{3}&6_{4}\\
\hline
2P&2c&2b&2c&2c&2f&2b&2f&2c&2b&2c&2h&2f&2h&2g&2g&2h&2h&2h&5a&3a&3b
&3c&3d\\
3P&4c&4d&4e&4f&4g&4h&4i&4j&4k&4l&4m&4n&4o&4p&4q&4r&4s&4t&5a&2a&2a
&2a&2a\\
5P&4c&4d&4e&4f&4g&4h&4i&4j&4k&4l&4m&4n&4o&4p&4q&4r&4s&4t&1a&6_{1}
&6_{2}&6_{3}&6_{4}\\
\end{array}
$$
} }

\newpage
Character table of $H(\Fi'_{24})$ (continued){
\renewcommand{\baselinestretch}{0.5}

{\tiny \setlength{\arraycolsep}{0,3mm}
$$
\begin{array}{r|rrrrrrrrrrrrrrrrrrrrrrrr}
2&8&9&10&5&11&8&8&8&8&8&8&8&8&8&8&6&7&5&5&5&5&7&5\\
3&5&4&3&6&2&3&3&3&3&3&3&3&3&3&3&4&3&4&4&4&4&2&3\\
5&.&.&.&.&.&.&.&.&.&.&.&.&.&.&.&.&.&.&.&.&.&.&.\\
7&.&.&.&.&.&.&.&.&.&.&.&.&.&.&.&.&.&.&.&.&.&.&.\\
\hline &6_{5}&6_{6}&6_{7}&6_{8}&6_{9}&6_{10}&6_{11}&6_{12}&6_{13}
&6_{14}&6_{15}&6_{16}&6_{17}&6_{18}&6_{19}&6_{20}&6_{21}&6_{22}&6_{23}
&6_{24}&6_{25}&6_{26}&6_{27}\\
\hline
2P&3e&3c&3b&3f&3b&3b&3a&3b&3b&3a&3c&3c&3e&3e&3b&3g&3e&3d&3d&3d&3d
&3b&3e\\
3P&2a&2b&2b&2a&2c&2f&2g&2e&2d&2f&2f&2g&2g&2f&2g&2a&2b&2d&2d&2e&2e
&2i&2d\\
5P&6_{5}&6_{6}&6_{7}&6_{8}&6_{9}&6_{10}&6_{11}&6_{12}&6_{13}
&6_{14}&6_{15}&6_{16}&6_{17}&6_{18}&6_{19}&6_{20}&6_{21}&6_{23}&6_{22}
&6_{25}&6_{24}&6_{26}&6_{27}\\
\end{array}
$$
} }

\newpage
Character table of $H(\Fi'_{24})$ (continued){
\renewcommand{\baselinestretch}{0.5}

{\tiny \setlength{\arraycolsep}{0,3mm}
$$
\begin{array}{r|rrrrrrrrrrrrrrrrrrrrrrrrrrrrr}
2&5&5&5&5&5&5&5&7&5&3&1&8&9&8&8&8&8&8&8&6&6&6&6&5&3&1&2&2\\
3&3&3&3&3&3&3&3&1&2&2&1&2&1&1&.&.&.&.&.&1&.&.&.&4&4&4&3&3\\
5&.&.&.&.&.&.&.&.&.&.&.&.&.&.&.&.&.&.&.&.&.&.&.&.&.&.&.&.\\
7&.&.&.&.&.&.&.&.&.&.&1&.&.&.&.&.&.&.&.&.&.&.&.&.&.&.&.&.\\
\hline
&6_{28}&6_{29}&6_{30}&6_{31}&6_{32}&6_{33}&6_{34}&6_{35}&6_{36}
&6_{37}&7a&8a&8b&8c&8d&8e&8f&8g&8h&8i&8j&8k&8l&9a&9b&9c&9d&9e\\
\hline
2P&3d&3e&3e&3d&3e&3f&3f&3e&3g&3g&7a&4a&4b&4b&4j&4j&4f&4j&4f&4g&4r
&4r&4i&9a&9b&9c&9e&9d\\
3P&2g&2g&2f&2f&2e&2g&2f&2h&2c&2i&7a&8a&8b&8c&8d&8e&8f&8g&8h&8i&8j
&8k&8l&3c&3c&3c&3d&3d\\
5P&6_{28}&6_{29}&6_{30}&6_{31}&6_{32}&6_{33}&6_{34}&6_{35}&6_{36}
&6_{37}&7a&8a&8b&8c&8d&8e&8f&8g&8h&8i&8j&8k&8l&9a&9b&9c&9e&9d\\
\end{array}
$$
} }

\newpage
Character table of $H(\Fi'_{24})$ (continued){
\renewcommand{\baselinestretch}{0.5}

{\tiny \setlength{\arraycolsep}{0,3mm}
$$
\begin{array}{r|rrrrrrrrrrrrrrrrrrrrrrr}
2&6&5&5&5&4&4&9&7&9&6&4&8&8&7&5&5&6&6&6&6&6&7\\
3&1&.&.&.&.&.&5&5&3&4&5&2&2&2&3&3&2&2&2&2&2&1\\
5&1&1&1&1&1&1&.&.&.&.&.&.&.&.&.&.&.&.&.&.&.&.\\
7&.&.&.&.&.&.&.&.&.&.&.&.&.&.&.&.&.&.&.&.&.&.\\
\hline &10a&10b&10c&10d&10e&10f&12_{1}&12_{2}&12_{3}&12_{4}&12_{5}
&12_{6}&12_{7}&12_{8}&12_{9}&12_{10}&12_{11}&12_{12}&12_{13}&12_{14}
&12_{15}&12_{16}\\
\hline
2P&5a&5a&5a&5a&5a&5a&6_{2}&6_{3}&6_{2}&6_{5}&6_{8}&6_{9}&6_{9}
&6_{3}&6_{20}&6_{20}&6_{7}&6_{14}&6_{7}&6_{5}&6_{14}&6_{9}\\
3P&10a&10c&10b&10d&10e&10f&4a&4a&4a&4a&4a&4e&4e&4b&4a&4a&4d&4g&4h
&4b&4i&4c\\
5P&2a&2c&2c&2b&2d&2e&12_{1}&12_{2}&12_{3}&12_{4}&12_{5}&12_{6}
&12_{7}&12_{8}&12_{9}&12_{10}&12_{11}&12_{12}&12_{13}&12_{14}&12_{15}
&12_{16}\\
\end{array}
$$
} }

\newpage
Character table of $H(\Fi'_{24})$ (continued){
\renewcommand{\baselinestretch}{0.5}

{\tiny \setlength{\arraycolsep}{0,3mm}
$$
\begin{array}{r|rrrrrrrrrrrrrrrrrrrr}
2&7&5&5&5&6&4&4&4&4&4&4&4&5&5&4&4&4&4&1\\
3&1&2&2&2&1&2&2&2&2&2&2&2&1&1&1&1&1&1&1\\
5&.&.&.&.&.&.&.&.&.&.&.&.&.&.&.&.&.&.&.\\
7&.&.&.&.&.&.&.&.&.&.&.&.&.&.&.&.&.&.&1\\
\hline
&12_{17}&12_{18}&12_{19}&12_{20}&12_{21}&12_{22}&12_{23}&12_{24}
&12_{25}&12_{26}&12_{27}&12_{28}&12_{29}&12_{30}&12_{31}&12_{32}&12_{33}
&12_{34}&14a\\
\hline
2P&6_{9}&6_{15}&6_{6}&6_{15}&6_{7}&6_{21}&6_{21}&6_{36}&6_{31}
&6_{31}&6_{36}&6_{8}&6_{35}&6_{35}&6_{17}&6_{31}&6_{11}&6_{31}&7a\\
3P&4f&4g&4h&4i&4k&4d&4h&4e&4g&4i&4e&4b&4o&4m&4q&4n&4p&4n&14a\\
5P&12_{17}&12_{18}&12_{19}&12_{20}&12_{21}&12_{22}&12_{23}&12_{24}
&12_{25}&12_{26}&12_{27}&12_{28}&12_{29}&12_{30}&12_{31}&12_{34}&12_{33}
&12_{32}&14a\\
\end{array}
$$
} }

\newpage
Character table of $H(\Fi'_{24})$ (continued){
\renewcommand{\baselinestretch}{0.5}

{\tiny \setlength{\arraycolsep}{0,3mm}
$$
\begin{array}{r|rrrrrrrrrrrrrrrrrrrrrrrrr}
2&1&5&5&3&1&4&2&2&2&2&2&2&4&3&3&3&1&1&5&5&3&3&4&4\\
3&1&.&4&4&4&2&3&3&2&2&2&2&.&.&.&.&1&1&2&2&2&2&1&1\\
5&1&.&.&.&.&.&.&.&.&.&.&.&1&1&1&1&.&.&.&.&.&.&.&.\\
7&.&.&.&.&.&.&.&.&.&.&.&.&.&.&.&.&1&1&.&.&.&.&.&.\\
\hline
&15a&16a&18a&18b&18c&18d&18e&18f&18g&18h&18i&18j&20a&20b&20c&20d
&21a&21b&24a&24b&24c&24d&24e&24f\\
\hline
2P&15a&8b&9a&9b&9c&9a&9e&9d&9e&9d&9d&9e&10a&10b&10d&10c&21b&21a
&12_{3}&12_{1}&12_{10}&12_{9}&12_{12}&12_{14}\\
3P&5a&16a&6_{3}&6_{3}&6_{3}&6_{6}&6_{4}&6_{4}&6_{22}&6_{24}&6_{23}
&6_{25}&20a&20d&20c&20b&7a&7a&8a&8a&8a&8a&8i&8c\\
5P&3a&16a&18a&18b&18c&18d&18f&18e&18i&18j&18g&18h&4a&4c&4d&4c&21a
&21b&24a&24b&24c&24d&24i&24g\\
\end{array}
$$
} }

\newpage
Character table of $H(\Fi'_{24})$ (continued){
\renewcommand{\baselinestretch}{0.5}

{\tiny \setlength{\arraycolsep}{0,3mm}
$$
\begin{array}{r|rrrrrrrrrr}
2&4&4&4&1&4&4&2&1&1\\
3&1&1&1&1&3&3&3&1&1\\
5&.&.&.&1&.&.&.&.&.\\
7&.&.&.&.&.&.&.&1&1\\
\hline
&24g&24h&24i&30a&36a&36b&36c&42a&42b\\
\hline
2P&12_{14}&12_{8}&12_{12}&15a&18a&18a&18b&21a&21b\\
3P&8c&8b&8i&10a&12_{2}&12_{2}&12_{2}&14a&14a\\
5P&24f&24h&24e&6_{1}&36a&36b&36c&42a&42b\\
\end{array}
$$
} }

\newpage
Character table of $H(\Fi'_{24})$ (continued){
\renewcommand{\baselinestretch}{0.5}

{\tiny \setlength{\arraycolsep}{0,3mm}
$$
\begin{array}{r|rrrrrrrrrrrrrrrrrr}
2&21&21&19&20&14&14&17&17&16&13&8&12&10&5&8&5&6\\
3&7&7&4&2&4&4&3&3&1&2&7&6&7&7&5&6&4\\
5&1&1&1&1&1&1&.&.&.&.&1&.&.&.&.&.&.\\
7&1&1&.&.&.&.&.&.&.&.&1&.&.&.&.&.&.\\
\hline
&1a&2a&2b&2c&2d&2e&2f&2g&2h&2i&3a&3b&3c&3d&3e&3f&3g\\
\hline
2P&1a&1a&1a&1a&1a&1a&1a&1a&1a&1a&3a&3b&3c&3d&3e&3f&3g\\
3P&1a&2a&2b&2c&2d&2e&2f&2g&2h&2i&1a&1a&1a&1a&1a&1a&1a\\
5P&1a&2a&2b&2c&2d&2e&2f&2g&2h&2i&3a&3b&3c&3d&3e&3f&3g\\
\end{array}
$$
} }

\newpage
Character table of $H(\Fi'_{24})$ (continued){
\renewcommand{\baselinestretch}{0.5}

{\tiny \setlength{\arraycolsep}{0,3mm}
$$
\begin{array}{r|rrrrrrrrrrrrrrrrrrrrrr}
2&15&15&13&11&13&14&11&11&11&14&11&12&10&10&10&9&9&10&10&9&6\\
3&5&2&1&2&2&1&2&2&2&.&1&.&1&1&1&1&1&.&.&.&1\\
5&1&.&1&1&.&.&.&.&.&.&.&.&.&.&.&.&.&.&.&.&1\\
7&.&.&.&.&.&.&.&.&.&.&.&.&.&.&.&.&.&.&.&.&.\\
\hline
&4a&4b&4c&4d&4e&4f&4g&4h&4i&4j&4k&4l&4m&4n&4o&4p&4q&4r&4s&4t&5a\\
\hline
2P&2a&2a&2c&2b&2c&2c&2f&2b&2f&2c&2b&2c&2h&2f&2h&2g&2g&2h&2h&2h&5a\\
3P&4a&4b&4c&4d&4e&4f&4g&4h&4i&4j&4k&4l&4m&4n&4o&4p&4q&4r&4s&4t&5a\\
5P&4a&4b&4c&4d&4e&4f&4g&4h&4i&4j&4k&4l&4m&4n&4o&4p&4q&4r&4s&4t&1a\\
\end{array}
$$
} }

\newpage
Character table of $H(\Fi'_{24})$ (continued){
\renewcommand{\baselinestretch}{0.5}

{\tiny \setlength{\arraycolsep}{0,3mm}
$$
\begin{array}{r|rrrrrrrrrrrrrrrrrrrrrr}
2&8&12&10&5&8&9&10&5&11&8&8&8&8&8&8&8&8&8&8&6&7\\
3&7&6&7&7&5&4&3&6&2&3&3&3&3&3&3&3&3&3&3&4&3\\
5&1&.&.&.&.&.&.&.&.&.&.&.&.&.&.&.&.&.&.&.&.\\
7&1&.&.&.&.&.&.&.&.&.&.&.&.&.&.&.&.&.&.&.&.\\
\hline
&6_{1}&6_{2}&6_{3}&6_{4}&6_{5}&6_{6}&6_{7}&6_{8}&6_{9}&6_{10}
&6_{11}&6_{12}&6_{13}&6_{14}&6_{15}&6_{16}&6_{17}&6_{18}&6_{19}&6_{20}
&6_{21}\\
\hline
2P&3a&3b&3c&3d&3e&3c&3b&3f&3b&3b&3a&3b&3b&3a&3c&3c&3e&3e&3b&3g&3e\\
3P&2a&2a&2a&2a&2a&2b&2b&2a&2c&2f&2g&2e&2d&2f&2f&2g&2g&2f&2g&2a&2b\\
5P&6_{1}&6_{2}&6_{3}&6_{4}&6_{5}&6_{6}&6_{7}&6_{8}&6_{9}&6_{10}
&6_{11}&6_{12}&6_{13}&6_{14}&6_{15}&6_{16}&6_{17}&6_{18}&6_{19}&6_{20}
&6_{21}\\
\end{array}
$$
} }

\newpage
Character table of $H(\Fi'_{24})$ (continued){
\renewcommand{\baselinestretch}{0.5}

{\tiny \setlength{\arraycolsep}{0,3mm}
$$
\begin{array}{r|rrrrrrrrrrrrrrrrrrrrrrrrrrr}
2&5&5&5&5&7&5&5&5&5&5&5&5&5&7&5&3&1&8&9&8&8&8&8&8&8&6\\
3&4&4&4&4&2&3&3&3&3&3&3&3&3&1&2&2&1&2&1&1&.&.&.&.&.&1\\
5&.&.&.&.&.&.&.&.&.&.&.&.&.&.&.&.&.&.&.&.&.&.&.&.&.&.\\
7&.&.&.&.&.&.&.&.&.&.&.&.&.&.&.&.&1&.&.&.&.&.&.&.&.&.\\
\hline
&6_{22}&6_{23}&6_{24}&6_{25}&6_{26}&6_{27}&6_{28}&6_{29}&6_{30}
&6_{31}&6_{32}&6_{33}&6_{34}&6_{35}&6_{36}&6_{37}&7a&8a&8b&8c&8d&8e&8f
&8g&8h&8i\\
\hline
2P&3d&3d&3d&3d&3b&3e&3d&3e&3e&3d&3e&3f&3f&3e&3g&3g&7a&4a&4b&4b&4j
&4j&4f&4j&4f&4g\\
3P&2d&2d&2e&2e&2i&2d&2g&2g&2f&2f&2e&2g&2f&2h&2c&2i&7a&8a&8b&8c&8d
&8e&8f&8g&8h&8i\\
5P&6_{23}&6_{22}&6_{25}&6_{24}&6_{26}&6_{27}&6_{28}&6_{29}&6_{30}
&6_{31}&6_{32}&6_{33}&6_{34}&6_{35}&6_{36}&6_{37}&7a&8a&8b&8c&8d&8e&8f
&8g&8h&8i\\
\end{array}
$$
} }

\newpage
Character table of $H(\Fi'_{24})$ (continued){
\renewcommand{\baselinestretch}{0.5}

{\tiny \setlength{\arraycolsep}{0,3mm}
$$
\begin{array}{r|rrrrrrrrrrrrrrrrrrrrrrrrr}
2&6&6&6&5&3&1&2&2&6&5&5&5&4&4&9&7&9&6&4&8&8&7&5&5\\
3&.&.&.&4&4&4&3&3&1&.&.&.&.&.&5&5&3&4&5&2&2&2&3&3\\
5&.&.&.&.&.&.&.&.&1&1&1&1&1&1&.&.&.&.&.&.&.&.&.&.\\
7&.&.&.&.&.&.&.&.&.&.&.&.&.&.&.&.&.&.&.&.&.&.&.&.\\
\hline
&8j&8k&8l&9a&9b&9c&9d&9e&10a&10b&10c&10d&10e&10f&12_{1}&12_{2}
&12_{3}&12_{4}&12_{5}&12_{6}&12_{7}&12_{8}&12_{9}&12_{10}\\
\hline
2P&4r&4r&4i&9a&9b&9c&9e&9d&5a&5a&5a&5a&5a&5a&6_{2}&6_{3}&6_{2}
&6_{5}&6_{8}&6_{9}&6_{9}&6_{3}&6_{20}&6_{20}\\
3P&8j&8k&8l&3c&3c&3c&3d&3d&10a&10c&10b&10d&10e&10f&4a&4a&4a&4a&4a
&4e&4e&4b&4a&4a\\
5P&8j&8k&8l&9a&9b&9c&9e&9d&2a&2c&2c&2b&2d&2e&12_{1}&12_{2}&12_{3}
&12_{4}&12_{5}&12_{6}&12_{7}&12_{8}&12_{9}&12_{10}\\
\end{array}
$$
} }

\newpage
Character table of $H(\Fi'_{24})$ (continued){
\renewcommand{\baselinestretch}{0.5}

{\tiny \setlength{\arraycolsep}{0,3mm}
$$
\begin{array}{r|rrrrrrrrrrrrrrrrrrrrr}
2&6&6&6&6&6&7&7&5&5&5&6&4&4&4&4&4&4&4&5&5\\
3&2&2&2&2&2&1&1&2&2&2&1&2&2&2&2&2&2&2&1&1\\
5&.&.&.&.&.&.&.&.&.&.&.&.&.&.&.&.&.&.&.&.\\
7&.&.&.&.&.&.&.&.&.&.&.&.&.&.&.&.&.&.&.&.\\
\hline
&12_{11}&12_{12}&12_{13}&12_{14}&12_{15}&12_{16}&12_{17}&12_{18}
&12_{19}&12_{20}&12_{21}&12_{22}&12_{23}&12_{24}&12_{25}&12_{26}&12_{27}
&12_{28}&12_{29}&12_{30}\\
\hline
2P&6_{7}&6_{14}&6_{7}&6_{5}&6_{14}&6_{9}&6_{9}&6_{15}&6_{6}&6_{15}
&6_{7}&6_{21}&6_{21}&6_{36}&6_{31}&6_{31}&6_{36}&6_{8}&6_{35}&6_{35}\\
3P&4d&4g&4h&4b&4i&4c&4f&4g&4h&4i&4k&4d&4h&4e&4g&4i&4e&4b&4o&4m\\
5P&12_{11}&12_{12}&12_{13}&12_{14}&12_{15}&12_{16}&12_{17}&12_{18}
&12_{19}&12_{20}&12_{21}&12_{22}&12_{23}&12_{24}&12_{25}&12_{26}&12_{27}
&12_{28}&12_{29}&12_{30}\\
\end{array}
$$
} }

\newpage
Character table of $H(\Fi'_{24})$ (continued){
\renewcommand{\baselinestretch}{0.5}

{\tiny \setlength{\arraycolsep}{0,3mm}
$$
\begin{array}{r|rrrrrrrrrrrrrrrrrrrrrrrrr}
2&4&4&4&4&1&1&5&5&3&1&4&2&2&2&2&2&2&4&3&3&3&1&1&5\\
3&1&1&1&1&1&1&.&4&4&4&2&3&3&2&2&2&2&.&.&.&.&1&1&2\\
5&.&.&.&.&.&1&.&.&.&.&.&.&.&.&.&.&.&1&1&1&1&.&.&.\\
7&.&.&.&.&1&.&.&.&.&.&.&.&.&.&.&.&.&.&.&.&.&1&1&.\\
\hline
&12_{31}&12_{32}&12_{33}&12_{34}&14a&15a&16a&18a&18b&18c&18d&18e
&18f&18g&18h&18i&18j&20a&20b&20c&20d&21a&21b&24a\\
\hline
2P&6_{17}&6_{31}&6_{11}&6_{31}&7a&15a&8b&9a&9b&9c&9a&9e&9d&9e&9d
&9d&9e&10a&10b&10d&10c&21b&21a&12_{3}\\
3P&4q&4n&4p&4n&14a&5a&16a&6_{3}&6_{3}&6_{3}&6_{6}&6_{4}&6_{4}
&6_{22}&6_{24}&6_{23}&6_{25}&20a&20d&20c&20b&7a&7a&8a\\
5P&12_{31}&12_{34}&12_{33}&12_{32}&14a&3a&16a&18a&18b&18c&18d&18f
&18e&18i&18j&18g&18h&4a&4c&4d&4c&21a&21b&24a\\
\end{array}
$$
} }

\newpage
Character table of $H(\Fi'_{24})$ (continued){
\renewcommand{\baselinestretch}{0.5}

{\tiny \setlength{\arraycolsep}{0,3mm}
$$
\begin{array}{r|rrrrrrrrrrrrrrr}
2&5&3&3&4&4&4&4&4&1&4&4&2&1&1\\
3&2&2&2&1&1&1&1&1&1&3&3&3&1&1\\
5&.&.&.&.&.&.&.&.&1&.&.&.&.&.\\
7&.&.&.&.&.&.&.&.&.&.&.&.&1&1\\
\hline
&24b&24c&24d&24e&24f&24g&24h&24i&30a&36a&36b&36c&42a&42b\\
\hline
2P&12_{1}&12_{10}&12_{9}&12_{12}&12_{14}&12_{14}&12_{8}&12_{12}
&15a&18a&18a&18b&21a&21b\\
3P&8a&8a&8a&8i&8c&8c&8b&8i&10a&12_{2}&12_{2}&12_{2}&14a&14a\\
5P&24b&24c&24d&24i&24g&24f&24h&24e&6_{1}&36a&36b&36c&42a&42b\\
\end{array}
$$
} }, where $ A = 12\zeta(3) + 2$, $B = -3\zeta(3) - 2$, $C =
-\zeta(3)$, $D = -12\zeta(3) - 6$, $E = -2\zeta(3) - 1$, $F =
2\zeta(21)_3\zeta(21)_7^4 + 2\zeta(21)_3\zeta(21)_7^2 +
2\zeta(21)_3\zeta(21)_7 + \zeta(21)_3 + \zeta(21)_7^4 +
\zeta(21)_7^2 + \zeta(21)_7 + 1$, $G = -F + 1$, $H = -3E$, $I =
6\zeta(3) + 1$, $J = -4\zeta(5)^3 - 4\zeta(5)^2 - 1$, $K =
4\zeta(5)^3 + 4\zeta(5)^2 + 3$, $L = 2\zeta(5)^3 + 2\zeta(5)^2 +
1$, $M = 18\zeta(3) + 9$, $N = -4\zeta(12)_4\zeta(12)_3 -
2\zeta(12)_4 $.
\end{ctab}

\newpage

\begin{ctab}\label{Fi_24 ct A_1}Character table of $A_1(\Fi'_{24}) = \langle y,q,s \rangle$
{ \setlength{\arraycolsep}{1mm}
\renewcommand{\baselinestretch}{0.5}

{\tiny \setlength{\arraycolsep}{0,3mm}
$$
\begin{array}{r|rrrrrrrrrrrrrrrrr}
2&19&19&17&14&19&19&14&17&10&9&8&5&11&14&14&13\\
3&9&9&6&6&4&4&4&3&7&9&7&7&4&4&1&3\\
5&2&2&1&2&1&1&1&.&1&.&.&.&1&.&1&.\\
7&1&1&1&1&.&.&.&.&1&.&.&.&1&.&.&.\\
11&1&1&1&.&.&.&.&.&.&.&.&.&.&.&.&.\\
13&1&1&.&.&.&.&.&.&.&.&.&.&.&.&.&.\\
\hline
&1a&2a&2b&2c&2d&2e&2f&2g&3a&3b&3c&3d&4a&4b&4c&4d\\
\hline
2P&1a&1a&1a&1a&1a&1a&1a&1a&3a&3b&3c&3d&2a&2e&2e&2e\\
\end{array}
$$
} }

\newpage
Character table of $A_1(\Fi'_{24})$ (continued){
\renewcommand{\baselinestretch}{0.5}

{\tiny \setlength{\arraycolsep}{0,3mm}
$$
\begin{array}{r|rrrrrrrrrrrrrrrrrrrrrr}
2&11&13&11&11&11&12&11&10&5&10&9&8&8&8&8&8&8&5&9&9&10\\
3&2&1&2&2&2&1&1&1&1&7&9&7&5&5&6&4&4&7&4&4&3\\
5&1&.&.&.&.&.&.&.&2&1&.&.&1&1&.&.&.&.&.&.&.\\
7&.&.&.&.&.&.&.&.&.&1&.&.&.&.&.&1&1&.&.&.&.\\
11&.&.&.&.&.&.&.&.&.&.&.&.&.&.&.&.&.&.&.&.&.\\
13&.&.&.&.&.&.&.&.&.&.&.&.&.&.&.&.&.&.&.&.&.\\
\hline
&4e&4f&4g&4h&4i&4j&4k&4l&5a&6_{1}&6_{2}&6_{3}&6_{4}&6_{5}&6_{6}
&6_{7}&6_{8}&6_{9}&6_{10}&6_{11}&6_{12}\\
\hline
2P&2d&2e&2d&2g&2g&2e&2d&2g&5a&3a&3b&3c&3a&3a&3b&3a&3a&3d&3b&3b&3a\\
\end{array}
$$
} }

\newpage
Character table of $A_1(\Fi'_{24})$ (continued){
\renewcommand{\baselinestretch}{0.5}

{\tiny \setlength{\arraycolsep}{0,3mm}
$$
\begin{array}{r|rrrrrrrrrrrrrrrrrrrrrrrrrr}
2&10&5&5&5&8&8&8&8&8&7&7&5&5&5&5&5&5&5&5&3&8&8&8&8&8\\
3&3&6&5&5&3&3&3&3&3&3&3&4&4&4&4&3&3&3&3&1&1&1&1&1&1\\
5&.&.&.&.&.&.&.&.&.&.&.&.&.&.&.&.&.&.&.&.&.&.&.&.&.\\
7&.&.&.&.&.&.&.&.&.&.&.&.&.&.&.&.&.&.&.&1&.&.&.&.&.\\
11&.&.&.&.&.&.&.&.&.&.&.&.&.&.&.&.&.&.&.&.&.&.&.&.&.\\
13&.&.&.&.&.&.&.&.&.&.&.&.&.&.&.&.&.&.&.&.&.&.&.&.&.\\
\hline
&6_{13}&6_{14}&6_{15}&6_{16}&6_{17}&6_{18}&6_{19}&6_{20}&6_{21}
&6_{22}&6_{23}&6_{24}&6_{25}&6_{26}&6_{27}&6_{28}&6_{29}&6_{30}&6_{31}
&7a&8a&8b&8c&8d&8e\\
\hline
2P&3a&3b&3c&3c&3c&3c&3b&3a&3a&3c&3c&3d&3d&3d&3d&3d&3d&3c&3c&7a&4f
&4d&4d&4d&4d\\
\end{array}
$$
} }

\newpage
Character table of $A_1(\Fi'_{24})$ (continued){
\renewcommand{\baselinestretch}{0.5}

{\tiny \setlength{\arraycolsep}{0,3mm}
$$
\begin{array}{r|rrrrrrrrrrrrrrrrrrrrrrrrr}
2&8&6&6&4&3&2&5&4&5&5&3&4&2&6&8&5&5&8&8&6&7&7&5&8\\
3&.&1&.&4&4&3&1&1&.&.&1&.&.&3&3&4&4&2&2&3&2&2&3&1\\
5&.&.&.&.&.&.&2&2&1&1&1&1&.&1&.&.&.&.&.&.&.&.&.&.\\
7&.&.&.&.&.&.&.&.&.&.&.&.&.&.&.&.&.&.&.&.&.&.&.&.\\
11&.&.&.&.&.&.&.&.&.&.&.&.&1&.&.&.&.&.&.&.&.&.&.&.\\
13&.&.&.&.&.&.&.&.&.&.&.&.&.&.&.&.&.&.&.&.&.&.&.&.\\
\hline
&8f&8g&8h&9a&9b&9c&10a&10b&10c&10d&10e&10f&11a&12_{1}&12_{2}
&12_{3}&12_{4}&12_{5}&12_{6}&12_{7}&12_{8}&12_{9}&12_{10}&12_{11}\\
\hline
2P&4f&4h&4i&9a&9b&9c&5a&5a&5a&5a&5a&5a&11a&6_{1}&6_{13}&6_{10}
&6_{2}&6_{13}&6_{13}&6_{10}&6_{13}&6_{13}&6_{22}&6_{13}\\
\end{array}
$$
} }

\newpage
Character table of $A_1(\Fi'_{24})$ (continued){
\renewcommand{\baselinestretch}{0.5}

{\tiny \setlength{\arraycolsep}{0,3mm}
$$
\begin{array}{r|rrrrrrrrrrrrrrrrrrrrr}
2&6&6&6&6&4&5&5&5&6&6&4&4&4&4&5&4&4&1&3&2\\
3&2&2&2&2&3&2&2&2&1&1&2&2&2&2&1&1&1&.&1&1\\
5&.&.&.&.&.&.&.&.&.&.&.&.&.&.&.&.&.&.&.&.\\
7&.&.&.&.&.&.&.&.&.&.&.&.&.&.&.&.&.&.&1&1\\
11&.&.&.&.&.&.&.&.&.&.&.&.&.&.&.&.&.&.&.&.\\
13&.&.&.&.&.&.&.&.&.&.&.&.&.&.&.&.&.&1&.&.\\
\hline
&12_{12}&12_{13}&12_{14}&12_{15}&12_{16}&12_{17}&12_{18}&12_{19}
&12_{20}&12_{21}&12_{22}&12_{23}&12_{24}&12_{25}&12_{26}&12_{27}&12_{28}
&13a&14a&14b\\
\hline
2P&6_{12}&6_{17}&6_{12}&6_{17}&6_{3}&6_{19}&6_{19}&6_{11}&6_{10}
&6_{12}&6_{23}&6_{23}&6_{28}&6_{28}&6_{22}&6_{28}&6_{28}&13a&7a&7a\\
\end{array}
$$
} }

\newpage
Character table of $A_1(\Fi'_{24})$ (continued){
\renewcommand{\baselinestretch}{0.5}

{\tiny \setlength{\arraycolsep}{0,3mm}
$$
\begin{array}{r|rrrrrrrrrrrrrrrrrrrrrrrrr}
2&2&3&5&5&4&3&3&3&3&3&4&4&2&2&2&2&2&3&3&3&2&2&2&2\\
3&.&1&.&.&4&4&3&3&3&3&2&2&3&3&3&2&2&1&.&.&1&.&.&.\\
5&.&1&.&.&.&.&.&.&.&.&.&.&.&.&.&.&.&1&1&1&.&.&.&.\\
7&1&.&.&.&.&.&.&.&.&.&.&.&.&.&.&.&.&.&.&.&1&.&.&.\\
11&.&.&.&.&.&.&.&.&.&.&.&.&.&.&.&.&.&.&.&.&.&1&1&1\\
13&.&.&.&.&.&.&.&.&.&.&.&.&.&.&.&.&.&.&.&.&.&.&.&.\\
\hline
&14c&15a&16a&16b&18a&18b&18c&18d&18e&18f&18g&18h&18i&18j&18k&18l
&18m&20a&20b&20c&21a&22a&22b&22c\\
\hline
2P&7a&15a&8f&8f&9a&9b&9a&9a&9a&9a&9a&9a&9b&9c&9b&9c&9c&10a&10c&10d
&21a&11a&11a&11a\\
\end{array}
$$
} }

\newpage
Character table of $A_1(\Fi'_{24})$ (continued){
\renewcommand{\baselinestretch}{0.5}

{\tiny \setlength{\arraycolsep}{0,3mm}
$$
\begin{array}{r|rrrrrrrrrrrrrrrrrrrrr}
2&5&5&5&5&4&4&4&4&1&2&3&2&2&3&3&2&2&2&2&2\\
3&1&1&1&1&1&1&1&1&.&.&1&1&1&2&2&2&1&1&1&1\\
5&.&.&.&.&.&.&.&.&.&.&1&1&1&.&.&.&.&.&.&1\\
7&.&.&.&.&.&.&.&.&.&1&.&.&.&.&.&.&1&1&1&.\\
11&.&.&.&.&.&.&.&.&.&.&.&.&.&.&.&.&.&.&.&.\\
13&.&.&.&.&.&.&.&.&1&.&.&.&.&.&.&.&.&.&.&.\\
\hline
&24a&24b&24c&24d&24e&24f&24g&24h&26a&28a&30a&30b&30c&36a&36b&36c
&42a&42b&42c&60a\\
\hline
2P&12_{5}&12_{6}&12_{6}&12_{5}&12_{13}&12_{20}&12_{13}&12_{20}&13a
&14a&15a&15a&15a&18g&18g&18b&21a&21a&21a&30a\\
\end{array}
$$
} }

\newpage
Character table of $A_1(\Fi'_{24})$ (continued){
\renewcommand{\baselinestretch}{0.5}

{\tiny \setlength{\arraycolsep}{0,3mm}
$$
\begin{array}{r|rrrrrrrrrrrrrrrr}
2&19&19&17&14&19&19&14&17&10&9&8&5&11&14&14\\
3&9&9&6&6&4&4&4&3&7&9&7&7&4&4&1\\
5&2&2&1&2&1&1&1&.&1&.&.&.&1&.&1\\
7&1&1&1&1&.&.&.&.&1&.&.&.&1&.&.\\
11&1&1&1&.&.&.&.&.&.&.&.&.&.&.&.\\
13&1&1&.&.&.&.&.&.&.&.&.&.&.&.&.\\
\hline
&1a&2a&2b&2c&2d&2e&2f&2g&3a&3b&3c&3d&4a&4b&4c\\
\hline
2P&1a&1a&1a&1a&1a&1a&1a&1a&3a&3b&3c&3d&2a&2e&2e\\
\end{array}
$$
} }

\newpage
Character table of $A_1(\Fi'_{24})$ (continued){
\renewcommand{\baselinestretch}{0.5}

{\tiny \setlength{\arraycolsep}{0,3mm}
$$
\begin{array}{r|rrrrrrrrrrrrrrrrrrrrr}
2&13&11&13&11&11&11&12&11&10&5&10&9&8&8&8&8&8&8&5&9\\
3&3&2&1&2&2&2&1&1&1&1&7&9&7&5&5&6&4&4&7&4\\
5&.&1&.&.&.&.&.&.&.&2&1&.&.&1&1&.&.&.&.&.\\
7&.&.&.&.&.&.&.&.&.&.&1&.&.&.&.&.&1&1&.&.\\
11&.&.&.&.&.&.&.&.&.&.&.&.&.&.&.&.&.&.&.&.\\
13&.&.&.&.&.&.&.&.&.&.&.&.&.&.&.&.&.&.&.&.\\
\hline
&4d&4e&4f&4g&4h&4i&4j&4k&4l&5a&6_{1}&6_{2}&6_{3}&6_{4}&6_{5}
&6_{6}&6_{7}&6_{8}&6_{9}&6_{10}\\
\hline
2P&2e&2d&2e&2d&2g&2g&2e&2d&2g&5a&3a&3b&3c&3a&3a&3b&3a&3a&3d&3b\\
\end{array}
$$
} }

\newpage
Character table of $A_1(\Fi'_{24})$ (continued){
\renewcommand{\baselinestretch}{0.5}

{\tiny \setlength{\arraycolsep}{0,3mm}
$$
\begin{array}{r|rrrrrrrrrrrrrrrrrrrrrrrr}
2&9&10&10&5&5&5&8&8&8&8&8&7&7&5&5&5&5&5&5&5&5&3&8\\
3&4&3&3&6&5&5&3&3&3&3&3&3&3&4&4&4&4&3&3&3&3&1&1\\
5&.&.&.&.&.&.&.&.&.&.&.&.&.&.&.&.&.&.&.&.&.&.&.\\
7&.&.&.&.&.&.&.&.&.&.&.&.&.&.&.&.&.&.&.&.&.&1&.\\
11&.&.&.&.&.&.&.&.&.&.&.&.&.&.&.&.&.&.&.&.&.&.&.\\
13&.&.&.&.&.&.&.&.&.&.&.&.&.&.&.&.&.&.&.&.&.&.&.\\
\hline
&6_{11}&6_{12}&6_{13}&6_{14}&6_{15}&6_{16}&6_{17}&6_{18}&6_{19}
&6_{20}&6_{21}&6_{22}&6_{23}&6_{24}&6_{25}&6_{26}&6_{27}&6_{28}&6_{29}
&6_{30}&6_{31}&7a&8a\\
\hline
2P&3b&3a&3a&3b&3c&3c&3c&3c&3b&3a&3a&3c&3c&3d&3d&3d&3d&3d&3d&3c&3c
&7a&4f\\
\end{array}
$$
} }

\newpage
Character table of $A_1(\Fi'_{24})$ (continued){
\renewcommand{\baselinestretch}{0.5}

{\tiny \setlength{\arraycolsep}{0,3mm}
$$
\begin{array}{r|rrrrrrrrrrrrrrrrrrrrrrrrrr}
2&8&8&8&8&8&6&6&4&3&2&5&4&5&5&3&4&2&6&8&5&5&8&8&6&7\\
3&1&1&1&1&.&1&.&4&4&3&1&1&.&.&1&.&.&3&3&4&4&2&2&3&2\\
5&.&.&.&.&.&.&.&.&.&.&2&2&1&1&1&1&.&1&.&.&.&.&.&.&.\\
7&.&.&.&.&.&.&.&.&.&.&.&.&.&.&.&.&.&.&.&.&.&.&.&.&.\\
11&.&.&.&.&.&.&.&.&.&.&.&.&.&.&.&.&1&.&.&.&.&.&.&.&.\\
13&.&.&.&.&.&.&.&.&.&.&.&.&.&.&.&.&.&.&.&.&.&.&.&.&.\\
\hline &8b&8c&8d&8e&8f&8g&8h&9a&9b&9c&10a&10b&10c&10d&10e&10f&11a
&12_{1}&12_{2}&12_{3}&12_{4}&12_{5}&12_{6}&12_{7}&12_{8}\\
\hline
2P&4d&4d&4d&4d&4f&4h&4i&9a&9b&9c&5a&5a&5a&5a&5a&5a&11a&6_{1}
&6_{13}&6_{10}&6_{2}&6_{13}&6_{13}&6_{10}&6_{13}\\
\end{array}
$$
} }

\newpage
Character table of $A_1(\Fi'_{24})$ (continued){
\renewcommand{\baselinestretch}{0.5}

{\tiny \setlength{\arraycolsep}{0,3mm}
$$
\begin{array}{r|rrrrrrrrrrrrrrrrrrrrr}
2&7&5&8&6&6&6&6&4&5&5&5&6&6&4&4&4&4&5&4&4\\
3&2&3&1&2&2&2&2&3&2&2&2&1&1&2&2&2&2&1&1&1\\
5&.&.&.&.&.&.&.&.&.&.&.&.&.&.&.&.&.&.&.&.\\
7&.&.&.&.&.&.&.&.&.&.&.&.&.&.&.&.&.&.&.&.\\
11&.&.&.&.&.&.&.&.&.&.&.&.&.&.&.&.&.&.&.&.\\
13&.&.&.&.&.&.&.&.&.&.&.&.&.&.&.&.&.&.&.&.\\
\hline
&12_{9}&12_{10}&12_{11}&12_{12}&12_{13}&12_{14}&12_{15}&12_{16}
&12_{17}&12_{18}&12_{19}&12_{20}&12_{21}&12_{22}&12_{23}&12_{24}&12_{25}
&12_{26}&12_{27}&12_{28}\\
\hline
2P&6_{13}&6_{22}&6_{13}&6_{12}&6_{17}&6_{12}&6_{17}&6_{3}&6_{19}
&6_{19}&6_{11}&6_{10}&6_{12}&6_{23}&6_{23}&6_{28}&6_{28}&6_{22}&6_{28}
&6_{28}\\
\end{array}%
$$
} }
where $ A = 16\zeta(3) + 8$, $B = 6\zeta(3) + 3$, $C =
-2\zeta(3) - 1$, $D = 4\zeta(3) + 2$, $E = -2\zeta(11)^9 -
2\zeta(11)^5 - 2\zeta(11)^4 - 2\zeta(11)^3 - 2\zeta(11) - 1$, $F =
18\zeta(3) + 9$, $G = 112\zeta(3) + 56$, $H = 12\zeta(3) + 6$, $I
= 4\zeta(12)_4\zeta(12)_3 + 2\zeta(12)_4$, $J = 96\zeta(3) + 48$,
$K = -128\zeta(3) - 64 $.
\end{ctab}
\normalsize

\setlength\textheight{39\baselineskip}
\addtolength\textheight{\topskip}

\end{appendix}



\newpage

\begin{thebibliography}{10}

\bibitem{magma}
John Cannon and Catherine Playoust.
\newblock {\em An Introduction to {\sc Magma}}.
\newblock School of Mathematics and Statistics, University of Sydney, 1993.

\bibitem{cannon}
John~J. Cannon, Derek~F. Holt.
\newblock Automorphism group computation and isomorphism
testing in finite groups.
\newblock {\em J. Symbolic Computat.}, 35:241-267, 2003.

\bibitem{carter}
R.W.\ Carter.
\newblock {\em Simple groups of Lie type}.
\newblock John Wiley and Sons, London, 1972.


\bibitem{atlas}
J.H.\ Conway, R.T.\ Curtis, S.P.\ Norton, R.A.\ Parker, and R.A.\
Wilson.
\newblock {\em Atlas of finite groups}.
\newblock Clarendon Press, Oxford, 1985.

\bibitem{fischer}
B. Fischer.
\newblock {Finite groups generated by $3$-transpositions}.
\newblock{\em Inventiones Math.}, 13:232-246, 1971.

\bibitem{fischer1}
B. Fischer.
\newblock {Finite groups generated by $3$-transpositions II}.
\newblock Lecture Notes, University of Warwick, Coventry, 1979.

\bibitem{hall}
J. Hall, L. H. Soicher
\newblock {Presenations of some $3$-transposition groups}.
\newblock{\em Comm. Algebra}, 23:232-246, 1995.

\bibitem{holt5} Holt, D. F., Cohomology and group extensions in
Magma, in W. Bosma, J. Cannon (eds), Discovering mathematics with
Magma, Springer, Berlin 2006, pp. 221--141.

\bibitem{james}
G.\ James.
\newblock The modular characters of the Mathieu groups.
\newblock {\em J. Algebra}, 27:57-111, 1973.

\bibitem{kim}
H. Kim, G. O. Michler.
\newblock Simultaneous constructions of the sporadic groups $\Co_2$ and $\Fi_{22}$.
\newblock in (L. -C. Kappe, A. Magidin, R. F. Morse, eds.)
\textit{Computational Group Theory and the Theory of Groups},
Contemporary Mathematics {\bf 470}, 141--234, Amer. Math. Soc, Providence,
RI., (2008).

\bibitem{kim1}
H. Kim.
\newblock Representation theoretic existence proof for Fischer's sporadic group $\Fi_{23}$.
\newblock Senior Thesis, Dept. Math. Cornell University, 2008. See also arXiv: 0904.0639v1.

\bibitem{michler}
G. O. Michler.
\newblock Theory of finite simple groups.
\newblock Cambridge University Press, Cambridge, 2006.

\bibitem{michler1}
G. O. Michler.
\newblock Constructing finite simple groups from irreducible subgroups of
$\GL_n(2)$,
\newblock in (L. -C. Kappe, A. Magidin, R. F. Morse, eds.)
\textit{Computational Group Theory and the Theory of Groups},
Contemporary Mathematics {\bf 470}, 235--262, Amer. Math. Soc, Providence,
RI., (2008).

\bibitem{michler2}
G. O. Michler.
\newblock Theory of finite simple groups II.
\newblock Cambridge University Press, Cambridge (to appear in 2009).


\bibitem{praeger} Ch. Praeger, L. Soicher.
\newblock Low rank representations and graphs for sporadic groups.
\newblock Cambridge University Press, Cambridge, 1997
\end{thebibliography}
\end{document}